%
%
%
%
\input amstex.tex
\documentstyle{amsppt}
%
%
%
%
\def\section#1{\par\bigpagebreak \csname subhead\endcsname #1 \endsubhead\par\medpagebreak}
\def\subsection#1{\csname subsubhead\endcsname #1 \endsubsubhead}
\def\Theorem#1#2{\csname proclaim\endcsname{Theorem #1} #2 \endproclaim}
\def\Corollary#1{\csname proclaim\endcsname{Corollary} #1 \endproclaim}
\def\Proposition#1#2{\csname proclaim\endcsname{Proposition #1} #2 \endproclaim}
\def\Lemma#1#2{\csname proclaim\endcsname{Lemma #1} #2 \endproclaim}
\def\Remark#1#2{\remark{Remark {\rm #1}} #2 \endremark}
\def\Proof#1{\demo{Proof} #1\qed\enddemo}
\def\Proofof#1#2{\demo{Proof of #1} #2\qed\enddemo}
%
%
\def\ld{\lambda}
\def\Ld{\Lambda}
\def\ep{\epsilon}
\def\vep{\varepsilon}
\def\Dt{\Delta}

\def\Mat{\text{\rm Mat}}
\def\GL{\text{\rm GL}}
\def\SL{\text{\rm SL}}
\def\SO{\text{\rm SO}}
\def\SP{\text{\rm Sp}}
\def\SU{\text{\rm SU}}
\def\gl{\frak{gl}}
\def\so{\frak{so}}
\def\sp{\frak{sp}}
\def\g{\frak{g}}
\def\k{\frak{k}}
\def\n{\frak{n}}
\def\kq{\frak{k}_q}
\def\N{\Bbb{N}}
\def\Z{\Bbb{Z}}
\def\Q{\Bbb{Q}}

\def\C{\Bbb{C}}
\def\K{\Bbb{K}}
\def\T{\Bbb{T}}
\def\Hom{\operatorname{Hom}_{\K}}
\def\End{\operatorname{End}_{\K}}

\def\otimesK{\otimes_{\K}}
\def\Uqgl{U_q(\gl(N))}
\def\AqGL{A_q(\GL(N))}
\def\Uq{U_q(\g)}
\def\Aq{A_q(G)}
\def\detq{{\operatorname{det}}_q}
\def\id{\operatorname{id}}
\def\pr{\text{\rm pr}}
\def\bv{\text{\bf v}}
\def\ox{\otimes}
\def\oxK{\otimes_{\K}}
\def\diag#1{\text{\rm diag}(#1)}
%
\leftheadtext{Masatoshi NOUMI}
\rightheadtext\nofrills{Macdonald polynomials on quantum homogeneous spaces}
\par\noindent
{\eightpoint {\it To appear in}: Adv.\ in Math.}
\par\vskip1cm
\topmatter
\title\nofrills
        Macdonald's symmetric polynomials \\
        as zonal spherical functions \\
        on some quantum homogeneous spaces
\endtitle
\author
        Masatoshi NOUMI
\endauthor
\affil
        Department of Mathematical Sciences\\
        University of Tokyo
\endaffil
\address
        Komaba 3-8-1, Meguro-Ku, Tokyo 153, Japan
\endaddress
\dedicatory
	Dedicated to Professor I.\,M.\,Gelfand for his eightieth birthday
\enddedicatory
\email
        noumi\@tansei.cc.u-tokyo.ac.jp
\endemail
\comment
\abstract
Quantum analogues of the homogeneous spaces 
$\GL(n)/\SO(n)$ and $\GL(2n)/\Sp(2n)$ are introduced.
The zonal spherical functions on these quantum homogeneous spaces 
are represented by Macdonald's symmetric polynomials 
$P_{\ld}=P_{\ld}(x_1,\cdots,x_n;q,t)$ 
with $t=q^{1 \over 2}$ or $t=q^2$. 
\endabstract
\endcomment
\endtopmatter
\document 
%
%
%
%
%
\par\bigpagebreak
%
%
\section{Introduction}
%
In this paper we introduce some quantum analogues of the homogeneous spaces 
$ \GL(n)/\SO(n) $ and $\GL(2n)/\SP(2n) $
in the framework of quantum general linear groups. 
On these ``quantum homogeneous spaces'', we investigate the zonal spherical 
functions associated with finite dimensional representations.
As a result, we will see that the zonal spherical functions in question are
represented by Macdonald's symmetric polynomials  
$ P_\ld=P_\ld(x_1,\cdots,x_n;q,t)$
(of type $A_{n-1}$) with  $t=q^{1\over 2}$ or $t=q^2$ (\cite{M2}). 
\par
This result can be regarded as a generalization of Koornwinder's
realization of the continuous $q$-Legendre polynomials by the quantum group
$\SU_q(2)$ (\cite{K1}). 
Our quantum analogue of $\GL(n)/\SO(n)$ is essentially the same as 
the one discussed by Ueno-Takebayashi \cite{UT}. 
As to the quantum analogue of $\GL(3)/\SO(3)$, 
it is already known by \cite{UT} that Macdonald's
symmetric polynomials arise as zonal spherical functions. 
Our result contains the affirmative answer to their conjecture for the case 
where $n>3$. 
Main results of this paper are announced in \cite{N3}. 
\medpagebreak\par
Throughout this paper, we will denote by $G$ the general linear group $\GL(N)$ and by $\g$ its Lie algebra $\gl(N)$. 
The $q$-deformation of the coordinate ring $A(G)$ of $G=\GL(N)$ and the universal enveloping algebra $U(\g)$ of $\g=\gl(N)$ will be denoted by $A_q(G)$ and by $U_q(\g)$, respectively. 
We consider the following closed subgroup $K$ of $G$ for ``quantization'':
$$
\align
\text{Case (SO):}& \ \  K=\SO(n)=\{g\in\GL(n)\,; \ g g^t =\id_n, \det(g)=1 \} \ \  (N=n)\\
\text{Case (Sp):}& \ \  K=\SP(2n)=\{g\in\GL(2n)\,; \ gJ_n g^t=J_n \} \ \  (N=2n),
\endalign 
$$
where $J_n=\sum_{k=1}^{n}{e_{2k-1,2k}-e_{2k,2k-1}}$. 
The corresponding Lie subalgebra of $\g$ will be denoted by $\k$.  
After the preliminaries on the quantum general linear groups ($A_q(G)$ and $U_q(\g)$), 
we introduce in Section 2 some {\it coideals} $\kq$ of $\Uq$, 
corresponding to the Lie subalgebras $\k=\so(n)\subset\gl(n)$ and 
$\k=\sp(2n)\subset\gl(2n)$. 
The construction of $\kq$ is carried out in the framework of $L$-operators 
as in \cite{RTF}, by using some constant solutions $J$ to the so-called 
{\it reflection equation}. 
In Case (SO), the coideals $\kq$ are closely related to the {\it twisted} 
$U_q(\frak{so}(n))$ of Gavrilik-Klimyk \cite{GK} (see Section 2.4). 
With these coideals, we investigate the quantum analogue of the homogeneous 
spaces $\GL(n)/\SO(n)$ and $\GL(2n)/\SP(2n)$ and their zonal spherical 
functions associated with finite dimensional representations. 
Actually we consider the invariant ring 
$$
        A_q(G/K) := \{\varphi \in \Aq;\ \kq.\varphi=0 \}
$$
as a $q$-deformation of the algebra of regular functions on the homogeneous space 
$G/K$. 
In order to consider the zonal spherical functions on this quantum homogeneous space 
$(G/K)_q$, we also investigate the subalgebra of $\kq$-biinvariants 
$$
 {\Cal H}=A_q(K\backslash G /K)
    :=\{\varphi\in \Aq;\, \kq.\varphi=\varphi.\kq=0 \}. 
$$
\par 
Representations with $\kq$-fixed vectors and the structure of these 
subalgebras of $\Aq$ will be studied in Sections 3 and 4, respectively.
In accordance with the classical case where $q=1$, 
the invariant ring $A_q(G/K)$ has a multiplicity free irreducible decomposition as 
a $\Uq$-module. 
The algebra ${\Cal H}=A_q(K\backslash G /K)$ is then decomposed into the direct sum of 
simultaneous eigenspaces ${\Cal H}(\ld)$ of the center of $\Uq$, 
which are parametrized by the dominant integral weights 
$\ld$ corresponding to the irreducible highest weight representations $V(\ld)$ with 
$\kq$-fixed vectors. 
Each simultaneous eigenspace ${\Cal H}(\ld)$ turns out to be one-dimensional 
just as in the case where $q=1$. 
We call a nonzero element in ${\Cal H}(\ld)$ with some normalization the 
{\it zonal spherical function \/} $\varphi(\ld)$ associated with $V(\ld)$. 
\par 
To describe the zonal spherical functions $\varphi(\ld)$ thus obtained, 
we will study in Section 5 their restriction $\varphi(\ld)|_\T$ 
to the diagonal subgroup $\T$ of $\GL_q(N)$. 
They are determined by analyzing the {\it radial component} of 
the central element $C_1$ of $\Uq$, presented in \cite{RTF}. 
This central element gives rise to a $q$-difference operator acting 
on a subalgebra of the Laurent polynomial ring.
As the eigenfunctions of the $q$-difference operator, $\varphi(\ld)|_\T$
are identified with Macdonald's symmetric functions $P_{\mu}(x_1,\cdots,x_n;q,t)$ (\cite{M2}) 
with $(q,t)$ replaced by $(q^4,q^2)$ in Case (SO), and by $(q^2,q^4)$ in Case (Sp). 
As an application of this realization, we will give in Section 6 
an evaluation of the ratio 
$\langle P_{\mu}, P_{\mu}\rangle'_{q,t}/\langle 1,1\rangle'_{q,t}$ 
of scalar products, for the special values $t=q^{1\over 2}$ and $t=q^2$, 
to see the coincidence with the formula proposed in \cite{M3}. 
This computation is based on the $q$-analogue of Schur's orthogonality of 
matrix elements of unitary representations of the unitary group 
$\text{\rm U}(n)$ and a result of Macdonald \cite{M2} on the principal 
specialization of $P_\mu$. 
\section{Contents} \par
\S1. Preliminaries on the quantum general linear groups. \par
\S2. Quantum analogue of some Lie subalgebras of $\gl(N)$.\par
\S3. Representations with $\kq$-fixed vectors. \par
\S4. Quantum homogeneous spaces and zonal spherical functions. \par
\S5. Macdonald's symmetric polynomials as zonal spherical functions. \par
\S6. Scalar product and orthogonality. \par
References
%
\section{\S1. Preliminaries on the quantum general linear groups}
%
%
In this section we give a review on the quantum general linear group 
$\GL_q(N)$. 
Main references for this section are: Jimbo \cite{J} and 
Reshetikhin-Takhtajan-Faddeev \cite{RTF} (see also \cite{NYM} and \cite{NUW1}). 
From now on, we fix the field $\K=\Q(q)$ of rational functions in the 
indeterminate $q$ as the ground field. 
Our presentation of the quantized universal enveloping algebra 
$\Uqgl$ is slightly different from those in \cite{J} and \cite{RTF}; 
modifications are made in order to avoid the fractional powers of $q$. 
In the arguments of this paper the ground field $\K$ can replaced 
by the field $\C$ of complex numbers, assuming that 
$q$ is a real number with $|q|\ne 0,1$. 
%
%
\subsection{1.1. Quantized coordinate ring $A_q(\GL(N))$}
%
%
        Let  $V$  be the $N$-dimensional vector space over  $\K$ with canonical
basis $(v_j)_{1\leq j\leq N}$.  We make use of the following quantum 
$R$-matrices $R^\pm$ in  $\End(V\otimesK V)$:
$$\align
        R^+ & =  {\sum}_{1\leq i,j\leq N} q^{\delta_{ij}} e_{ii}\otimes e_{jj}
                +(q-q^{-1}) {\sum}_{1\leq i<j\leq N} e_{ij}\otimes e_{ji}, 
\tag{1.1}\\
         R^-& =  {\sum}_{1\leq i,j\leq N} q^{-\delta_{ij}} e_{ii}\otimes e_{jj}
                -(q-q^{-1}) {\sum}_{1\leq j<i\leq N} e_{ij}\otimes e_{ji},
\endalign$$
where $e_{ij} \in \End(V) \ (1\leq i,j\leq N) $ are the matrix units with 
respect to the basis $(v_j)_{1\leq j\leq N}$.
It is well-known that these matrices satisfy the Yang-Baxter equation
$$
        R^\ep_{12} R^\ep_{13} R^\ep_{23}
        =R^\ep_{23} R^\ep_{13} R^\ep_{12}\quad(\epsilon =\pm)
\tag{1.2}
$$
in $\End(V_1\otimesK V_2\otimesK V_3)$ with $V_a=V$ for $a=1,2,3$. 
Here the subscripts $a,b$ for $R^\epsilon_{ab}$ indicates the pair of 
components it should act on.
We recall that $R^{+}-R^{-}=(q-q^{-1})P$ with the matrix 
$P={\sum}_{ij}e_{ij}\otimes e_{ji}$ representing the flip 
$v\otimes w \mapsto w\otimes v$. 
Note also that $ (R^+_{12})^{-1} = R^-_{21} $ and 
$(R^{\pm}_{12})^t=R^{\pm}_{21}$  with double signs in the same order. 
\par
        Following \cite{RTF}, we define the {\it coordinate ring} 
$A_q(\Mat(N))$ of the quantum matrix space of rank $N$  
to be the $\K$-algebra generated by the 
{\it canonical coordinates} $t_{ij}$ ($1\leq i,j\leq N$) with the fundamental 
relations 
$$\align
\text{(i)}\quad& t_{ki}t_{kj}=qt_{kj}t_{ki}, t_{ik}t_{jk}=qt_{jk}t_{ik}
\quad (i<j),
\tag{1.3}\\
\text{(ii)}\quad& t_{i\ell}t_{kj}=t_{kj}t_{i\ell},
t_{ij}t_{k\ell}-t_{k\ell}t_{ij}=(q-q^{-1})t_{i\ell}t_{kj}\,
\quad(i<k;\  j<\ell).
\endalign$$
In terms of the matrix $T=(t_{ij})_{1\leq i,j\leq N}$ in 
$A_q(\Mat(N))\otimesK\End(V)$,
the commutation relations above are equivalently written as the Yang-Baxter 
equation
$$
R^+_{12} T_2 T_1 =T_1 T_2 R^+_{12} \quad\text{in}
\quad A_q(\Mat(N))\otimesK\End(V\otimesK V).
\tag{1.4}
$$
This algebra $A_q(\Mat(N))$ has a distinguished central element
$$
\detq(T)={\sum}_{w \in \frak{S}_N} (-q)^{\ell(w)} t_{w(1)1}\cdots t_{w(N)N},
\tag{1.5}
$$
called the {\it quantum determinant}.
Here  $\frak{S}_N$ is the permutation group of the indexing set 
$\{1,2,\cdots,N\}$ and, for each $w \in \frak{S}_N$, $\ell(w)$ stands for the 
number of inversions in $w$. 
The {\it coordinate ring}  $A_q(\GL(N))$ of the quantum general linear group  
$\GL_q(N)$  is then defined by adjoining the inverse of the quantum determinant 
$\detq(T)$  to $A_q(\Mat(N))$:  $A_q(\GL(N))=A_q(\Mat(N))[\detq(T)^{-1}]$.  
This algebra has a structure of Hopf algebra such that 
$$
\Delta(t_{ij})={\sum}_{k=1}^N t_{ik}\otimes t_{kj},\quad\vep(t_{ij})=\delta_{ij}
\tag{1.6}
$$
for $1\leq i,j\leq N$.
In the matrix notation, these formulas will be written as 
$\Delta(T)=T\dot{\otimes}T$ and $\vep(T)=\id_V$.
Note also that $\detq(T)$ is a group-like element, namely, 
$\Delta(\detq(T))=\detq(T)\otimes\detq(T)$ and $\vep(\detq(T))=1$.
The antipode $S$ of $\AqGL$ is the $\K$-algebra anti-automorphism such that
$$
TS(T)=S(T)\,T=\id_V,
\tag{1.7}
$$
where $S(T)=(S(t_{ij}))_{1\leq i,j\leq N}$.
%
%
\subsection{1.2. Quantized universal enveloping algebra $\Uqgl$}
%
%
Let $P$ be the weight lattice for $\GL(N)$; it is the free 
$\Z$-module of rank $N$ with canonical basis 
$(\ep_j)_{1\leq j\leq N}$ and we fix a symmetric bilinear form 
$\langle \,,\,\rangle : P\times P\rightarrow\Z $ such that $\langle\ep_i,\ep_j\rangle=\delta_{ij}$ for 
$1\leq i,j\leq N$.
Through this pairing, we will frequently identify $P$ with its dual 
$P^*={\text{\rm Hom}}_{\Z}(P,\Z)$.
We will also use the notation of {\it simple roots}: 
$\alpha_k=\ep_k-\ep_{k+1}$ for $1\leq k\leq N-1$.
\par
        The {\it quantized universal enveloping algebra} $\Uqgl$ is the $\K$-algebra 
generated by the symbols $q^h (h\in P^*)$ and $e_k, f_k (1\leq k\leq N-1) $ 
with the following fundamental relations:
$$\align
 \text{(i)}\quad& q^0=1, \ q^h.q^{h'}=q^{h+h'} , 
\tag{1.8}\\
 \text{(ii)}\quad& q^h e_i q^{-h} = q^{\langle h,\alpha_i\rangle}e_i, \ 
        q^h f_i q^{-h} = q^{-\langle h,\alpha_i\rangle}f_i , \\
 \text{(iii)}\quad& e_i f_j-f_j e_i = 
        \delta_{ij} \frac{q^{\ep_i-\ep_{i+1}}-q^{-\ep_i+\ep_{i+1}}}{q-q^{-1}}, \\
 \text{(iv)}\quad&  e_i^2 e_j-(q+q^{-1})e_i e_j e_i+e_j e_i^2 = 0 \,(|i-j|=1); 
        \quad e_i e_j=e_j e_i \,(|i-j|>1),\\
 \text{(v)}\quad& f_i^2 f_j-(q+q^{-1})f_i f_j f_i+f_j f_i^2 = 0 \,(|i-j|=1);
        \quad f_i f_j=f_j f_i \,(|i-j|>1), 
\endalign$$
where $h,h'\in P^* $ and $1\leq i,j\leq N-1$.
We take the following Hopf algebra structure of $\Uqgl$:
$$\align
\text{(i)}\quad&\Delta(q^h)=q^h\otimes q^h, \ \vep(q^h)=1, \ S(q^h)=q^{-h}, 
\tag{1.9}\\
\text{(ii)}\quad&\Delta(e_k)=e_k \otimes 1 + q^{\ep_k-\ep_{k+1}}\otimes e_k, \ 
        \vep(e_k)=0, \ S(e_k)=-q^{-\ep_k+\ep_{k+1}}e_k,\\
\text{(iii)}\quad&\Delta(f_k)=f_k \otimes q^{-\ep_k+\ep_{k+1}} + 1\otimes f_k, \ 
        \vep(f_k)=0, \ S(f_k)=-f_k q^{\ep_k-\ep_{k+1}},
\endalign$$
where $h\in P^*$ and $1\leq k\leq N-1$. 
We use the notation $U_q(\frak{t})$ to refer the Hopf subalgebra 
$\K[q^h\,;\,(h\in P^*)]$ of $\Uqgl$ corresponding to the diagonal 
Lie subalgebra $\frak{t}\subset\frak{gl}(N)$. 
\par
        In this paper, we will extensively use the $L$-{\it operators} 
$L^+_{ij},L^-_{ij} \in \Uqgl$ as in Reshetikhin-Takhtajan-Faddeev \cite{RTF}; 
these play the role of root vectors of $\frak{gl}(N)$.
As is stated in \cite{J}, there is a unique family of elements $E_{ij} 
\ (1\leq i,j\leq N, i\neq j)$ in $\Uqgl$ such that
$$\align
\text{(i)}\quad& E_{i,i+1}=e_i, 
E_{ij}=E_{ik}E_{kj} -qE_{kj}E_{ik}\quad(i<k<j),
\tag{1.10}\\
\text{(ii)}\quad& E_{i+1,i}= f_i, 
E_{ij}=E_{ik}E_{kj} -q^{-1}E_{kj}E_{ik}\quad(i>k>j).
\endalign$$
With these elements, define the elements  $L^{\pm}_{ij}\in\Uqgl$ by
$$\align
\text{(i)}\quad& L^+_{ii}=q^{\ep_i}, \ L^+_{ij}=(q-q^{-1})q^{\ep_i}E_{ji} \ (i<j), \ 
L^+_{ij}=0 \, (i>j) 
\tag{1.11}\\
\text{(ii)}\quad& L^-_{ii}=q^{-\ep_i}, \  L^-_{ij}=-(q-q^{-1})E_{ji}q^{-\ep_j} \  (i>j), \ 
L^-_{ij}=0 \, (i<j).
\endalign$$
Then it is known by Jimbo \cite{J} that the the matrices 
$L^{\pm}={\sum}_{i,j} e_{ij}\otimes L^{\pm}_{ij}$ 
in $\End(V)\otimes_{\K}\Uqgl$ satisfy the following Yang-Baxter equations:
$$
        R^+_{12} L^{\ep}_1 L^{\ep}_2
        =L^{\ep}_2 L^{\ep}_1 R^+_{12} \quad (\ep=\pm)\quad\text{and}\quad
        R^+_{12} L^+_1 L^-_2 = L^-_2 L^+_1 R^+_{12}.
\tag{1.12}
$$
We also remark that, for the Hopf algebra structure of $\Uqgl$, 
the matrices $L^{\pm}$ satisfy 
$$
\Delta(L^\ep)=L^\ep\dot{\otimes} L^\ep \quad\text{and}\quad\vep(L^\ep)=\id_V 
\quad (\ep=\pm). 
\tag{1.13}
$$
This fact was fundamental in the framework of Reshetikhin-Takhtajan-Faddeev
\cite{RTF}.
We remark that the Yang-Baxter equations (1.12) assure that there exists a
$\K$-algebra homomorphism $\rho_V : \Uqgl \rightarrow \End(V)$  such that 
$$
        R^{\pm} = {\sum}_{1\le i,j \le N}  e_{ij}\otimes \rho_V(L^{\pm}_{ij}). 
\tag{1.14}
$$
The vector space $V$,  regarded as a left $\Uqgl$-module, is called 
the {\it vector representation} of $\Uqgl$. 
\par
	The square $S^2$ of the antipode of $\Uqgl$ is an automorphism of 
the Hopf algebra $\Uqgl$. 
For any $a\in\Uqgl$, we have
$$
	S^2(a)= q^{-2\rho} a q^{2\rho}\quad\text{for any}\ \ a\in \Uqgl,
\tag{1.15}
$$
where $q^{2\rho}$ is the group-like element of $\Uq$ corresponding to 
the sum of positive roots
$$
		2\rho= {\sum}_{k=1}^N  2(N-k)\ep_k. 
\tag{1.16}
$$
%
%
\subsection{1.3. Pairing between $\AqGL$ and $\Uqgl$}
%
        Let $U$ and $A$ be two Hopf algebras over $\K$.
We say that a $\K$-bilinear form $(\ ,\ ): U\times A \rightarrow \K$ is a
{\it pairing of Hopf algebras} if it satisfies the following three conditions:
$$
\align
\text{(i)}\quad&
(a.b,\varphi)=(a\otimes b,\Delta_A(\varphi))\quad\text{and}\quad(1_U,\varphi)=\vep_A(\varphi), 
\tag{1.17}\\
\text{(ii)}\quad&  (a,\varphi.\psi)=(\Delta_U(a),\varphi\otimes\psi)\quad\text{and}\quad(a,1_A)=\vep_U(a),\\
\text{(iii)}\quad& (S_U(a),\varphi)=(a,S_A(\varphi)),
\endalign
$$
for all $a,b\in U$ and $\varphi,\psi\in A$. 
Through such a pairing, one can define a $U$-bimodule structure on $A$ 
by setting
$$ 
a.\varphi = (\id_A\otimes \hat{a})\circ\Delta_A(\varphi)\text{\quad and \quad} 
\varphi.a = (\hat{a}\otimes\id_A)\circ\Delta_A(\varphi),
\tag{1.18}
$$
for any $a\in U$ and $\varphi\in A$.
In the right-hand side, the symbol $\hat{a}$ stands for the linear functional on 
$A$ induced from $a\in U$ by the pairing. 
The algebra $A$ then becomes an {\it algebra with two-sided $U$-symmetry }, 
in the sense that both the multiplication 
$ A\otimesK A \rightarrow A$ and the unit homomorphism 
$ \K \rightarrow A$ are homomorphisms of $U$-bimodules. 
As to the $U$-bimodule structure of $A$, we have
$$
	a.S_A(\varphi)=S_A(\varphi.S_U(a))\quad \text{and} \quad 
	S_A(\varphi).a=S_A(S_U(a).\varphi),
\tag{1.19}
$$
for all $\varphi\in A$ and $a\in U$.
\par 
        We now take the Hopf algebras $\Uqgl$ and $\AqGL$ for $U$ and $A$ above. 
As to these Hopf algebras, 
it is known that there exists a unique pairing of Hopf algebras 
$(\ , \ ) : \Uqgl\times\AqGL \rightarrow \K$
such that 
$$
        (L^{\pm}_1, T_2) = R^{\pm}_{12} \quad\text{and} 
	\quad (L^{\pm},\detq(T))=q^{\pm 1} \id_V. 
\tag{1.20}
$$
By this pairing, the algebra $\AqGL$  becomes an algebra with two-sided symmetry over
$\Uqgl$.  
In terms of the $L$-operators, the $\Uqgl$-bimodule structure of $\AqGL$ is described as follows:
$$
        L^{\ep}_1. T_2 = T_2 R^{\ep}_{12} \quad\text{and}\quad 
        T_2. L^{\ep}_1 = R^{\ep}_{12} T_2 \quad(\ep=\pm). 
\tag{1.21}
$$
By means of the $\Uqgl$-bimodule structure, the square $S^2$ 
of the antipode of $\AqGL$ is described as 
$S^2(\varphi) = q^{2\rho}.\varphi.q^{-2\rho}$ for any $\varphi\in\AqGL$. 
\par
	We say that a left $\Uqgl$-module is $P$-{\it weighted} if it
has a $K$-basis consisting of weight vectors with weights in $P$.  
Let $P^+$ be the set of all dominant integral weights in $P$:
$$
        P^+=\{ \ld={\sum}_{k=1}^N \ld_k \ep_k \in P \ ;\ \ld_1\ge\ld_2\ge\cdots\ge\ld_N \}. 
\tag{1.22}
$$
For each $\ld\in P^+$, we denote by $V(\ld)$ the unique irreducible finite 
dimensional left $\Uqgl$-module with highest weight $\ld$; 
it is characterized as the unique irreducible left $\Uqgl$-module 
generated by an element $u(\ld)$ such that 
$q^h.u(\ld)=q^{\langle h,\ld\rangle}u(\ld)$ 
for $h\in P^*$ and $e_k.u(\ld)=0$ for $1\le k\le N-1$. 
It is well-known that any finite dimensional $P$-weighted $\Uqgl$-module 
is completely reducible and that any finite dimensional irreducible 
$P$-weighted $\Uqgl$-module is isomorphic to $V(\ld)$ for some 
dominant integral weight $\ld\in P^+$ (see \cite{L, R}). 
Furthermore, 
each irreducible $\Uqgl$-modules $V(\ld)$ is obtained as the ``differential 
representation'' from the underlying right $\AqGL$-comodule structure 
of $V(\ld)$. 
All these $V(\ld)\ (\ld\in P^+)$ are realized as right $\AqGL$-subcomodules of 
$\AqGL$ by means of standard monomials of quantum minor determinants 
(see \cite{TT, NYM}, for instance). 
\par
We denote by $W(\ld)$ the $\K$-vector subspace of $\AqGL$ spanned by the 
matrix elements of the right $\AqGL$-comodule $V(\ld)$. 
Then $W(\ld)$ is an irreducible $\Uqgl$-bimodule isomorphic to the 
tensor product $\Hom(V(\ld),\K)\otimesK V(\ld)$. 
Furthermore, the regular representation $\AqGL$ has the irreducible 
decomposition
$$
        \AqGL={\bigoplus}_{\ld\in P^+}\ W(\ld) 
\tag{1.23}
$$
as a $\Uqgl$-bimodule, 
which corresponds to the Peter-Weyl Theorem for the quantum unitary 
group $\text{\sl U}_q(N)$ (see \cite{H1, NYM, W}, for instance). 
We also remark that the $\K$-subspace $W(\ld)$ of $\AqGL$ is 
characterized as the simultaneous eigenspace of the center of $\Uqgl$:
$$
W(\ld)=\{\varphi\in\AqGL;\; C.\varphi=\chi_{\ld}(C)\varphi
\ \text{for any central element}\,  C\in\Uqgl\},
\tag{1.24}
$$
where $\chi_{\ld}(C)$ stands for the eigenvalue of the central element 
$C$ acting on the irreducible representation $V(\ld)$.
%
%
\subsection{1.4. Involutions on $\AqGL$ and $\Uqgl$}
%
In this subsection, we recall some involutions on $\AqGL$ and $\Uqgl$, related to the 
quantum unitary group $\text{\rm U}_q(N)$. 
\par 
	For the moment, let  be $\K$ an arbitrary field and fix an involutive automorphism 
$c\mapsto \overline{c} $, which we call the conjugation of $\K$. 
A Hopf algebra $A$ over $\K$ is called a {\it Hopf $\ast$-algebra} if it has a pair $(\iota, \tau)$ 
of involutive, conjugate-linear mappings $\iota, \tau: A \to A$ such that 
\roster 
	\item $\iota$ is an algebra anti-automorphism and a coalgebra automorphism.
	\item $\tau$ is an algebra automorphism and a coalgebra anti-automorphism.
	\item The antipode $S$ of $A$ is expressed as $S=\iota\circ\tau$. 
\endroster
We will call the involution $\iota$ the {\it $\ast$-operation} of $A$ and write 
$ \iota(a)=a^*$  for $a\in A$. 
Note that, 
if $A$ has an involutive, conjugate linear mapping $\iota$ satisfying (1) 
and the condition $(\iota\circ S)^2=\id_A$, 
it becomes a Hopf $\ast$-algebra together with the involution $\tau=\iota\circ S$. 
The condition $(\iota\circ S)^2=\id_A$ is sometimes referred to as 
{\it Woronowicz's condition} (see \cite{W}). 
\par
	Returning to the previous setting, we take the field $\Q(q)$ of rational functions in $q$ 
as the ground field $\K$. 
In this case, we set $\overline{c}=c$ for any $c\in\Q(q)$. 
When we take $\K=\C$ instead, we denote by $\overline{c}$ the complex conjugation 
and assume that $q$ is a real number with $|q|\ne 0,1$. 
\par
	Either $\AqGL$ or $\Uqgl$ has a structure of Hopf $\ast$-algebra 
corresponding to the real form $\text{\rm U}(N)$ or $\frak{u}(N)$. 
As to $\AqGL$, we can take 
 the $\ast$-operation such that 
$$
	t_{ij}^*= S(t_{ji}) \ \ (1\le i,j\le N) \ \ \text{and}\ \ \detq(T)^*=\detq(T)^{-1}. 
\tag{1.25}
$$
The corresponding $\tau$ is given by the ``transposition":
$$
	\tau(t_{ij})= t_{ji} \ \ (1\le i,j\le N) \ \ \text{and}\ \ \tau(\detq(T))=\detq(T). 
\tag{1.26}
$$
On the other hand, as to $\Uqgl$, we can take the $\ast$-operation such 
$$
	{L^\pm_{ij}}^*= S(L^{\mp}_{ji}) \ \ (1\le i,j\le N). 
\tag{1.27}
$$
At the Chevalley generators, this involution takes the values
$$
(q^h)^*= q^h\ \ (h\in P^*), \ \ {e_k}^*=q^{-1}f_kt_k, \ \ 
{f_k}^*=qt^{-1}_ke_k \ \ (1\le k\le N-1),
\tag{1.28}
$$
where $t_k=q^{\ep_k-\ep_{k+1}}$ for $1\le k\le N-1$. 
The corresponding $\tau$ is given by 
$$
	\tau(L^\pm_{ij})= L^{\mp}_{ji} \ \ (1\le i,j\le N). 
\tag{1.29}
$$
\par
	These Hopf $\ast$-algebra structures on $\AqGL$ and $\Uqgl$ are compatible with 
the pairing, in the sense that 
$$
	(a^*, \varphi)=\overline{(a,\tau(\varphi))}, \quad (\tau(a), \varphi)=\overline{(a,\varphi^*)}
\tag{1.30}
$$
for all $a\in \Uqgl$ and $\varphi\in\AqGL$. 
From this compatibility, we have the following formulas concerning the $\Uq$-bimodule 
structure of $\AqGL$:
$$
	a. \varphi^* =(\tau(a).\varphi)^*\quad\text{and}\quad
	a.\tau(\varphi)=\tau(\varphi.a^*), 
\tag{1.31}
$$
for all $\varphi\in \AqGL$ and $a\in\Uqgl$.
%
\section{\S2. Quantum analogue of some Lie subalgebras of $\gl(N)$}
%
In this section, we introduce quantum analogues of the Lie subalgebras 
$\k=\so(n)\subset\gl(n)$ and $\k=\sp(2n)\subset\gl(2n)$. 
We will define a family of {\it coideals\/} $\kq$ of 
the quantized universal enveloping algebra $\Uqgl$ 
(with $N=n$ or $2n$) which ``tends'' to $\k$ as $q \rightarrow 1$. 
From this section on, we denote by $G$ the general linear group $\GL(N)$ and 
by $\g$ its Lie algebra $\gl(N)$. 

%
\subsection{2.1. Definition of the coideal $\kq$} 
%
        Keeping the notation in the previous section, we denote by $V$ 
the vector representation of $\Uq=\Uqgl$ with canonical basis 
$(v_j)_{1\leq j\leq N}$.
We consider the following two types of Lie subalgebras  $\k$ of 
$\gl(N)=\End(V)$ for quantization:
$$
\align
\text{Case (SO):}& \quad \k=\so(n)=\{X\in\End(V); X +X^t=0\} \ \  \text{with} \ \  N=n
\tag{2.1}\\
\text{Case (Sp):}& \quad \k=\sp(2n)=\{X\in\End(V); X J_n+J_n X^t =0\} \ \  \text{with} \ \  N=2n, 
\endalign 
$$
where $J_n=\sum_{k=1}^{n}{e_{2k-1,2k}-e_{2k,2k-1}}$.
In order to ``quantize'' this setting, we define the matrix $J(a)\in\End(V)$,
depending on the parameters $a=(a_1,\cdots,a_n)$ in the algebraic torus 
$(\K^*)^n$, by
$$
\align
\text{Case (SO):}& \quad J(a):={\sum}_{k=1}^{n}{e_{kk} a_k}, 
\tag{2.2} \\
\text{Case (Sp):}& \quad J(a):={\sum}_{k=1}^{n}{(e_{2k-1,2k}-qe_{2k,2k-1}) a_k}.
\endalign
$$
Note that $J(a)$ is an invertible matrix and that  $J(a)^{-1}=J(a^{-1})$ 
in Case (SO) and $J(a)^{-1}=-q^{-1}J(a^{-1})$ in Case (Sp).
\par
        Fixing the parameter $a\in(\K^*)^n$, we introduce a coideal $\kq=\kq(a)$ 
of $\Uq$ by using the matrix $J=J(a)$. 
Define a matrix  $M=M(a)$ in 
$\End(V)\otimesK \Uq$ by the formula 
$$
 M:=L^+ -JS(L^-)^t\ J^{-1}.     \tag{2.3}
$$
Writing the matrix $M$ in the form $M=\sum_{i,j} e_{ij}\otimes M_{ij}$, 
we denote by $\kq=\kq(a)$ the vector subspace of $\Uq$ spanned by 
the matrix elements $M_{ij}$ ($1\leq i,j\leq N$):
$$
  \kq:={\sum}_{1\leq i,j\leq N} M_{ij} \subset \Uq.     \tag{2.4}
$$
\par 
        We show first that the vector subspace $\kq$ is actually a coideal of $\Uq$. 
\Proposition{2.1}{
For any $a\in(\K^*)^n$, the $\K$-vector subspace $\kq=\kq(a) $ is a coideal of $\Uq$. 
To be more precise, the matrix $M=M(a)$ satisfies 
$$
\Delta(M) =  L^+ \dot{\otimes} M + M \dot{\otimes} J S(L^-)^t J^{-1}
        \quad \text{and} \quad \vep(M)=0. 
\tag{2.5}
$$
}
\Proof{
By $\Delta(L^+)=L^+\dot{\otimes}L^+$ and 
$\Delta(S(L^-)^t)=S(L^-)^t\dot{\otimes} S(L^-)^t $, 
one has
$$\align
\Dt(M)J&=L^+\dot{\otimes}L^+J + JS(L^-)^t \dot{\otimes}S(L^-)^t \tag{2.6}\\
     &=L^+\dot{\otimes}(L^+J-JS(L^-)^t)
        + (L^+J-J S(L^-)^t)\dot{\otimes}S(L^-)^t\\
     &=L^+\dot{\otimes}MJ + MJ\dot{\otimes} S(L^-)^t.
\endalign$$
The assertion $\vep(M)=0$ is clear, since $\vep(L^+)=\vep(S(L^-)^t)=\id_V$. 
}
%
\subsection{2.2. Reflection equation for the matrix $J$}
%
We now explain the reason why the matrices $J=J(a) \ (a\in(\K^*)^n)$ in (2.2) 
are chosen for the quantization of the Lie subalgebras $\so(n)$ and $\sp(2n)$. 
As for these matrices $J=J(a)$, the following lemma is fundamental. 
\Lemma{2.2}{
For any $a\in(\K^*)^n$, the matrix $J=J(a)$ defined above satisfies the reflection equation
$$
        R^+_{12} J_2 {R^+_{12}}^{t_2} J_1 = J_1 {R^+_{12}}^{t_2} J_2 R^+_{12}
\tag{2.7}
$$
in $\End(V\otimesK V)$, where ${R^+_{12}}^{t_2}$ stands for the matrix obtained
from $R^+_{12}$ by transposition in the second component. 
}
Lemma 2.2 can be checked by direct calculations.  For the reflection equations and related topics, we refer the reader to \cite{Ku}. 
\par
In our context, the meaning of the reflection equation above for $J=J(a)$ can 
be formulated as follows. 
\Proposition{2.3}{
Given a matrix $J={\sum}_{1\le i,j\le N} J_{ij}e_{ij}$ in $\End(V)$, define an element 
$w_J$ in the tensor product $V\otimesK V$ by
$$
w_J = {\sum}_{1\le i,j\le N}  v_i J_{ij} \otimes v_j  \ \ 
\in V\otimesK V. \tag{2.8} 
$$
On the other hand, define a family of elements $M_{ij} (1\le i,j\le N)$ in $\Uq$ as in (2.3).
Then we have  $M_{ij}. w_J =0$ for $1\le i,j \le N$ if and only if the matrix $J$ satisfies
the reflection equation (2.7). 
}
\Proof{
Let us denote by $\bv=(v_1,\cdots,v_N)$ the row vector representing the canonical 
basis for $V$. 
Then the action of the $L$-operators on $\bv$ is described as follows:
$$
        L^+_1 . \bv_2= \bv_2 R^+_{12},  \quad
        S(L^-)^t_1. \bv_2 = \bv_2 {R^+_{12}}^{t_2}. 
\tag{2.9}
$$
The second formula is obtained from 
$(R^-_{12})^{-1} = R^+_{21}$ and ${R^+_{21}}^{t_1}={R^+_{12}}^{t_2}$.  
Noting that $w_J=\bv J \otimes \bv^t$, we compute
$$
\align
        L^+J. w_J &= L^+_1. (\bv_2J_2\otimes \bv^t_2) J_1 
        =  L^+_1.\bv_2J_2 \otimes L^+_1. \bv^t_2 J_1 
\tag{2.10}\\
        & = \bv_2 R^+_{12} J_2\otimes{R^+_{12}}^{t_2} \ \bv^t_2 J_1 
        = \bv_2 R^+_{12} J_2 {R^+_{12}}^{t_2} J_1 \otimes \bv^t_2. 
\endalign
$$
Similarly we have 
$$
\align
        J S(L^-)^t . w_J &= J_1 S(L^-)^t_1 . (\bv_2J_2\otimes\bv^t_2) 
        =  J_1 S(L^-)^t_1.\bv_2J_2 \otimes S(L^-)^t_1. \bv^t_2 
\tag{2.11}\\ 
        & = J_1 \bv_2 {R^+_{12}}^{t_2} J_2\otimes R^+_{12} \bv^t_2
        = \bv_2 J_1 {R^+_{12}}^{t_2} J_2 R^+_{12}\otimes \bv^t_2. 
\endalign
$$
Hence we have  $L^+J. w_J=J S(L^-)^t . w_J$ if and only if the matrix $J$ satisfies the 
reflection equation (2.7). 
}
\par 
        Recall that the tensor product $V\otimesK V$ has the irreducible decomposition 
$ V\otimesK V = V_+ \oplus V_-$ as a left $\Uq$-module into the ``symmetric part"  
$V_+$, isomorphic to the highest weight module $V(2\ep_1) $, 
and the ``anti-symmetric part" $V_-$, isomorphic to $V(\ep_1+\ep_2)$.  
These two components are explicitly given as
$$
\align
        V_+&= {\sum}_{1\le k\le N} \K v_k\otimes v_k +   
        {\sum}_{1\le i<j \le N} \K (q v_i\otimes v_j+ v_j\otimes v_i), \ \text{and}
\tag{2.12}\\
        V_-&= {\sum}_{1\le i<j \le N} \K (v_i\otimes v_j - q v_j\otimes v_i),
\endalign
$$
respectively.  
Note that the element $w_J$ for $J=J(a)$, defined in Proposition 2.3, takes the form 
$$
\align
\text{Case (SO):}& \ \  w_J={\sum}_{k=1}^n  v_k\otimes v_k a_k, 
\tag{2.13}\\
\text{Case (Sp):}&  \ \  w_J=
{\sum}_{k=1}^n  (v_{2k-1}\otimes v_{2k}-q v_{2k}\otimes v_{2k-1})a_k. 
\endalign
$$
The ``quadratic form"  $w_J$  belongs to  $V_+$ in Case (SO), and to $V_-$ in Case (Sp).  
Furthermore, the coideal $\kq$ is so chosen that $\kq$ should annihilate the element $w_J$.  
%
\subsection{2.3. Some remarks on $\/\kq$} 
%
        In the limit as $q\rightarrow 1$, the elements $(q-q^{-1})^{-1} M_{ij}$ recover a $\K$-basis for the Lie subalgebra $\k=\so(n)$ or $\k=\sp(2n)$ , if $a_k=q^{s_k}$ for some 
$s_k\in\Z$ ($1\leq k\leq n$).
We give here explicit formulas of the elements $M_{ij}$ for the comparison with the
case $q=1$.
In Case (SO), they are written as
$$
  M_{ij}=L^+_{ij} - a_i a_j^{-1} S(L^-_{ji})  \quad \text{for}\quad  1\leq i<j\leq n,
\tag{2.14.a}
$$
and $M_{ij}=0$ for $i\geq j$.
Note that the elements $M_{k,k+1}$ are also written as 
$M_{k,k+1}=(q-q^{-1})q^{\ep_k}(f_k - a_k a^{-1}_{k+1} t^{-1}_k e_k)$ 
for $1\le k \le N-1$, where $t_k=q^{\ep_k-\ep_{k+1}}$.  
In Case (Sp), the nonzero elements among $M_{ij}$ are classified 
into the following four groups:
$$
\align
\text{(i)}\quad &M_{2r-1,2r-1}=-M_{2r,2r}=q^{\ep_{2r-1}}-q^{\ep_{2r}},
        \quad (1\leq r\leq n), \tag{2.14.b}\\
\text{(ii)}\quad &M_{2r-1,2r}=L^+_{2r-1,2r},\quad M_{2r,2r-1}=qS(L^-_{2r,2r-1})
        \quad (1\leq r\leq n), \\
\text{(iii)}\quad &M_{2r-1,2s-1}=L^+_{2r-1,2s-1}-a_ra_s^{-1}S(L^-_{2s,2r}) 
        \quad \text{and} \\
&M_{2r,2s}=L^+_{2r,2s}-a_ra_s^{-1}S(L^-_{2s-1,2r-1}) 
        \quad (1\leq r<s\leq n),\\
\text{(iv)}\quad &M_{2r-1,2s}=L^+_{2r-1,2s}+q^{-1}a_ra_s^{-1}S(L^-_{2s-1,2r}) 
        \quad \text{and} \\
&M_{2r,2s-1}=L^+_{2r,2s-1}+qa_ra_s^{-1}S(L^-_{2s,2r-1}) 
        \quad (1\leq r<s\leq n).
\endalign
$$
In order to see what happens as $q$ tends to $1$, 
one has only to note that 
$L^{\pm}_{ij}/(q-q^{-1}) \rightarrow X_{ji}$ 
for $i\lessgtr j$ and that 
$(q^{\ep_i}-q^{-\ep_i})/(q-q^{-1})\rightarrow X_{ii}$,
where $(X_{ij})_{ij}$ is the basis for $\gl(N)$ corresponding to the 
matrix units.
\par
From Proposition 2.1, it follows that the left ideal $\Uq\kq$ and 
the right ideal $\kq\Uq$ are both coideals of $\Uq$. 
As for generator systems of these ideals, we have 
\Proposition{2.4}{
Both the left ideal $\Uq\kq$ and the right ideal $\kq\Uq$ have the following 
generator system:
$$
\align
\text{Case (SO):}\quad
  & M_{k,k+1}=L^+_{k,k+1}-a_ka_{k+1}^{-1}S(L^-_{k+1,k}) \quad (1\leq k\leq n-1).\\
\text{Case (Sp):}\quad 
  & M_{2r-1,2r-1}=q^{\ep_{2r-1}}-q^{\ep_{2r}} \quad  (1\leq r\leq n),\\
  & M_{2r-1,2r}=L^+_{2r-1,2r},\quad M_{2r,2r-1}=qS(L^-_{2r,2r-1}) 
     \quad  (1\leq r\leq n),\\
  & M_{2r,2r+1}=L^+_{2r,2r+1} + qa_ra_{r+1}^{-1}S(L^-_{2r+2,2r-1}) 
     \quad (1\leq r\leq n-1).
\endalign
$$
}
Proposition 2.4 can be proved by using commutation relations between the 
elements $L^+_{ij}$ and $S(L^-_{ij})$:
$$
  S(L^-_2)R^+_{12}L^+_1 = L^+_1R^+_{12}S(L^-_2), \tag{2.15}
$$
which is a direct consequence of (1.12).  
One of the commutation relations is
$$
  qS(L^-_{ji})L^+_{ij}+(q-q^{-1})\sum_{\mu>i}S(L^-_{j\mu})L^+_{\mu j} 
  =qL^+_{ij}S(L^-_{ji})+(q-q^{-1})\sum_{\nu<j}L^+_{i\nu}S(L^-_{\nu i}). 
\tag{2.16}
$$
This relation is related with the recurrence relation (5.33) which will play 
the key role in the computation of radial components of a central element of 
$\Uq$ in Section 5.
Since we will not explicitly use Proposition 2.4 hereafter, 
we omit the detail of its proof. 
We remark that, in Case (SO), 
the generators for the left ideal $\Uq\kq$ in Proposition 2.4 
are the same as Ueno and Takebayashi used to define the quantum analogue of $\GL(n)/\SO(n)$ in \cite{UT}, 
while our coideal $\kq$ gives the whole set of root vectors. 
%
%
\subsection{2.4.  Relation to twisted quantized universal enveloping algebras}
%
It is natural to ask which subalgebra of $\Uq$ is appropriate 
as an object that should play the role of the subalgebra $U(\k)$ 
of $U(\frak{g})$. 
It seems to be a common understanding that the quantized universal 
enveloping algebra $\Uq$ does {\it not} have as many Hopf subalgebras
as the classical $U(\frak{g})$ does. 
This point is discussed in \cite{H2}, in its dual version. 
In this context, it is necessary to take subalgebras of $\Uq$ that are 
{\it not} closed under the coproduct into consideration, 
to recover the degrees of freedom of 
subgroups in quantum groups. 
The arguments in this paper can be reformulated from this point of view, 
namely in the framework of {\it twisted} quantized universal enveloping 
algebras. 
Although we will not use this structure explicitly, 
it should be noted that, 
in the discussion of invariant rings which is a subject from the next section 
on, the central role is played by the left or right ideal generated by $\kq$, not by the coideal $\kq$ itself. 
In Case (SO),
our definition of the coideal $\kq$ is closely related to 
the $q$-deformation of $U(\so(n))$ due to Gavrilik-Klimyk \cite{GK}. 
\par 
Recall that our coideal $\kq$ is defined by using the matrix $M$ formed in (2.3) 
{\it additively} from $L^+$ and $S(L^-)$. 
As a generator system for $\Uq\kq$, this matrix $M$ can be replaced by any 
of the following four matrices:
$$
\align	
-S(L^+) M J  
&= S(L^+)JS(L^-)^t - J \tag{2.17}\\
L^- J^t M^t 
&= L^- J^t (L^+)^t - J^t\\
(L^-)^t QJ^{-1}M 
&= (L^-)^t QJ^{-1} L^+ - QJ^{-1}\\
-S(L^+)^t Q^{-1} M^t (J^{-1})^t
&= S(L^+)^t Q^{-1}(J^{-1})^tS(L^-)-Q^{-1}(J^{-1})^t,
\endalign
$$
where $Q=\diag{q^{2(N-1)},q^{2(N-2)},\cdots,1}$ is the representation
matrix of the group-like element $q^{2\rho}$ on the vector representation. 
The last two equalities in (2.17) are obtained by using the description 
(1.15) of $S^2$. 
We can use any of the four matrices 
$$
S(L^+)JS(L^-)^t, \ L^- J^t (L^+)^t,\ (L^-)^t QJ^{-1} L^+ \ 
\text{and} \ S(L^+)^t Q^{-1}(J^{-1})^t S(L^-). 
\tag{2.18}
$$ 
to define {\it multiplicatively} a subalgebra of $\Uq$ corresponding 
to the subalgebra $U(\frak{k})$ of $U(\frak{g})$. 
To fix the idea, let us take the matrix
$$
	K=S(L^+)J S(L^-)^t	\tag{2.19}
$$
and denote by $U^{\text{tw}}_q(\frak{k})$ the $\K$-subalgebra of 
$\Uq$ generated by the matrix elements $K_{ij}\ (1\le i,j\le N)$ of $K$. 
Then the left ideal $\Uq\kq$ is described as 
$$
\Uq \kq= {\sum}_{a\in U^{\text{tw}}_q(\frak{k})} \Uq(a-\vep(a))
	= {\sum}_{1\le i,j\le N}\Uq(K_{ij}-\vep(K_{ij})). 
\tag{2.20}
$$ 
\par 
This type of subalgebras $U^{\text{tw}}_q(\frak{k})$ defined above 
are analogue of the ``twisted Yangians'' introduced for Yangians by 
G.I.\,Olshanski \cite{O}. 
These ``twisted'' subalgebras have an advantage in the point 
that the commutation relations for the generators are described neatly, 
again by {\it reflection equations}.  
In fact one can show that the matrix $K$ satisfies a reflection equation 
similar to (2.7) (cf. the proof of Proposition 4.4 in Section 4). 
Furthermore, $U^{\text{tw}}_q(\frak{k})$ becomes a {\it coideal} of $\Uq$ 
(except for the condition on the counit).  
\par
In Case (SO), the twisted quantized universal enveloping algebra  $U^{\text{tw}}_q(\frak{so}(n))$ is already found in the work of 
Gavrilik-Klimyk \cite{GK}, although the connection with reflection 
equations is not apparent in their presentation. 
One can show that the subalgebra $U^{\text{tw}}_q(\frak{so}(n))$ is generated 
by the elements $K_{j,j+1}\ (1\le j\le n-1)$ on the subdiagonal:
$$
	K_{j,j+1}=(q-q^{-1})(a_j t^{-1}_je_j-a_{j+1}f_j). \tag{2.21}
$$
Assuming that $a=(q^{n-1},q^{n-2},\cdots,1)$, take the elements 
$$
	\theta_j=f_j - qt^{-1}_j e_j\quad(1\le j\le n-1), \tag{2.22}
$$
so that $K_{j,j+1}=-(q-q^{-1})q^{n-j-1}\theta_j$.  
Then one can check that the generators $\theta_1,\cdots,\theta_{n-1}$ 
of the subalgebra $U^{\text{tw}}_q(\frak{so}(n))$ satisfy the commutation 
relations
$$
\align
\text{\rm(i)}&\ \ 
\theta^2_i\theta_j-(q+q^{-1})\theta_i\theta_j\theta_i
+\theta_j\theta^2_i
=-\theta_j \ \ \text{if}\ \ |i-j|=1, \tag{2.23}\\
\text{\rm(ii)}&\ \ \theta_i\theta_j = \theta_j\theta_i\ \ 
\text{if} \ \ |i-j|>1,
\endalign
$$
as in the definition of $q$-deformation of $U(\frak{so}(n))$ of 
Gavrilik-Klimyk. 
We finally remark that this algebra $U^{\text{tw}}_q(\frak{so}(n))$ 
arises naturally as the commutant of a $q$-analogue of the oscillator 
representation (see \cite{NUW2}). 
%
\section{\S3. Representations with $\kq$-fixed vectors}
%
We now investigate the irreducible representations with $\kq$-fixed vectors. 
As in Section 2, we fix the parameter $a=(a_1,\cdots,a_n)\in (\K^*)^n$ 
and set $J=J(a)$ and $\kq=\kq(a)$.  
%
\subsection{3.1. $\kq$-fixed vectors}
%
Recall that, 
for each dominant integral weight $\ld\in P^+$, 
there exists an irreducible left $\Uq$-module $V(\ld)$ 
with highest weight $\ld$,
uniquely determined up to isomorphism.  
For the coideal $\kq=\kq(a)$ of $\Uq$, defined in Section 2, 
we denote by $V(\ld)_{\kq}$ the vector subspace of all $\kq$-fixed vectors 
in $V(\ld)$:
$$
        V(\ld)_{\kq}:=\{v\in V(\ld); \kq.v=0\}.
\tag{3.1}
$$
In this section, we will prove the following theorem.
\Theorem{3.1}{
{\rm (1)} For any $\ld\in P^+$, one has  $\dim_{\K} V(\ld)_{\kq}\le 1$. 
\newline
{\rm (2)} The left $\Uq$-module $V(\ld)$ has a nonzero 
$\kq$-fixed vector if and only if the dominant integral weight 
$\ld=\sum_{k=1}^N \ld_k\ep_k$ satisfies the following condition:
$$
\align
\text{Case (SO):} & \quad \ld_k-\ld_{k+1}\in 2\Z 
        \quad (1\le k\le n-1),\\
\text{Case (Sp):} & \quad \ld_{2k-1}=\ld_{2k} 
        \quad (1\le k\le n).
\endalign
$$
}
We remark that Theorem 3.1 for Case (SO) is already announced by 
Ueno-Takebayashi \cite{UT}. 
We first prove statement (1) of Theorem 3.1. 
\Lemma{3.2}{ 
Let $v$ be a nonzero $\kq$-fixed vector in $V(\ld)$ $(\ld\in P^+)$.
Decompose $v$ into the sum of weight vectors
$v={\sum}_{\mu\in P} v_{\mu}$, 
so that $ q^h.v_{\mu}=v_{\mu}q^{\langle h,\mu\rangle} $ for all $h\in P^*$.
Then one has $v_{\ld} \neq 0$.
}
\Proof{
By setting, $\widetilde{M}=-J^{-1}MJ$, we take the generators
$$
	\widetilde{M}_{ij}=S(L^-_{ji}) - (J^{-1}L^+ J)_{ij}\quad(1\le i,j\le N)
\tag{3.2}
$$
for the coideal $\kq$.  
Note that, if $i<j$, the element $\widetilde{M}_{ij}$ is nonzero and its 
leading term $S(L^-_{ji})$ has weight $\ep_i-\ep_j$.  
Under the lexicographic order of $P$, let $\mu_0$ be the maximum 
of all $\mu\in P$ such that $v_\mu\neq 0$.  
In the equation $\widetilde{M}_{ij}v=0$ for $i<j$, 
we take the component of weight $\mu_0+\ep_i-\ep_j$,
to obtain $S(L^-_{ji})v_{\mu_0}=0$ for all $i<j$.  
This means that $v_{\mu_0}$ is a highest weight vector of $V(\ld)$. 
Hence we have $\mu_0=\ld$,  namely $v_{\ld}\neq 0$.  
}
Suppose that the $\Uq$-module $V(\ld)$ ($\ld\in P^+$) 
has a nonzero $\kq$-fixed vector $v$.   
Then Lemma 3.2 implies $v_{\ld}\neq 0$.
If $w$ is another $\kq$-fixed vector, there is a constant $c\in\K$ 
such that  $w_{\ld}=cv_{\ld}$ since $\dim_{\K}V(\ld)_{\ld}=1$.
Then the difference $w-cv$ is a $\kq$-fixed vector with $(w-cv)_{\ld}=0$. 
By Lemma 3.2 again, one has $w-cv=0$, i.e., $w=cv$.
This means that $\dim_{\K}V(\ld)_{\kq}\le 1$, as desired.
%
\subsection{3.2. The rank-one case} 
%
Before the proof of statement (2) of Theorem 3.1, we consider 
the case of quantum analogue of the Lie subalgebra $\so(2)$ of $\gl(2)$.  
Although Theorem 3.1 for this case is already known 
by Koornwinder \cite{K1}, we give here a direct proof of this 
statement for completeness (see also \cite{N2}, \cite{NM3}). 
\par
In this case, the coideal $\kq$ is generated by a single element
$$
	M_{12}=(q-q^{-1})q^{\ep_1}(f-a t^{-1}e), \tag{3.3}
$$
where we set $f=f_1,\ e=e_1,\ t=q^{\ep_1-\ep_2}$ and $a=a_1/a_2$. 
Consider the dominant integral weight 
$$
\ld=\ld_1\ep_1+\ld_2\ep_2=\ell\ep_1+\ld_2(\ep_1+\ep_2)
  \ \ \text{with} \ \ \ell=\ld_1-\ld_2\in\N.
\tag{3.4} 
$$
Then $V(\ld)$ is an $(\ell+1)$ dimensional representation and one can 
take a basis $\{u_0,u_1,\cdots,u_\ell\}$ for $V(\ld)$ such that 
$$
  q^h.u_j=q^{\langle h,\ld-j\alpha\rangle} u_j, \  e.u_j=[j] u_{j-1}, \  f.u_j=[\ell-j] u_{j+1}
\tag{3.5}
$$
for $0\le j\le \ell$, where $[j]=(q^j-q^{-j})/(q-q^{-1})$.  
One can easily show that an element 
$$
	v={\sum}_{j=0}^\ell u_j c_j \in V(\ld) \tag{3.6}
$$
satisfies the equation $(f-at^{-1}e).v=0$ if and only if the 
coefficients satisfy the recurrence formula 
$$
	aq^{-\ell+2j+2}[j+1]c_{j+1}=[\ell-j+1]c_{j-1} \tag{3.7}
$$
for $0\le j\le \ell$ with boundary condition $c_{-1}=c_{\ell+1}=0$. 
It is immediately seen that the equation (3.7) has no solution 
if $\ell$ is odd, and that, if $\ell$ is even, the solutions of (3.7) are 
explicitly given by
$$
\align
&c_{2k}=(-1)^k a^{-k} q^{2k(\ell-k)}
\frac{(q^{-2\ell};q^4)_k}{(q^4;q^4)_k} c_0 \ \ (0\le k\le \ell/2),
\tag{3.8}\\
&c_{2k+1}=0 \ \ (0\le k<\ell/2),
\endalign
$$
where $(a;q)_k=(1-a)(1-aq)\cdots(1-aq^{k-1})$.  
This proves Theorem 3.1 for the Case (SO) with $N=n=2$. 
%
\subsection{3.3. Proof of Theorem 3.1.(2)}
%
Next we prove the ``only if'' part of statement (2) of Theorem 3.1 
for the general case. 
\par\medpagebreak\noindent
Case (SO): This assertion is reduced to the rank-one case. 
For each $k$ with $1\le k\le n-1$, we consider the subalgebra 
$U_q(\frak{g}_k)$ of $\Uq$ generated by $q^{\pm\ep_k}$, $q^{\pm\ep_{k+1}}$, 
$e_k$ and $f_k$. 
Note that $M_{k,k+1}\in U_q(\frak{g}_k)$ for $1\le k\le n-1$.  
For a fixed $k$, the $\Uq$-module $V(\ld)$ ($\ld\in P^+$) is 
decomposed into a direct sum
$$
  V(\ld)=W_1\oplus\cdots\oplus W_m \tag{3.9}
$$
of irreducible $U_q(\frak{g}_k)$-submodules. 
We denote by 
$\mu^{(j)}=\mu^{(j)}_k\ep_k + \mu^{(j)}_{k+1}\ep_{k+1}$
the highest weight of $W_j$ for $1\le j\le m$. 
In the decomposition above, we may also assume that 
each $W_j$ is stable under the action of the subalgebra
$U_q(\frak{t})=\K[q^h (h\in P^*) ]$ of $\Uq$; 
decompose the kernel of the operator 
$ e_k : V(\ld) \rightarrow V(\ld)$ by the action of $U_q(\frak{t})$, 
if necessary. 
Let now $v$ be a nonzero $\kq$-fixed vector in $V(\ld)$ 
and decompose it in the form
$$
 v=w_1+\cdots+w_m \ \  \text{with}\ \  w_j\in W_j\,\,(1\le j\le m).
\tag{3.10}
$$
Then the summands $w_j$ are annihilated by $M_{k,k+1}$,
since $W_j \ (1\le j \le m)$ are all $U_q(\frak{g}_k)$-submodules. 
Hence one sees that $\langle\mu^{(j)},\alpha_k\rangle\in 2\Z$ if $w_j\neq 0$,
by the result of the rank-one case.
Taking the component of weight $\ld$ of $v$, one has 
$ v_{\ld}=(w_1)_{\ld}+\cdots+(w_m)_{\ld}.$
Since $v_{\ld}\neq 0$ by Lemma 3.2, one has $(w_j)_{\ld}\neq 0$
and $e_k. (w_j)_{\ld}=0$ for some $j$. 
Hence, $\ld_k\ep_k+\ld_{k+1}\ep_{k+1}$ is 
the highest weight of $W_j$. 
This means that 
$\ld_k\ep_k+\ld_{k+1}\ep_{k+1}=\mu^{(j)}$. 
Hence, one has 
$\langle\ld,\alpha_k\rangle\in 2\Z$ for $1\le k\le n-1$.
\par\medpagebreak\noindent
Case (Sp): 
Let $v$ be a nonzero $\kq$-fixed vector in $V(\ld)$ and decompose it 
into the sum of weight vectors $v={\sum}_{\mu} v_{\mu}$, where 
$v_{\ld}\neq 0$ by Lemma 3.2. 
Recall that the coideal $\kq$ contains the elements 
$M_{2r-1,2r-1}=q^{\ep_{2r-1}}-q^{\ep_{2r}}$ ($1\le r\le n$).
Since $(q^{\ep_{2r-1}}-q^{\ep_{2r}}).v_{\ld}
=(q^{\ld_{2r-1}}-q^{\ld_{2r}})v_{\ld}=0$, one has 
$\ld_{2r-1}=\ld_{2r}$ as desired.
\par\medpagebreak\noindent
Thus we have proved the ``only if'' part of Theorem 3.1.(2).
\par\medpagebreak
\par
        The ``if'' part of Theorem 3.1.(2) is proved in 
a constructive manner. 
From now on, we denote by $P^+_{\k}$ the set of all dominant integral
weights satisfying the condition of Theorem 3.1.(2):
$$
\align
\text{Case (SO):}
  &\quad P^+_{\k}:=\{\ld\in P^+; \langle\ld,\alpha_k\rangle\in 2\Z\quad(1\le k\le n-1)\},
  \tag{3.11}\\
\text{Case (Sp):}
  &\quad P^+_{\k}:=\{\ld\in P^+; \langle\ld,\alpha_{2k-1}\rangle=0\quad(1\le k\le n) \}.
\endalign
$$
Denoting the fundamental weights by 
$\Ld_r={\sum}_{k=1}^r \ep_k\quad (1\le r\le N)$,
we have the following alternative expression of $P^+_\k$:
$$\align
\text{Case (SO):}\quad
  & P^+_{\k}={\sum}_{r=1}^{n-1}2\N\Ld_r + \Z\Ld_n,
  \tag{3.12}\\
\text{Case (Sp):}\quad 
  & P^+_{\k}={\sum}_{r=1}^{n-1}\N\Ld_{2r} + \Z\Ld_{2n},
\endalign
$$
where $\N=\{0,1,2,\cdots \}$.
It is clear that the one dimensional representations $V(\ell\Ld_N)$($\ell\in\Z$) 
have nonzero $\kq$-fixed vectors. 
\par 
        We start with constructing nonzero $\kq$-fixed vectors in 
$V(2\Ld_r)$ for Case (SO) and in $V(\Ld_{2r})$ in Case (Sp), for $1\le r\le n$. 
For this purpose we make use of the $q$-exterior algebra 
${\bigwedge}_q(V) $; 
it is the quotient algebra of the tensor algebra 
$T(V)={\bigoplus}_{d=0}^\infty V^{\otimes d} $
modulo the two-sided ideal generated by the ``symmetric part'' $V_+$ 
of $V\otimesK V$.  
Let us denote the multiplication in this $q$-exterior algebra 
by $\wedge$.
Then ${\bigwedge}_q(V)$ is the $\K$-algebra generated by the 
elements $v_1,\cdots,v_N$ with fundamental relations 
$$
  v_k\wedge v_k=0 \quad(1\le k\le N)\quad\text{and}\quad 
  qv_i\wedge v_j+v_j\wedge v_i=0\quad(1\le i<j\le N).
\tag{3.13}
$$
Note that the $q$-exterior algebra ${\bigwedge}_q(V)$ is an algebra 
with $\Uq$-symmetry
generated by the vector representation $V={\bigoplus}_{k=0}^N \K v_k$.  
Namely the multiplication 
${\bigwedge}_q(V)\otimesK{\bigwedge}_q(V)\to {\bigwedge}_q(V)$ and the unit
homomorphism $\K\to{\bigwedge}_q(V)$ 
are $\Uq$-homomorphisms. 
Furthermore, it decomposes as 
$$
  {\bigwedge}_q(V) = {\bigoplus}_{r=0}^N {\bigwedge}^r_q(V) 
\quad\text{with}\quad {\bigwedge}^r_q(V)@<\sim<<V(\Ld_r),
\tag{3.14}
$$
into irreducible components. 
Note also that, for each $0\le r \le N$, the $\Uq$-submodule ${\bigwedge}^r_q(V)$ 
has the basis 
$$
	v_{k_1}\wedge v_{k_2}\wedge\cdots\wedge v_{k_r}\quad (1\le k_1<k_2<\cdots<k_r\le n). 
\tag{3.15}
$$
\Lemma{3.3.A}{
In Case (SO), the $\Uq$-module 
${\bigwedge}^r_q(V)\otimesK {\bigwedge}^r_q(V)$ has the $\kq$-fixed
vector
$$
  w_r:={\sum}_{1\le k_1<\cdots<k_r\le N} 
v_{k_1}\wedge \cdots\wedge v_{k_r}\otimes 
v_{k_1}\wedge \cdots\wedge v_{k_r} a_{k_1}\cdots a_{k_r} 
\tag{3.16}
$$
for $1\le r\le n$.
Hence $V(2\Ld_r)$ has a nonzero $\kq$-fixed vector 
for $1\le r \le n$.
}
\Proof{
In order to obtain a nonzero $\kq$-fixed vector in
${\bigwedge}^r_q(V)\otimesK {\bigwedge}^r_q(V)$,
we construct an intertwining operator 
$$
\Phi: (V\otimesK V)^{\otimes r} @>\sim>> V^{\otimes r}\otimesK V^{\otimes r}. 
\tag{3.17}
$$
Recall that there is a $\Uq$-isomorphism 
$s_{12}: V_1\otimesK V_2 @>\sim>> V_2\otimes V_1$,
for $V_1=V_2=V$, whose matrix representation is given by
$\check{R}_{12}=R^+_{12}P_{12}$. 
In the notation $\bv=(v_1,\cdots,v_n)$ as in the proof of 
Proposition 2.3, this isomorphism can be described as
$$
s_{12}(\bv_1\otimes\bv_2)=\bv_2\otimes\bv_1 R^+_{21}, 
\ \ \text{or equivalently,}\ \ 
s_{12}(\bv^t_1\otimes\bv_2)=\bv_2\otimes{R^+_{12}}^{t_2}\bv^t_1 .   
\tag{3.18}
$$
By composing isomorphisms of this type repeatedly, we obtain 
an isomorphism
$$
\Phi: (V_1\otimes V_{1'})\otimes \cdots\otimes(V_r\otimes V_{r'})
@>\sim>>(V_1\otimes\cdots\otimes V_r)\otimes (V_{1'}\otimes\cdots\otimes V_{r'})
\tag{3.19}
$$
for $V_k=V_{k'}=V\ (1\le k\le r)$; here we take 
$$
\Phi=s_{1'r}\circ s_{2'r}\circ\cdots\circ s_{(r-1)',r}\circ s_{1',r-1}\circ\cdots
	\circ s_{1'3}\circ s_{2'3}\circ s_{1'2}.
\tag{3.20}
$$
Note that the element $(w_J)^{\otimes r}$ in $(V\otimesK V)^{\otimes r}$
is a $\kq$-fixed vector, since $\kq$ is a coideal. 
By this isomorphism (3.20), the $\kq$-fixed vector  $(w_J)^{\otimes r}$ is transformed into 
$$
\align
\Phi(w_J\otimes\cdots\otimes w_J)
&= \Phi((\bv_1J_1\otimes\bv^t_1)\otimes\cdots\otimes(\bv_rJ_r\otimes\bv^t_r))
\tag{3.21}\\
&= \bv_1\otimes\cdots\otimes\bv_r
K_{1\cdots r}\otimes\bv^t_1\otimes\cdots\otimes\bv^t_r,
\endalign
$$
where
$$
K_{1\cdots r}=J_1 {R^+_{12}}^{t_2} J_2 {R^+_{13}}^{t_3} {R^+_{23}}^{t_3} J_3
\cdots J_{r-1}{R^+_{1r}}^{t_r}\cdots{R^+_{r-1,r}}^{t_r} J_r. 
\tag{3.22}
$$
Note here that, as to the coefficients of the $R$-matrix, 
we have   
$$
	({R^+_{12}}^{t_2})^{ik}_{j\ell} = (R^+_{12})^{i\ell}_{jk} = \delta_{ij}\delta_{k\ell}
\quad\text{if}\ \ i\neq k. 
\tag{3.23}
$$
Hence we see that, if the indices $i_1,\cdots,i_r$ are mutually distinct, then 
$$
	(K_{1\cdots r})^{i_1\cdots i_r}_{j_1\cdots j_r} 
	= \delta_{i_1j_1}\cdots\delta_{i_rj_r} a_{i_1}\cdots a_{i_r}. 
\tag{3.24}
$$
Denoting by $\pr_r : V^{\otimes r}\rightarrow {\bigwedge}^r_q(V)$ the canonical
projection, we now consider the $\Uq$-homomorphism
$$
\Psi=(\pr_r\otimes\pr_r)\circ\Phi : 
(V\otimesK V)^{\otimes r} \rightarrow {\bigwedge}^r_q(V)\otimesK {\bigwedge}^r_q(V). 
\tag{3.25}
$$
Then, by (3.21) and (3.24), we get
$$
\align
\Psi(w_J\otimes\cdots\otimes w_J)
&= \bv_1\wedge\cdots\wedge\bv_r
K_{1\cdots r}\otimes\bv^t_1\wedge\cdots\wedge\bv^t_r
\tag{3.26}\\
&={\sum}_{1\le k_1, \cdots, k_r\le n} 
v_{k_1}\wedge \cdots\wedge v_{k_r}\otimes 
v_{k_1}\wedge \cdots\wedge v_{k_r} a_{k_1}\cdots a_{k_r} \\
& = [r]_{q^2}! w_r,
\endalign
$$
where $[r]_{q^2}!=(q^2;q^2)_r/(1-q^2)^r$.  
This shows that $w_r$ of (3.16) is a $\kq$-fixed vector. 
Since ${\bigwedge}^r_q(V)$ is isomorphic to $V(\Ld_r)$, there is a nontrivial 
$\Uq$-homomorphism 
${\bigwedge}^r_q(V)\otimesK{\bigwedge}^r_q(V)\rightarrow V(2\Ld_r)$. 
The image of $w_r$ by this homomorphism gives a nonzero $\kq$-fixed vector 
in $V(2\Ld_r)$ since $(w_r)_{2\Ld_r}\neq 0$.
}
We remark that Lemma 3.3.A can be proved also by chasing directly 
the action of $M_{k,k+1} (1\le k\le n-1)$ on $w_r$ (see Proposition 2.4).  
The intertwining operators $\Phi$ and $\Psi$ above will be used again later in the 
discussion of the invariant ring in Section 4. 
\Lemma{3.3.B}{
In Case (Sp), the element
$$
 w_r={\sum}_{1\le k_1<\cdots<k_r\le n}
        v_{2k_1-1}\wedge v_{2k_1}\wedge \cdots\wedge v_{2k_r-1}\wedge v_{2k_r}
        a_{k_1}\cdots a_{k_r},
\tag{3.27}
$$
gives a $\kq$-fixed vector in $\bigwedge^{2r}_q(V)$. 
Hence, $V(\Ld_{2r})$ has a nonzero $\kq$-fixed vector
for $1\le r \le n$.
}
\Proof{
Note first that the projection $\pr_2: V\otimesK V\rightarrow {\bigwedge}^2_q(V)$ 
maps the $\kq$-fixed vector $w_J$ of (2.13) to 
$$
	\pr_2(w_J) = (1+q^2) {\sum}_{k=1}^n v_{2k-1}\wedge v_{2k}a_k = (1+q^2)w_1. 
\tag{3.28}
$$
Hence $w_1$ is a $\kq$-fixed vector in ${\bigwedge}^r_q(V)$. 
We now compute the $r$-th power of $w_1$ in the $q$-exterior algebra ${\bigwedge}_q(V)$
using the $q$-binomial theorem, to obtain
$$
\align
 (w_1)^{\wedge r} &=({\sum}_{k=1}^n v_{2k-1}\wedge v_{2k}a_k)^{\wedge r} \tag{3.29}\\
	&= [r]_{q^4}! {\sum}_{1\le k_1<\cdots<k_r\le n} 
        v_{2k_1-1}\wedge v_{2k_1}\wedge\cdots\wedge v_{2k_r-1}\wedge v_{2k_r}
        a_{k_1}\cdots a_{k_r}\\
	&=[r]_{q^4}! w_r. 
\endalign
$$
Since $\kq$ is a coideal, equality (3.29) shows that $w_r$ is a $\kq$-fixed vector in ${\bigwedge}^r_q(V)$ for $1\le r \le n$. 
}
        In each case, we have now a system of generators 
of the monoid $P^+_{\k}$, consisting of weights $\ld$ with 
$V(\ld)_{\k}\neq0$. 
In order to complete the proof of Theorem 3.1.(2), we have only 
to show that, for any $\ld,\mu\in P^+$, the tensor product 
$V(\ld)\otimesK V(\mu)$ has a nonzero $\kq$-fixed vector, if both 
$V(\ld)$ and $V(\mu)$ do. 
This follows from the fact that there exists a nontrivial 
$\Uq$-homomorphism $V(\ld)\otimesK V(\mu)\rightarrow V(\ld+\mu)$,
up to a scalar multiple.
If $v$ and $w$ are nonzero $\kq$-fixed vectors of $V(\ld)$ and
$V(\mu)$ respectively, the tensor product $v\otimes w$ 
is annihilated by $\kq$ since $\kq$ is a coideal. 
Since $(v\otimes w)_{\ld+\mu}=v_{\ld}\otimes w_{\mu}\neq 0$,
the image of $v\otimes w$ in $V(\ld+\mu)$ then gives 
a nonzero $\kq$-fixed vector.
\par
        We have thus proved that $V(\ld)$ has a nonzero $\kq$-fixed
vector if and only if $\ld\in P^+_{\k}$ and that 
${\dim}_{\K}V(\ld)_{\kq}=1$ for each $\ld\in P^+_{\k}$.
%
\subsection{3.4. Passage from the left to the right}
%
In this section, we have discussed $\kq$-fixed vectors in left $\Uq$-modules. 
The same argument naturally applies to right $\Uq$-modules. 
The passage from left $\Uq$-modules to right $\Uq$-modules
can be described, functorially, by the $\ast$-operation of $\Uq$ 
explained in Section 1.4.
\par
	For a left $\Uq$-module $M$, let us denote by $M^\circ$ 
the right $\Uq$-module obtained from $M$ by regarding it as a right $\Uq$-module 
through this involution: 
$$
	x.a = a^*.x \quad\text{ for all}\ \ a\in\Uq\ \ \text{ and }\ \ x\in M.
\tag{3.30}
$$ 
Then, for each dominant integral weight $\ld\in P^+$, the right $\Uq$-module 
$V(\ld)^\circ$ gives rise to the irreducible right $\Uq$-module 
of highest weight $\ld$. 
\par
We now look at the coideal $\kq(a)\  (a\in (\K^*)^n) $ of $\Uq$. 
One can easily compute the action of the involution on the matrix
$M(a)=L^+ - J(a)S(L^-)^t J(a)^{-1}$, 
by using $J(a)^{-1}=J(a^{-1})$ in Case (SO), 
or $J(a)^{-1}=-q^{-1}J(a^{-1})$ in Case (Sp). 
In fact we have 
$$
			M(a)^* = -(J(a) M(a^{-1}) J(a)^{-1})^t, 
\tag{3.31} 
$$
where the left-hand side denotes the matrix $({M(a)_{ji}}^*)_{1\le i,j\le N}$. 
Hence we have $\kq(a)^* = \kq(a^{-1})$ for any $a\in(\K^*)^n$. 
This implies that each $\kq(a)$-fixed vector in a left $\Uq$-module $M$ 
can be read as a $\kq(a^{-1})$-fixed vector in the right $\Uq$-module $M^\circ$. 
%
\section{\S4. Quantum homogeneous spaces and zonal spherical functions}
%
This section is devoted to the study of quantum analogue of the coset spaces 
$G/K$ and $K\backslash G/K$, 
for the closed subgroup $K$ of $G=\GL(N)$, corresponding to
the Lie subalgebra $\k$ of $\g=\gl(N)$ of (2.1). 
We fix the parameter $a\in(\K^*)^n$ involved in the definition of 
the coideal $\kq(a)$ and set $J=J(a), M=M(a), \kq=\kq(a)$. 
%
\subsection{4.1. Quantum analogue of the homogeneous space $G/K$}
%
In terms of the coideal $\kq=\kq(a) \ (a\in (\K^*)^n)$, we can study the quantum analogue 
of the homogeneous space $G/K$ for the following closed subgroup $K$: 
$$
\align
\text{Case (SO):}& \ \  K=\SO(n)=\{g\in\GL(n)\ ; \ g g^t =\id_n, \det(g)=1 \} \ \  (N=n)
\tag{4.1}\\
\text{Case (Sp):}& \ \  K=\SP(2n)=\{g\in\GL(2n)\ ; \ gJ_n g^t=J_n \} \ \  (N=2n). 
\endalign 
$$
Recall that the ``coordinate ring" $\Aq$ of the quantum general linear group $\GL_q(N)$ is
a $\K$-algebra with two-sided $\Uq$-symmetry. 
By means of the bimodule structure over $\Uq$, we consider 
the $\K$-vector subspace of $\kq$-invariant elements in $\Aq$
under the left action of $\Uq$:
$$
        A_q(G/K) := \{\varphi \in \Aq;\ \kq.\varphi=0 \}. 
\tag{4.2}
$$
Thanks to the fact that $\kq$ is a coideal (Proposition 2.1), 
this subspace actually becomes $\K$-subalgebra of $\Aq$.
Note also that the subalgebra $A_q(G/K)$ is 
a left $\Aq$-subcomodule, hence a right $\Uq$-submodule of $\Aq$. 
The algebra $A_q(G/K)$ is a $\K$-algebra with right $\Uq$-symmetry, 
and is regarded as the algebra of 
regular functions on the left quantum $G_q$-space $(G/K)_q$. 
\par
Recall that $\Aq$ has the irreducible decomposition 
$$
        \Aq={\bigoplus}_{\ld\in P^+}\ W(\ld),
\tag{4.3}
$$
as a $\Uq$-bimodule. 
Here $W(\ld)$ is the $\K$-subspace of $\Aq$ spanned by the matrix elements of 
the irreducible right $\Aq$-comodule $V(\ld)$. 
Note also that the $\Uq$-bimodule $W(\ld)$ is isomorphic to the tensor product 
of $\Hom(V,\K)$, regarded as a right $\Uq$-module, and the left $\Uq$-module $V(\ld)$.
We denote by $V(\ld)^\circ$ the irreducible right $\Uq$-module with highest weight 
$\ld$, so that $W(\ld) @<\sim<<V(\ld)^\circ\oxK V(\ld)$. 
From the irreducible decomposition (4.3) and Theorem 3.1, we obtain the multiplicity free decomposition of the invariant ring 
$A_q(G/K)$. 
\Proposition{4.1}{
The $\K$-subalgebra $A_q(G/K)$ of left $\kq$-invariants in $\Aq$ decomposes into the form 
$$
  A_q(G/K)\;@<\sim<<  {\bigoplus}_{\ld\in P^+_\k}V(\ld)^\circ
\tag{4.4}
$$
as a right $\Uq$-module. 
}
\Proof{
For each dominant integral weight $\ld\in P^+$, 
we have $W(\ld)@<\sim<<V(\ld)^\circ\oxK V(\ld)$ as a $\Uq$-bimodule. 
Hence $W(\ld)_{\kq} @<\sim<<V(\ld)^\circ\oxK V(\ld)_{\kq}$ as a right 
$\Uq$-module. 
On the other hand, we know that $\dim_\K V(\ld)_{\kq}=1$ if $\ld \in P^+_\k$ and 
$V(\ld)_{\kq}=0$ otherwise, from Theorem 3.1. 
Hence we have the multiplicity free decomposition (4.4) by the irreducible decomposition 
(4.3) of $\Aq$. 
}
For the description of the invariant ring $A_q(G/K)$, we define a family of 
quadratic elements $x_{ij}$\, $(1\le i,j\le N)$ in $\Aq$ 
by the formula
$$
    X=(x_{ij})_{1\le i,j\le N},\quad X=TJ T^t.
\tag{4.5}
$$
These elements $x_{ij}$ are explicitly written as
$$
        x_{ij} = {\sum}_{k=1}^n\; t_{ik} t_{jk} a_k \quad(1 \le i,j\le n) 
\tag{4.6.a}
$$
in Case (SO) and
$$
        x_{ij} = {\sum}_{k=1}^n \;(t_{i,2k-1} t_{j,2k}-q t_{i,2k} t_{j,2k-1}) a_k
\quad(1\le i,j\le 2n)
\tag{4.6.b}
$$
in Case (Sp). 
\Lemma{4.2}{
The elements $x_{ij} \in \Aq$ defined as above are invariant under the left action of the coideal 
$\kq$. 
Namely, $x_{ij}\in A_q(G/K)$ for all $1 \le i,j\le N$. 
}
\Proof{
Note that from (1.21) we have 
$$
        L^+_1. T_2 = T_2 R^+_{12},  \quad
        S(L^-)^t_1.T_2  = T_2 {R^+_{12}}^{t_2}, 
\tag{4.7}
$$
since $(R^-_{12})^{-1}=R^+_{21}$ and ${R^+_{21}}^{t_1}={R^+_{12}}^{t_2}$. 
Hence we have 
$$
\align
        L^+_1. X_2&=  L^+_1. T_2  J_2 L^+_1.T^t_2
        =T_2 R^+_{12} J_2 {R^+_{12}}^{t_2} T^t_2
\tag{4.8}\\
        S(L^-)^t_1. X_2 &= S(L^-)^t_1.T_2 J_2 S(L^-)^t_1.T^t_2
        =T_2 {R^+_{12}}^{t_2}J_2 R^+_{12} T^t_2
\endalign
$$
Since the matrix $J$ satisfies the reflection equation (2.7), we have
$$
        (L^+_1J_1- J_1S(L^-)^t_1).X_2=0. \tag{4.9}
$$
This means that $M_1.X_2=0$,  namely,  $\kq.X=0$.
}
As to the structure of this invariant ring, we have
\Theorem{4.3}{
The $\K$-subalgebra $A_q(G/K)$ of left $\kq$-invariants in $\Aq$ 
is generated by the quadratic elements 
$x_{ij}$\,$(1\le i,j\le N)$ defined by (4.5), 
together with $\detq(T)^{\pm 1}$. 
}
The proof of Theorem 4.3 will be given later in Section 4.4. 
As we will see later, the square $\detq(T)^2$ of the quantum determinant 
is actually represented by $x_{ij}$'s in Case (SO), and, in Case (Sp), 
so is $\detq(T)$ (see Remark 4.12 below). 
Note that, in each case, $\detq(T)$ lies in the center of $A_q(G/K)$. 
For the present, we will show that the commutation relations among 
the quadratic elements $x_{ij}\ (1\le i,j\le N) $ are described again by the 
reflection equation. 
\Proposition{4.4} {
The elements $x_{ij}$\,$(1\le i,j\le N)$ defined above satisfy 
$$\align
\text{Case (SO)}:
&\quad x_{ij}=qx_{ji} \;(1\le i<j\le n)
\tag{4.10}\\
\text{Case (Sp)}:
&\quad x_{ii}=0 \; (1\le i\le 2n),\ \ qx_{ij}+x_{ji}=0 \; (1\le i<j\le 2n).
\endalign
$$
In each case, they have the commutation relations
$$ 
R^+_{12}\,X_2\,{R^+_{12}}^{t_2}\,X_1=
X_1\,{R^+_{12}}^{t_2}\,X_2R^+_{12}.
\tag{4.11}
$$
}
\Proof{
In each case, formula (4.10) can be shown directly by using 
the expression (4.6) and 
the commutation relations (1.3) of 
the ``coordinates" $t_{ij}$. 
We show that the matrix $X$ satisfies the reflection equation (4.11). 
Note that commutation relations (1.4) for $T$ implies 
$$
        T^t_2 {R^+_{12}}^{t_2} T_1=T_1 {R^+_{12}}^{t_2} T^t_2, \ 
        T_2 {R^+_{12}}^{t_2} T^t_1=T^t_1 {R^+_{12}}^{t_2} T_2, \ 
        R^+_{12} T^t_2 T^t_1=T^t_1 T^t_2 R^+_{12}. 
\tag{4.12}
$$
(Note that ${R^+_{12}}^t=R^+_{21}$ and ${R^+_{12}}^{t_2}={R^+_{21}}^{t_1}$.)
By using (1.4), (4.12) and the reflection equation for $J$, we obtain 
$$
\align
        R^+_{12}X_2 {R^+_{12}}^{t_2} X_1 
        &=R^+_{12}T_2 J_2T^t_2 {R^+_{12}}^{t_2}\ T_1 J_1T^t_1 
                =R^+_{12}T_2 J_2T_1 {R^+_{12}}^{t_2}\ T^t_2 J_1T^t_1\tag{4.13}\\
        &=R^+_{12}T_2 T_1 J_2{R^+_{12}}^{t_2} J_1T^t_2 T^t_1 
                =T_1 T_2 R^+_{12}J_2{R^+_{12}}^{t_2} J_1T^t_2 T^t_1\\
        &=T_1 T_2 J_1{R^+_{12}}^{t_2} J_2R^+_{12}T^t_2 T^t_1 
                =T_1 T_2 J_1{R^+_{12}}^{t_2} J_2T^t_1 T^t_2 R^+_{12}\\
        &=T_1 J_1T_2 {R^+_{12}}^{t_2}\ T^t_1 J_2T^t_2R^+_{12} 
                =T_1 J_1T^t_1 {R^+_{12}}^{t_2}\ T_2 J_2T^t_2R^+_{12}\\
        &=X_1 {R^+_{12}}^{t_2} X_2R^+_{12}, 
\endalign
$$
as desired. 
}
This description of the invariant ring corresponds to the realization of 
the homogeneous space $\SL(n)/\SO(n)$ 
and $\SL(2n)/\SP(2n)$ as an orbit of symmetric and skew-symmetric matrices,
respectively. 
The quadratic elements $x_{ij}\ (1\le i,j\le N)$ above can be thought of as the quantum 
analogue of the coordinates for the space of symmetric or skew-symmetric matrices. 
\Remark{4.5}{
By the right action of $\Uq$ on $\Aq$, we can also consider the quantum analogue of the 
right homogeneous space $K\backslash G$. 
The subalgebra of right $\kq$-invariants 
$$
        A_q(K\backslash G):=\{\varphi\in\Aq ; \varphi.\kq=0 \} 
\tag{4.14}
$$
of $\Aq$ is a $\K$-algebra with left $\Uq$-symmetry, and has the multiplicity free 
decomposition 
$$
  A_q(K\backslash G)\;@<\sim<<  {\bigoplus}_{\ld\in P^+_\k}V(\ld)
\tag{4.15}
$$
as a left $\Uq$-module, similarly to (4.4). 
In this case we define the quadratic elements $y_{ij}\ \ (1\le i,j\le N)$ by the formula
$$
        Y=(y_{ij})_{1\le i,j\le N},\quad Y=T^t J^{-1} T. 
\tag{4.16}
$$
For these elements $y_{ij}$, 
we have exactly the same statement as Proposition 4.4.   
This follows from that fact that  $J^{-1}=J(a)^{-1}$ 
satisfies the same reflection equation (2.7); 
$J(a)^{-1}$ is a constant multiple of $J(a^{-1})$. 
The invariant ring $A_q(K\backslash G)$ is then generated by the quadratic elements 
$y_{ij}\ (1\le i,j\le N)$ and $\detq(T)^{\pm 1}$. 
We remark that all these properties of $A_q(K\backslash G)$
are obtained from those of $A_q(G/K)$ by using involutions 
as we explained in Section 3.4.
}
%
\subsection{4.2. Quantum analogue of the double coset space $K\backslash G/K$}
%
The next step is to study the quantum analogue of the double coset space 
$K\backslash G/K$. 
We consider the following $\K$-subalgebra of $\kq$-biinvariant 
elements in $\Aq$ :
$$
 {\Cal H}=A_q(K\backslash G /K)
    :=\{\varphi\in \Aq;\, \kq.\varphi=\varphi.\kq=0 \}. 
\tag{4.17}
$$
From the irreducible decomposition (4.3) of $\Aq$ again, 
we have the direct decomposition of the $\K$-subalgebra $\Cal{H}$ 
of $\kq$-biinvariants
$$
  \Cal{H}={\bigoplus}_{\ld\in P^+_{\k}} \Cal{H}(\ld), 
  \quad\text{with}\ \Cal{H}(\ld)=\Cal{H}\cap W(\ld). 
\tag{4.18}
$$
We remark that $\dim_\K\Cal{H}(\ld)=1$ if $\ld\in P^+_{\kq}$ and $\Cal{H}(\ld)=0$ otherwise,
as is immediately seen from Theorem 3.1 and its  right $\Uq$-module version.  
We say that a nonzero element $\varphi$ of 
$\Cal{H}(\ld)\,(\ld\in P^+_{\k})$ 
is a {\it zonal spherical function\/} associated with the representation 
$V(\ld)$. 
Recall that the subspace $W(\ld)$ of matrix elements of the right $\Aq$-comodule 
$V(\ld)$ is characterized as the simultaneous eigenspace of the center 
of $\Uq$:
$$
W(\ld)=\{\varphi\in\Aq;\; C.\varphi=\chi_{\ld}(C)\varphi
\ \ \text{for any central }\,  C\in\Uq\},
\tag{4.19}
$$
where $\chi_{\ld}(C)$ denotes the eigenvalue of the central element $C\in\Uq$ on the 
irreducible representation $V(\ld)$. 
From this, we see that a nonzero element $\varphi \in \Aq$ is a zonal spherical function 
associated with $V(\ld)$ if and only if 
\roster
        \item $\kq.\varphi =\varphi.\kq = 0$, and 
        \item $C.\varphi=\chi_{\ld}(C)\varphi\ $ for any central element $C\in\Uq$. 
\endroster
Summarizing these remarks, we have
\Proposition{4.6}{
The subalgebra $\Cal{H}=A_q(K\backslash G/K)$ of $\kq$-biinvariant elements in
$\Aq$ has the simultaneous eigenspace decomposition 
$\Cal{H}={\bigoplus}_{\ld\in P^+_\k }\Cal{H}(\ld)  $
under the action of the center of $\Uq$.  
Furthermore, the simultaneous eigenspace $\Cal{H}(\ld)$ is one-dimensional for each 
$\ld\in P^+_\k$. 
}
\par
For the description of the zonal spherical functions $\varphi$, we consider 
their ``restriction'' $\varphi|_{\T}$ to the diagonal subgroup $\T=(\K^*)^N$ 
of the quantum group $GL_q(N)$. 
Recall that the quantum general linear group $\GL_q(N)$ ``contains" the $N$-dimensional
algebraic torus $\T=(\K^*)^N$ on its diagonal. 
Let $z=(z_1,\cdots,z_N)$ be the canonical coordinates of $\T$; 
the coordinate ring $A(\T)$ is the $\K$-algebra of Laurent 
polynomials $\K[z^{\pm 1}]=\K[{z_1}^{\pm 1},\cdots,{z_N}^{\pm 1}]$. 
Then there exists a unique Hopf algebra homomorphism 
$\varphi\mapsto \varphi|_{\T}: \Aq \rightarrow A(\T)$ 
such that $t_{ij}|_{\T}=\delta_{ij}z_j$ for $1\le i,j\le N$. 
\par
The restriction of  $\kq$-biinvariant ``functions'' on $\GL_q(N)$ 
to the diagonal subgroup $\T$
is described by the composition of $\K$-algebra homomorphisms 
$$
\Cal{H}=A_q(K\backslash G/K) \hookrightarrow \Aq \rightarrow A(\T). 
\tag{4.20}
$$ 
\Theorem{4.7}{
The restriction mapping $\Cal{H}=A_q(K\backslash G/K) \rightarrow A(\T)$ is 
an injective $\K$-algebra homomorphism; hence, $\Cal{H}$ is a commutative 
$\K$-subalgebra of $\Aq$. 
Furthermore the image $\Cal{H}|_{\T}$ of $\Cal{H}$ is given by 
$$
\align
\text{Case (SO):}\quad&\Cal{H}|_{\T}=
\K[{z_1}^2,\cdots,{z_n}^2]^{{\frak S}_n}[(z_1\cdots z_n)^{-1}]\tag{4.21}\\
\text{Case (Sp):}\quad&\Cal{H}|_{\T}= 
\K[z_1z_2,\cdots,z_{2n-1}z_{2n}]^{{\frak S}_n}[(z_1z_2\cdots z_{2n})^{-1}]. 
\endalign
$$
}
The proof of Theorem 4.7 will be given in Section 4.5. 
Theorem 4.7 for Case (SO) is stated also in \cite{UT}. 
%
\subsection{4.3. $\kq$-invariant matrix elements}
%
In order to investigate the invariant rings $A_q(G/K)$ and $A_q(K\backslash G/K)$, 
we study the $\kq$-invariance of matrix elements of irreducible representations $V(\ld)$. 
For this purpose, we will make use of the unitarizability of $\Aq$-comodules, 
with respect to the Hopf $\ast$-algebra structure of $\Aq$ 
explained in Section 1.4. 
In our setting $\K=\Q(q)$, we take the conjugation on $\K$ to be the identity mapping, which means that $q$ is ``real''.
\par
        Let $M$ be an arbitrary finite dimensional right $\Aq$-comodule $M$, 
and $\rho_G : M\rightarrow M\oxK \Aq$ its right coaction.  
We use the subscript $G$ for the coaction 
to remember that this structure corresponds to 
a group representation. 
It is known that there exists a nondegenerate hermitian form 
$ \langle\ ,\ \rangle  : M\times M \to \K$, conjugate linear in the first argument, 
which is invariant with respect to the quantum unitary group 
$\text{\rm U}_q(N)$ in the sense
$$
\langle\rho_G(u), \rho_G(v) \rangle  = \langle u,v\rangle .1 \quad\text{for any}\ \ u,v\in M. 
\tag{4.22}
$$
Here, we use the same notation $\langle\ ,\ \rangle $ to refer to the hermitian form 
$(M\oxK\Aq)\times(M\oxK\Aq) \to \Aq$, naturally defined on
$M\oxK\Aq$; namely,
$$
 \langle u\ox\varphi, v\ox\psi\rangle  = \langle u,v\rangle \varphi^* \psi \in \Aq,
\tag{4.23}
$$
for $u,v\in M$ and $\varphi, \psi\in \Aq$. 
Furthermore, $\langle\ ,\ \rangle$ can be chosen so that it should be positive definite
when $q$ is specialized to a real number with $|q|\ne 0,1$; 
then one has $\langle u, u\rangle\ne 0$ for any nonzero vector $u\in M$. 
As to the left $\Uq$-module structure of $M$, the $\text{\rm U}_q(N)$-invariance 
of $\langle \ ,\ \rangle $ implies 
$$
                \langle u, a.v\rangle =\langle a^*.u, v\rangle \quad\text{for}\ \ u,v\in M \ \ \text{and}\ \ a\in\Uq, 
\tag{4.24}
$$
under the $\ast$-operation of $\Uq$ (see Section 1.4). 
Note that, if $M$ is irreducible, a $\text{\rm U}_q(N)$-invariant 
hermitian form is determined uniquely up to a scalar multiple. 
\par
        Fixing a hermitian form on $M$ as above, we define the matrix element 
$\phi_M(u,v) \in \Aq$ of  $M$  associated with a pair $(u,v)$ of elements in $M$ by 
$$
                \phi_M(u,v):=\langle u, \rho_G(v)\rangle  \in \Aq. 
\tag{4.25}
$$
By the $\text{\rm U}_q(N)$-invariance of the hermitian form, one can easily show 
\Lemma{4.8}{
(1)  $\phi_M(v,u) = \tau(\phi_M(u,v)) $ for all $u,v\in M$. \newline\noindent
(2)  $a.\phi_M(u,v) = \phi_M(u, a.v)$ and $\phi_M(u,v).a=\phi_M(a^*.u, v)$ for all $u,v\in M$ and $a\in\Uq$. 
}
When $M=V(\ld)\ (\ld\in P^+)$, we will write $\phi_\ld(u,v)=\phi_{V(\ld)}(u,v)$ for short. 
Denoting $V(\ld)^\circ$ the right $\Uq$-module obtained from $V(\ld)$, we can regard 
the hermitian form $\langle \ ,\ \rangle $ as a $\K$-bilinear form $V(\ld)^\circ \times V(\ld)\to\K$. 
Accordingly, the mapping $(u,v)\mapsto\phi_{\ld}(u,v)$ give rise to
a $\Uq$-bimodule homomorphism
$$ 
        \phi_\ld : V(\ld)^\circ\oxK V(\ld) \to \Aq, 
\tag{4.26}
$$
and its image coincide with the subspace $W(\ld)$ of matrix elements of $\Aq$. 
This gives a description of the isomorphism $V(\ld)^\circ\ox V(\ld) @>\sim>>W(\ld)$. 
\par
        We fix a highest weight vector $u(\ld)$ in $V(\ld)$. 
We now assume that $\ld\in P^+_\k$ and 
denote by $w(\ld)$ the $\kq$-fixed vector in $V(\ld)$ 
normalized so that the highest weight component of $w(\ld)$ should give $u(\ld)$; 
this normalization makes sense by Lemma 3.2. 
Recalling that $\kq=\kq(a)$ implies ${\kq}^*=\kq(a^{-1})$, we also take the 
${\kq}^*$-fixed vector $w^*(\ld)$ in $V(\ld)$ such that 
$w^*(\ld)_\ld=u(\ld)$. 
By using these vectors, we define the matrix elements $\varphi_0(\ld)$ and $\varphi(\ld)$ by 
$$
\align
        \varphi_0(\ld)&=\phi_{\ld}(u(\ld), w(\ld))/\langle u(\ld), u(\ld)\rangle, 
\tag{4.27}\\ 
        \varphi(\ld)&=\phi_{\ld}(w^*(\ld), w(\ld))/\langle u(\ld), u(\ld)\rangle   
\endalign
$$
\Lemma{4.9}{
Let  $\ld={\sum}_{k=1}^N \ \ld_k \ep_k$ be an element of $P^+_\k$. \newline\noindent
{\rm (1)} 
The matrix element $\varphi_0(\ld)$ gives a highest weight vector of the right 
$\Uq$-module $W(\ld)_{\kq}$. 
Furthermore, its restriction to the diagonal subgroup $\T$ is given by
$$
        \varphi_0 (\ld)|_\T=z^\ld={z_1}^{\ld_1}\cdots {z_N}^{\ld_N}. 
\tag{4.28}
$$
\newline\noindent
{\rm (2)} The matrix element $\varphi(\ld)$ gives a $\kq$-biinvariant element in $W(\ld)$, 
i.e., $\varphi(\ld)\in\Cal{H}(\ld)$. 
Furthermore, its restriction to the diagonal subgroup $\T$ is a homogeneous polynomial 
in the form 
$$
        \varphi (\ld)|_\T=z^\ld + {\sum}_{\mu<\ld} a_{\ld\mu} z^\mu 
\quad(a_{\ld\mu}\in \K), 
\tag{4.29}
$$
where $<$ denotes the dominance order of weights in $P$. 
}
\Proof{
From Lemma 4.8, it follows directly 
that $\varphi_0(\ld)$ is a highest weight vector 
of $W(\ld)_{\kq}$ and that $\varphi(\ld)$ is a $\kq$-biinvariant element 
in $W(\ld)$. 
Setting $u_0=u(\ld)$, take a basis $\{u_0,u_1,\cdots,u_m\}$ for $V(\ld)$, 
consisting of weight vectors $u_j$ of weight $\mu^{(j)}$, so that $\mu^{(0)}=\ld$. 
Note that $\langle u_0, u_j\rangle=0$ for $j=1,\cdots,m$. 
Write the vectors $w(\ld)$ and $w^*(\ld)$ in the form
$$
\align
	w(\ld) & = u_0 + {\sum}_{j=1}^m c_j u_j, 
\tag{4.30}\\
	w^*(\ld) & = u_0 + {\sum}_{j=1}^m d_j u_j, 
\endalign
$$
where $c_j, d_j\in\K$.  
Note that the restriction of $\varphi_0(\ld)$ and $\varphi(\ld)$ to $\T$ 
can be written as follows by using the coaction 
$\rho_\T : V(\ld)\to V(\ld)\oxK A(\T)$ :
$$
\align
\varphi_0(\ld)|_\T &= \langle u(\ld), \rho_\T(w(\ld))\rangle/\langle u(\ld), u(\ld)\rangle, 
\tag{4.31}\\
\varphi(\ld)|_\T &= \langle w^*(\ld), \rho_\T(w(\ld))\rangle/\langle u(\ld), u(\ld)\rangle. 
\endalign
$$
Since 
$\rho_\T(w(\ld))={\sum}_{j=0}^m c_j u_j\ox z^{\mu^{(j)}}$,
we compute
$$
\langle u(\ld), u(\ld)\rangle\varphi_0(\ld)|_\T = {\sum}_{j=0}^m c_j \langle u(\ld),u_j\rangle  z^{\mu^{(j)}} 
	= \langle u(\ld),u(\ld)\rangle z^\ld. 
\tag{4.32}
$$
Similarly, we have
$$
\align
\langle u(\ld), u(\ld)\rangle\varphi(\ld)|_\T &= {\sum}_{j=0}^m c_j \langle w^*(\ld),u_j\rangle  z^{\mu^{(j)}} 
\tag{4.33}\\
	&= \langle u(\ld),u(\ld)\rangle z^\ld + {\sum}_{j=1}^m c_j \langle w^*(\ld),u_j\rangle  z^{\mu^{(j)}},
\endalign
$$
as desired. 
}
For the proof of Theorem 4.3 and Theorem 4.7, we need to determine the explicit form 
of $\varphi_0(\ld)$ and $\varphi(\ld)$ for the fundamental weights $\ld$ in $P^+_\k$. 
Let $1\le i_1<i_2<\cdots<i_r\le N$  and $1\le j_1<j_2<\cdots<j_r\le N$ be two increasing 
sequence of indices of length $r \ (1\le r\le N)$.  
Then we denote by $\xi^{i_1\cdots i_r}_{j_1\cdots j_r}$ the quantum minor determinant 
of the matrix $T=(t_{ij})_{1\le i,j\le N}$, with row indices $i_1,\cdots,i_r$ and column 
indices $j_1,\cdots,j_r$: 
$$
        \xi^{i_1\cdots i_r}_{j_1\cdots j_r}
        ={\sum}_{w\in\frak{S}_r}\ (-q)^{\ell(w)} 
	t_{i_{w(1)}j_1}t_{i_{w(2)}j_2}\cdots t_{i_{w(r)}j_r}.  
\tag{4.34}
$$
When $(i_1,i_2,\cdots,i_r)=(1,2,\cdots,r)$, we write 
$\xi_{j_1\cdots j_r}=\xi^{1\ \cdots r}_{j_1\cdots j_r}$ for short. 
It is known that the fundamental representation $V(\Ld_r)={\bigwedge}^r_q(V)$ has a $\text{\rm U}_q(N)$-invariant 
hermitian form such that the basis 
$\{v_{j_1}\wedge\cdots\wedge v_{j_r}\}_{1\le j_1<\cdots<j_r\le N}$
is an orthonormal basis (see \cite{NYM}). 
Under this hermitian form, we have
$$
\xi^{i_1\cdots i_r}_{j_1\cdots j_r}=
\phi_{\Ld_r}(v_{i_1}\wedge\cdots\wedge v_{i_r},v_{j_1}\wedge\cdots\wedge v_{j_r}),
\tag{4.35}
$$
for $i_1<\cdots<i_r$ and $j_1<\cdots<j_r$.
\Lemma{4.10.A}{
In Case (SO), the matrix elements $\varphi_0(2\Ld_r)$ and $\varphi(2\Ld_r)$ are determined
as follows:
$$
\align
\varphi_0(2\Ld_r)&= {\sum}_{j_1<\cdots<j_r} 
(\xi_{j_1\cdots j_r})^2 a^{-1}_1\cdots a^{-1}_r a_{j_1}\cdots a_{j_r},
\tag{4.36}\\
\varphi(2\Ld_r)&= {\sum}_{i_1<\cdots<i_r; j_1<\cdots<j_r} 
(\xi^{i_1\cdots i_r}_{j_1\cdots j_r})^2 a^{-1}_{i_1}\cdots a^{-1}_{i_r}a_{j_1}\cdots a_{j_r}.
\endalign
$$
Furthermore, for any $\ell\in\Z$, we have 
$\varphi_0(\ell\Ld_n)=\varphi(\ell\Ld_n)=\detq(T)^\ell.$ 
}
\Proof{
It is directly checked that the $\kq$-fixed vector $w_r$ of Lemma 3.3.A lies in the 
$\Uq$-submodule of ${\bigwedge}^r_q(V)\oxK{\bigwedge}^r_q(V)$, generated by 
the highest weight vector $v_1\wedge\cdots\wedge v_r\ox v_1\wedge\cdots\wedge v_r$.  
In fact, each summand 
$v_{k_1}\wedge\cdots\wedge v_{k_r}\ox v_{k_1}\wedge\cdots\wedge v_{k_r}$ 
of $w_r$ is obtained from $v_1\wedge\cdots\wedge v_r\ox v_1\wedge\cdots\wedge v_r$ 
by applying the elements ${f_j}^2\ (1\le j\le n-1)$ repeatedly. 
Hence, we can compute the matrix elements 
$\varphi_0(2\Ld_r)$ and $\varphi(2\Ld_r)$ by means of the $\kq$-fixed vector in 
${\bigwedge}^r_q(V)\oxK{\bigwedge}^r_q(V)$.  
Setting 
$u(2\Ld_r)=v_1\wedge\cdots\wedge v_r\ox v_1\wedge\cdots\wedge v_r$, we have 
$$
\align
w(2\Ld_r)&=a^{-1}_1\cdots a^{-1}_r w_r 
\tag{4.37}\\
	 &={\sum}_{k_1<\cdots<k_r} 
v_{k_1}\wedge\cdots\wedge v_{k_r}\ox v_{k_1}\wedge\cdots\wedge v_{k_r}
a^{-1}_1\cdots a^{-1}_r a_{k_1}\cdots a_{k_r}, 
\\
w^*(2\Ld_r)&={\sum}_{k_1<\cdots<k_r} 
v_{k_1}\wedge\cdots\wedge v_{k_r}\ox v_{k_1}\wedge\cdots\wedge v_{k_r}
a_1\cdots a_r a^{-1}_{k_1}\cdots a^{-1}_{k_r}. 
\endalign
$$
This implies 
$$
\align
a_1\cdots a_r\rho_G(w(2\Ld_r)) &= \rho_G(w_r)
\tag{4.38}\\
&={\sum}_{i_1<\cdots<i_r;j_1<\cdots<j_r} 
v_{i_1}\wedge\cdots\wedge v_{i_r}\ox v_{j_1}\wedge\cdots\wedge v_{j_r}\\
&\quad\quad\ox 
{\sum}_{k_1<\cdots<k_r}
\xi^{i_1\cdots i_r}_{k_1\cdots k_r}
\xi^{j_1\cdots j_r}_{k_1\cdots k_r} 
a_{k_1}\cdots a_{k_r}.
\endalign
$$
Hence, under the induced hermitian form on 
${\bigwedge}^r_q(V)\oxK{\bigwedge}^r_q(V)$, we compute 
$$
\phi_{2\Ld_r}(u(2\Ld_r),w(2\Ld_r))=
{\sum}_{k_1<\cdots<k_r} (\xi_{k_1\cdots k_r})^2
a^{-1}_1\cdots a^{-1}_r a_{k_1}\cdots a_{k_r},
\tag{4.39}
$$
and
$$
\align
&\phi_{2\Ld_r}(w^*(2\Ld_r),w(2\Ld_r))
\tag{4.40}\\
&\quad\quad={\sum}_{i_1<\cdots<i_r; j_1<\cdots<j_r} 
(\xi^{i_1\cdots i_r}_{j_1\cdots j_r})^2 a^{-1}_{i_1}\cdots a^{-1}_{i_r}a_{j_1}\cdots a_{j_r}.
\endalign
$$
The last statement of Lemma is clear since $V(\ell\Ld_n)$ is one dimensional and its matrix 
element is given by $\detq(T)^\ell$. 
}
\Lemma{4.10.B}{
In Case (Sp), the matrix elements $\varphi_0(\Ld_{2r})$ and $\varphi(\Ld_{2r})$ are determined
as follows:
$$
\align
  \varphi_0(\Ld_{2r})&= {\sum}_{1\le j_1<\cdots<j_r\le n}  
  \xi_{2j_1-1,2j_1,\cdots, 2j_r-1,2j_r} a^{-1}_1\cdots a^{-1}_r a_{j_1}\cdots a_{j_r},
\tag{4.41}\\
  \varphi(\Ld_{2r})&= {\sum}_{1\le i_1<\cdots<i_r\le n; 1\le j_1<\cdots<j_r\le n} 
\xi^{2i_1-1,2i_1,\cdots,2i_r-1, 2i_r}_{2j_1-1,2j_1,\cdots,2j_r-1, 2j_r} a^{-1}_{i_1}\cdots a^{-1}_{i_r}a_{j_1}\cdots a_{j_r}.
\endalign
$$
Furthermore, for any $\ell\in\Z$, we have 
$\varphi_0(\ell\Ld_{2n})=\varphi(\ell\Ld_{2n})=\detq(T)^\ell.$ 
}
\Proof{
As to ${\bigwedge}^{2r}_q(V)$, we can take
$$
\align
	u(\Ld_{2r}) &= v_1\wedge v_2\wedge\cdots\wedge v_{2r},
\tag{4.42}\\
	w(\Ld_{2r}) &= a^{-1}_1\cdots a^{-1}_r w_r\\
	&=\sum_{1\le k_1<\cdots<k_r\le n} 
	v_{2k_1-1}\wedge v_{2k_1}\wedge\cdots\wedge v_{2k_r-1}\wedge v_{2k_r}
	a^{-1}_1\cdots a^{-1}_r a_{k_1}\cdots a_{k_r},\\
	w^*(\Ld_{2r}) &=
	\sum_{1\le k_1<\cdots<k_r\le n} 
	v_{2k_1-1}\wedge v_{2k_1}\wedge\cdots\wedge v_{2k_r-1}\wedge v_{2k_r}
	a_1\cdots a_r a^{-1}_{k_1}\cdots a^{-1}_{k_r}.
\endalign
$$
Then we have
$$
\align
	a_1\cdots a_r \rho_G(w(\Ld_{2r}))&= \rho_G(w_r) \tag{4.43}\\
	&={\sum}_{i_1<i_2<\cdots<i_{2r}} 
	v_{i_1}\wedge v_{i_2}\wedge\cdots\wedge v_{i_{2r}}\\
	&\quad\quad\ox {\sum}_{1\le k_1<\cdots<k_r\le n} 
	\xi^{i_1 \ i_2 \ \cdots\cdots \ i_{2r}}_{2k_1-1,2k_1,\cdots,2k_r-1,2k_r}
	a_{k_1}\cdots a_{k_r}.
\endalign
$$
With the vectors in (4.42), we can easily see that the elements 
$\phi_{\Ld_{2r}}(u(\Ld_{2r}),w(\Ld_{2r}))$
and
$\phi_{\Ld_{2r}}(w^*(\Ld_{2r}),w(\Ld_{2r}))$
are written in the form (4.41). 
}
%
%
\subsection{4.4. Proof of Theorem 4.3}
%
We now prove that the quadratic elements $x_{ij}\ (1\le i,j\le N)$ and 
$\detq(T)^{\pm 1}$
generate the invariant ring $A_q(G/K)$. 
For this purpose, we look at the highest weight vectors in $A_q(G/K)$. 
Consider the subalgebra
$$
        A_q(N\backslash G/K)=\{\varphi\in A_q(G/K); \varphi.\n_q =0 \}
\tag{4.44}
$$
of highest weight vectors in $A_q(G/K)$, 
where $\n_q$ denotes the coideal of $\Uq$ spanned by $L^+_{ij}\ (i<j)$,
corresponding to the nilpotent Lie subalgebra of lower triangular matrices. 
By the irreducible decomposition of Proposition 4.1, we see that the matrix 
elements $\varphi_0(\ld)\ (\ld\in P^+_\k)$ form a $\K$-basis for 
$A_q(N\backslash G/K)$:
$$
        A_q(N\backslash G/K)={\bigoplus}_{\ld\in P^+_\k} \K \varphi_0(\ld).  
\tag{4.45} 
$$
Then by Lemma 4.9.(1), we see that this algebra is isomorphic to the following 
subalgebra of $A(\T)=\K[z^{\pm 1}]$, through the restriction mapping 
$A_q(N\backslash G/K) \to A(\T)$:
$$
        A_q(N\backslash G/K)|_\T ={\bigoplus}_{\ld\in P^+_\k} \K z^\ld. 
\tag{4.46} 
$$
This implies:
$$
\align
  \text{Case (SO):}&\quad
  A_q(N\backslash G/K)|_\T = 
  \K[z^{2\Ld_1},z^{2\Ld_2},\cdots,z^{2\Ld_{n-1}}, z^{\pm\Ld_n}], 
\tag{4.47}\\
  \text{Case (Sp):}&\quad
  A_q(N\backslash G/K)|_\T = 
  \K[z^{\Ld_{2}},z^{\Ld_{4}},\cdots,z^{\Ld_{2(n-1)}}, z^{\pm\Ld_{2n}}].  
\endalign
$$
Hence we have
\Lemma{4.11}{
The algebra $A_q(N\backslash G/K)$ is a commutative $\K$-algebra generated by the following matrix elements:
$$
\align
\text{Case (SO):} &\quad \varphi_0(2\Ld_r) \  (1\le r\le n-1),\ \  \detq(T)^{\pm 1},
\tag{4.48}\\
\text{Case (Sp):} &\quad \varphi_0(\Ld_{2r}) \  (1\le r\le n-1),\ \  \detq(T)^{\pm 1}.
\endalign
$$
}
Note that the algebra $A_q(G/K)$ is generated by the subalgebra 
$A_q(N\backslash G/K)$ as a right $\Uq$-module. 
Hence, 
in order to prove that the quadratic elements $x_{ij}\ (1\le i,j\le N)$ and $\detq(T)^{\pm 1}$ generate 
the algebra $A_q(G/K)$, we have only to show that the generators of 
$A_q(N\backslash G/K)$ above are actually contained in the algebra 
$\K[x_{ij} \ (1\le i,j\le N), \detq(T)^{\pm}]$. 
It should be noted that the quadratic elements $x_{ij}$ arise from 
the coaction of $\Aq$ at the $\kq$-invariant element 
$w_J\in V\oxK V $ of Proposition 2.3:
$$
\align
\rho_G(w_J)
&={\sum}_{1\le i,j\le N} v_i\ox v_j
\ox {\sum}_{1\le k,\ell\le N}t_{ik}J_{k\ell}t_{j \ell} \tag{4.49}\\
&={\sum}_{1\le i,j\le N} v_i\ox v_j\ox x_{ij}
\endalign
$$
\par
\medpagebreak\par\noindent
Case (SO):  We show that the matrix elements $\varphi_0(2\Ld_r)\ (1\le r \le n)$ 
are contained in the subalgebra $\K[x_{ij}\ (1\le i,j\le n)]$ of $A_q(G/K)$. 
We make use of the intertwining operators 
$
\Phi : (V\oxK V)^{\ox r} @>\sim>> V^{\ox r}\oxK V^{\ox r}
$
and 
$
\Psi : (V\oxK V)^{\ox r} \to {\bigwedge}^r_q(V)\oxK {\bigwedge}^r_q(V)
$
as in the proof of Lemma 3.3.A. 
We take the following matrix representation of $\Phi$:
$$
\Phi(\bv_1\ox\bv_{1'}\ox\cdots\ox\bv_r\ox\bv_{r'})
=\bv_1\ox\cdots\ox\bv_r\ox\bv_{1'}\ox\cdots\ox\bv_{r'} 
\Phi_{1\cdots r;1'\cdots r'},
\tag{4.50}
$$
where
$$
\Phi_{1\cdots r;1'\cdots r'}=R^+_{r1'}\cdots R^+_{r(r-1)'}\cdots 
R^+_{31'} R^+_{32'} R^+_{21'}. 
\tag{4.51}
$$
Then we have 
$$
\align
&\Psi(v_{i_1}\ox v_{j_1}\ox\cdots\ox v_{i_r}\ox v_{j_r})
\tag{4.52} \\
&\quad={\sum}_{\mu_1<\cdots<\mu_r;\nu_1<\cdots<\nu_r} 
v_{\mu_1}\wedge\cdots\wedge v_{\mu_r}\ox v_{\nu_1}\wedge\cdots\wedge v_{\nu_r} 
\Psi^{\mu_1\cdots\mu_r;\nu_1\cdots\nu_r}_{i_1 \cdots \ i_r ;j_1 \cdots j_r}, 
\endalign
$$
where the matrix coefficients of $\Psi$ are given by 
$$
\Psi^{\mu_1\cdots\mu_r;\nu_1\cdots\nu_r}_{i_1 \cdots i_r ;j_1 \cdots  j_r}
={\sum}_{\sigma,\tau\in\frak{S}_r} (-q)^{\ell(\sigma)+\ell(\tau)}
\Phi^{\mu_{\sigma(1)}\cdots\mu_{\sigma(r)};\nu_{\tau(1)}\cdots\nu_{\tau(r)}}
_{i_1\ \cdots\cdots\ i_r\ ; \ j_1\ \cdots\cdots\ j_r}. 
\tag{4.53}
$$
Since $\Psi$ is a homomorphism of right $\Aq$-comodules, 
the equality $ \Psi({w_J}^{\ox r}) = [r]_{q^2}! w_r $ in (3.26) implies
$$
\Psi(\rho_G(w_J^{\ox r})) = [r]_{q^2}! \rho_G(w_r).  
\tag{4.54}
$$
As the right-hand side is already given in (4.38), 
we now compute the left-hand side of (4.54) explicitly. 
From $\rho_G(w_J)={\sum}_{i,j} v_i\ox v_j \ox x_{ij}$, we have 
$$
\rho_G({w_J}^{\ox r})=
{\sum}_{i_1,\cdots i_r;j_1,\cdots,j_r} 
v_{i_1}\ox v_{j_1}\ox\cdots\ox v_{i_r}\ox v_{j_r}
\ox x_{i_1j_1}\cdots x_{i_rj_r}. 
\tag{4.55}
$$
Hence we have 
$$
\align
\Psi(\rho_G(w_J^{\ox r}))
&={\sum}_{\mu_1<\cdots<\mu_r; \nu_1<\cdots<\nu_r}
v_{\mu_1}\wedge\cdots\wedge v_{\mu_r}\ox
v_{\nu_1}\wedge\cdots\wedge v_{\nu_r} 
\tag{4.56}\\
&\quad\ox {\sum}_{i_1,\cdots,i_r;j_1,\cdots,j_r}
\Psi^{\mu_1\cdots \mu_r;\nu_1\cdots \nu_r}_{i_1 \cdots i_r;j_1 \cdots j_r}
x_{i_1j_1}\cdots x_{i_rj_r}. 
\endalign
$$
by (4.52). 
Comparing this formula with (4.38), we obtain
$$
\align
&[r]_{q^2}! {\sum}_{k_1<\cdots<k_r} (\xi_{k_1\cdots k_r})^2 a_{k_1}\cdots a_{k_r}
\tag{4.57}\\
&\quad ={\sum}_{i_1,\cdots,i_r:j_1,\cdots,j_r} 
\Psi^{1 \  \cdots r ;1 \  \cdots \ r}
_{i_1\cdots i_r;j_1\cdots j_r}
x_{i_1j_1}\cdots x_{i_rj_r},
\endalign
$$
which gives an expression of $\varphi_0(2\Ld_r)$ in terms of 
$x_{ij}\ (1\le i,j\le n)$. 
This completes the proof of Theorem 4.3 for Case (SO). 
\par
\medpagebreak\par\noindent
Case (Sp): We show that the matrix elements $\varphi_0(\Ld_{2r})\ (1\le r \le n)$ 
are contained in the subalgebra $\K[x_{ij}\ (1\le i,j\le 2n)]$ of $A_q(G/K)$. 
With the notation of Lemma 3.3.B, we consider the $\kq$-fixed vector $w_r$ in 
${\bigwedge}^r_q(V)$.
From (4.49), we obtain
$$
	\rho_G(w_1)={\sum}_{1\le i<j\le 2n} v_i\wedge v_j\ox x_{ij}. 
\tag{4.58}
$$
Since the coaction $\rho_G: {\bigwedge}_q(V) \to {\bigwedge}_q(V)\oxK\Aq$ 
is an algebra homomorphism, we have 
$$
\rho_G({w_1}^{\wedge r})
\quad={\sum}_{i_1<j_1;\cdots;i_r<j_r} 
v_{i_1}\wedge v_{j_1}\wedge\cdots\wedge
v_{i_r}\wedge v_{j_r}\ox x_{i_1j_1}\cdots x_{i_rj_r}. 
\tag{4.59}
$$
On the other hand, we know that $w_1^{\wedge r}=[r]_{q^4}! w_r$; hence 
equality (4.59) gives a formula for $ [r]_{q^4}! \rho_G(w_r)$.
Comparing this with (4.43), we have
$$
\align
&[r]_{q^4}!{\sum}_{k_1<\cdots<k_r} \xi_{2k_1-1,2k_1,\cdots,2k_r-1,2k_r}
a_{k_1}\cdots a_{k_r} \tag{4.60}\\
&\quad
={\sum}_{w\in\frak{S}_{2r}: w(1)<w(2);\cdots;w(2r-1)<w(2r)}
(-q)^{\ell(w)} x_{w(1)w(2)}\cdots x_{w(2r-1)w(2r)}.  
\endalign
$$
This gives the expression of $\varphi_0(\Ld_{2r})$ in terms of 
$x_{ij}\ (1\le i,j\le 2n)$. 
This completes the proof of Theorem 4.3 for Case (Sp). 
\Remark{4.12}{
In Case (Sp), (4.60) contains the following formula for the quantum 
determinant $\detq(T)$:
$$
\align
&[n]_{q^4}! \detq(T) a_1\cdots a_n \tag{4.61}\\
&\quad={\sum}_{w\in\frak{S}_{2n}: w(1)<w(2);\cdots;w(2n-1)<w(2n)}
(-q)^{\ell(w)} x_{w(1)w(2)}\cdots x_{w(2n-1)w(2n)}, 
\endalign
$$
which can be regarded as the ``quantum Pfaffian'' of the $q$-skew-symmetric 
matrix $X=TJT^t$.
In Case (SO), (4.57) implies a formula representing the square 
$\detq(T)^2$ of the quantum determinant in terms of the matrix elements 
of $X=TJT^t$.  
We do not know, however, whether (4.57) reduces to a simple formula, even in the 
case when $r=n$. 
}
%
%
\subsection{4.5. Proof of Theorem 4.7}
%
By Lemma 4.9.(2), we see that the zonal spherical functions 
$\varphi(\ld)$ ($\ld\in P^+_{\k}$) form a $\K$-basis of the subalgebra
${\Cal H}=A_q(K\backslash G/K)$ of $\kq$-biinvariant elements.  
Furthermore, the restriction of $\varphi(\ld)$ to the diagonal subgroup 
$\T$ has the leading term $z^{\ld}$ under the lexicographic order 
of the monomials in $A(\T)=\K[z^{\pm 1}]$. 
This shows that the restriction mapping 
${\Cal H}=A_q(K\backslash G/K) @>>> A(\T) $
is injective. 
\par
From Lemma 4.10.A, it follows that 
$$
  \varphi(2\Ld_r)|_\T = e_r(z_1^2,\cdots,z_n^2)\quad(1\le r\le n)
  \tag{4.62.a}
$$
in Case (SO).
Here we denoted by $e_r(x_1,\cdots,x_n)$ the elementary symmetric
function of degree $r$ in the variables $(x_1,\cdots,x_n)$. 
Similarly, from Lemma 4.10.B it follows that 
$$
  \varphi(\Ld_{2r})|_\T = e_r(z_1z_2,\cdots,z_{2n-1}z_{2n})\quad(1\le r\le n), 
  \tag{4.62.b}
$$
in Case (Sp). 
Note also that 
$$
  \varphi(\ell\Ld_N)|_\T = \detq(T)^\ell|_\T = (z_1z_2\cdots z_N)^\ell\quad(\ell\in\Z), 
\tag{4.63}
$$
in each case. 
Hence the statement concerning the image ${\Cal H}|_\T$ in Theorem 4.7 is equivalent
to the following lemma.
\Lemma{4.13}{ The algebra ${\Cal H}=A_q(K\backslash G/K)$ of $\kq$-biinvariants 
is generated by the following matrix elements:
$$
\align
\text{Case (SO):}& \quad\varphi(2\Ld_r) \ (1\le r\le n-1), \ \detq(T)^{\pm 1}
\tag{4.64}\\
\text{Case (Sp):}& \quad\varphi(\Ld_{2r}) \ (1\le r\le n-1), \ \detq(T)^{\pm 1}
\endalign
$$
}
\Proof{
Recall that the ``coordinate ring'' $A_q(\Mat(N))$ of the quantum matrix space $\Mat_q(N)$ 
has the following irreducible decomposition as a $\Uq$-bimodule:
$$
	A_q(\Mat(N))={\bigoplus}_{\ld\in P^+ \cap L} W(\ld),
\tag{4.65}
$$
where $L=\{\ld\in P ; \langle \ld,\ep_k\rangle\ge 0 \ (1\le r\le N)\}$ denotes the first 
quadrant of the weight lattice $P$. 
Let us denote by ${\Cal H}_{\ge 0}=A_q(K\backslash \Mat(N)/K)$ the 
the subalgebra of all $\kq$-biinvariants in $A_q(\Mat(N))$. 
Then from the decomposition (4.65) we have 
$$
 {\Cal H}_{\ge 0} = {\bigoplus}_{\ld\in P^+_\k \cap L} \K \varphi(\ld) 
\tag{4.66}
$$
just as in the case of ${\Cal H}$. 
It is clear that 
$$
\varphi(\ld) \detq(T)^\ell = \varphi(\ld+\ell\Ld_N) \quad\text{for any} \ \ 
\ld\in P^+_\k,  \ \ \ell\in\Z.
\tag{4.67}
$$
Hence we have only to prove that ${\Cal H}_{\ge 0}$ is generated by the following 
matrix elements:
$$
\align
\text{Case (SO):}& \quad\varphi(2\Ld_r)  \ (1\le r\le n-1), \ \detq(T), 
\tag{4.68}\\
\text{Case (Sp):}& \quad\varphi(\Ld_{2r})  \ (1\le r \le n).
\endalign
$$
Recall that any $\ld\in P^+_\k \cap L$ can be written in the form 
$$
\align
\text{Case (SO):}& \quad\ld={\sum}_{r=1}^{n-1} 2m_r\Ld_r +\ell\Ld_n 
\tag{4.69}\\
\text{Case (Sp):}& \quad\ld={\sum}_{r=1}^{n} m_r\Ld_{2r} 
\endalign
$$
for some $m_r\in\N \ (1\le r\le n)$ and $\ell\in\N$. 
In view of (4.69), we define the element $\tilde{\varphi}(\ld)$ 
in ${\Cal H}_{\ge 0}$ by 
$$
\align
\text{Case (SO):}& \quad\widetilde{\varphi}(\ld)=
\varphi(2\Ld_1)^{m_1}\cdots\varphi(2\Ld_{n-1})^{m_{n-1}}\varphi(\Ld_n)^{\ell}
\tag{4.70}\\
\text{Case (Sp):}& \quad\widetilde{\varphi}(\ld)=
\varphi(\Ld_2)^{m_1}\cdots\varphi(\Ld_{2n})^{m_n}.
\endalign
$$
Note that, 
under the lexicographic order $\preceq$ of $L$, 
the polynomials 
$\varphi(\ld)|_\T$ and $\widetilde{\varphi}(\ld)|_\T$ in $\K[z]$ 
have the common leading term $z^\ld$.  
From this fact, one can easily show that ${\Cal H}_{\ge 0}$ is generated
by the elements in (4.68), by the induction with respect to the 
well-ordering $\preceq$ of $L$. 
If $\varphi$ is a nonzero elements in ${\Cal H}_{\ge 0}$, its 
restriction $\varphi|_\T$ has the leading term $cz^\mu$ for some 
$\mu\in P^+_\k \cap L$ and $c\in\K$ by (4.66). 
By the induction hypothesis,
the element $\psi=\varphi-c\widetilde{\varphi}(\mu)$ is 
expressed as a linear combination of elements in (4.70), 
since $\psi$ has the leading exponent strictly less than $\mu$ 
under $\preceq$. 
Hence $\varphi=\psi+c\widetilde{\varphi}(\mu)$ also lies 
in the subalgebra of ${\Cal H}_{\ge 0}$ generated by the elements 
in (4.68). 
}
This completes the proof of Theorem 4.7. 
\Remark{4.14}{
As we have seen in the proof of Lemma 4.13, we have
$$
{\bigoplus}_{\ld\in P^+_\k \cap L} \K \varphi(\ld)|_\T
=\K[z_1z_2, \cdots,z_{2n-1}z_{2n}]^{\frak{S}_n}
\tag{4.71.b}
$$
in Case (Sp). 
In Case (SO), we have 
$$
{\bigoplus}_{\ld\in P^+_\k \cap 2L} \K \varphi(\ld)|_\T
=\K[z^2_1 \cdots,z^2_n]^{\frak{S}_n}. 
\tag{4.71.a}
$$
To show (4.71.a), we have to prove that, if $\ld\in P^+_\k \cap 2L$, then 
$\varphi(\ld)|_\T$ belongs to the algebra $\K[z^2_1 \cdots,z^2_n]$. 
If $\ld\in P^+_\k \cap 2L$, we can set $\ell=2m_n$ in (4.69) and we have
$$
\varphi_0(\ld)= \varphi_0(2\Ld_1)^{m_1}\cdots\varphi_0(2\Ld_n)^{m_n}
\tag{4.72}
$$
By Lemma 4.10.A, it is clear that, in the weight decomposition of $\varphi_0(2\Ld_r)$,  
nonzero components occurs only for the weights in $2 L$, for any $1\le r\le n$. 
Accordingly, $\varphi_0(\ld)$ also has the same property by (4.72). 
In view of (4.33), we conclude that $\varphi(\ld)|_\T$ lies in $\K[z^2_1 \cdots,z^2_n]$
as desired. 
}
%
\section{\S5. Macdonald's symmetric polynomials as zonal spherical functions}
%
%
In this section we investigate the restriction of the zonal spherical 
functions $\varphi(\ld)\ (\ld\in P^+_\k)$ to the diagonal subgroup $\T$ 
of $\GL_q(N)$.
They are expressed by Macdonald's symmetric polynomials 
$P_\mu(x;q,t)$ in $n$ variables with a special value of $(q,t)$. 
This result will be established by computing the radial component of a central 
element of $\Uq$. 
%
%
\subsection{5.1. Macdonald's symmetric polynomials}
%
%
We begin with a recall on Macdonald's symmetric polynomials
(\cite{M2, M3}). 
{\it Macdonald's symmetric polynomials\/} $P_\mu(x;q,t)=P_\mu(x_1,\cdots,x_n;q,t)$ 
are a family of symmetric polynomials in $\Q(q,t)[x_1,\cdots,x_n]$, 
homogeneous of degree $\sum_{k=1}^n \mu_k$, parametrized by partitions 
$\mu=(\mu_1,\cdots,\mu_n)\in \N^n\ (\mu_1\ge\cdots\ge\mu_n\ge 0)$. 
Among many characterizations of $P_\mu(x;q,t)$, we recall now the 
following two properties. 
\par\noindent
i) For each $\mu$, the polynomial $P_{\mu}(x;q,t)$ has an expression
$$
P_{\mu}(x;q,t)=m_{\mu}(x)+\sum_{\nu<\mu} c_{\mu\nu} m_{\nu}(x) 
\quad( c_{\mu\nu}\in \Q(q,t) ), 
\tag{5.1}
$$
where $m_{\mu}(x)$ stands for the monomial symmetric function of 
monomial type $\mu$ and the symbol $<$ denotes the dominance order of 
partitions. 
\par\noindent
ii) For each $\mu$, $P_{\mu}(x;q,t)$ satisfies the $q$-difference equation 
$$ 
\sum_{k=1}^n 
\frac{\Dt(x_1,\cdots,tx_k,\cdots,x_n)}{\Dt(x_1,\cdots,x_n)}
T_{q,x_k}\,P_{\mu}(x;q,t)=(\sum_{k=1}^n t^{n-k}q^{\mu_k})P_{\mu}(x;q,t),
\tag{5.2} 
$$
where $\Dt(x_1,\cdots,x_n)$ stands for the difference product 
$$
\Dt(x_1,\cdots,x_n)=\prod_{1\le i<j\le n} (x_j-x_i) 
\tag{5.3}
$$
and $\Dt(x_1,\cdots,tx_k,\cdots,x_n)$ for $\Dt(x_1,\cdots,x_n)$ 
with $x_k$ replaced by $tx_k$. 
The symbol $T_{q,x_k}$ denotes the $q$-shift operator in $x_k$ defined by 
$$
(T_{q,x_k}f)(x_1,\cdots,x_n)=f(x_1,\cdots,qx_k,\cdots,x_n).
\tag{5.4}
$$
\par
We remark that property i) implies that the symmetric polynomials $P_\mu(x;q,t)$ 
form a $\Q(q,t)$-basis for the algebra $\Q(q,t)[x]^{\frak{S}_n}$, as $\mu$ runs 
over the set of all partitions in $\N^n$. 
This means that the $q$-difference operator in the left-hand side of (5.2) 
is diagonalizable on $\Q(q,t)[x]^{\frak{S}_n}$ and that its eigenspaces 
are all one-dimensional. 
Note that the coefficients of the $q$-difference operator are also written in the form 
$$
\frac{\Dt(x_1,\cdots,tx_k,\cdots,x_n)}{\Dt(x_1,\cdots,x_n)}
=\prod_{1\le j\le n ; j\ne k}\frac{tx_k-x_j}{x_k-x_j}.
%
\tag{5.5}
$$
\par
Returning to the setting of Section 4, we consider the zonal spherical function 
$\varphi(\ld)$ associated with the representation $V(\ld)\ (\ld\in P^+_\k)$. 
For the description of $\varphi(\ld)|_\T$, we use the following 
parametrization of $\ld$ by partitions $\mu$:
$$
\align
\text{Case (SO):}&\quad \ld=\sum_{k=1}^n 2\mu_k\ep_k + \ell\Ld_n, 
\tag{5.6} \\
\text{Case (Sp):}&\quad \ld=\sum_{k=1}^n \mu_k(\ep_{2k-1}+\ep_{2k}) +\ell\Ld_{2n},
\endalign
$$
where $\mu=(\mu_1,\cdots,\mu_n)$ stands for a partition in $\N^n$ and $\ell\in\Z$. 
\Theorem{5.1}{
For each $\ld\in P^+_\k$, the restriction of the zonal spherical function 
$\varphi(\ld)$ to the diagonal subgroup $\T$ is expressed in terms of Macdonald's 
symmetric polynomial $P_\mu$. 
To be more precise, we have 
$$
\align
\text{Case (SO):}&\quad 
\varphi(\ld)|_\T = P_\mu(z^2_1,\cdots,z^2_n;q^4,q^2) (z_1\cdots z_n)^\ell,
\tag{5.7}\\
\text{Case (Sp):}&\quad 
\varphi(\ld)|_\T = P_\mu(z_1z_2,\cdots,z_{2n-1}z_{2n};q^2,q^4) 
(z_1z_2\cdots z_{2n})^\ell, 
\endalign
$$
under the parametrization of $\ld$ as in (5.6). 
}
For the identification of $\varphi(\ld)|_\T$ with Macdonald's symmetric polynomial 
$P_\mu$, we will show that $\varphi(\ld)|_\T$ satisfies a $q$-difference equation 
corresponding to (5.2).  
Such a $q$-difference equation arises as the radial component of a central 
element of $\Uq$. 
%
%
\subsection{5.2. Radial component of a central element in $\Uq$}
%
%
Let $C$ be a central element of $\Uq$. 
Then its action on $\Aq$ preserves the subalgebra ${\Cal H}=A_q(K\backslash G/K)$.
Hence, the action of $C$ on ${\Cal H}$ induces a $\K$-linear operator 
on its image ${\Cal H}|_\T$ by the restriction ${\Cal H}\to A(\T)$. 
This operator acting on ${\Cal H}|_\T$ will be called the {\it radial component} 
of $C$ and denoted by $C|_\T$. 
Note that, since $C$ is central, the two actions of $C$ on $\Aq$, 
one from the left and the other from the right, eventually coincide. 
In what follows, we define the elements $x_1,\cdots,x_n$ in 
$A(T)=\K[z^{\pm 1}]$ by 
$$
\align 
\text{Case (SO):}&\quad x_1=z^2_1, \cdots, x_n=z^2_n\ \ ,
\tag{5.8}\\
\text{Case (Sp):}&\quad x_1=z_{1}z_{2}, \cdots, x_n=z_{2n-1}z_{2n}. 
\endalign
$$
With these elements, the image ${\Cal H}|_\T$ is written as 
$$
\align 
\text{Case (SO):}&\quad 
{\Cal H}|_\T=\K[x_1,\cdots,x_n]^{\frak{S}_n}[(z_1\cdots z_n)^{-1}],
\tag{5.9}\\
\text{Case (Sp):}&\quad 
{\Cal H}|_\T=\K[x_1,\cdots,x_n]^{\frak{S}_n}[(x_1\cdots x_n)^{-1}].
\endalign
$$
\par
We now recall on the central elements $C_r\ (r=1,2,\cdots)$ of $\Uq$ 
proposed by \cite{RTF}.  
They are defined as 
$$
  C_r = \text{\rm tr}_q((L^+ S(L^-))^r), 
\tag{5.10}
$$
where we use the notation of $q$-trace 
$$
  \text{\rm tr}_q(A)=\sum_{k=1}^N q^{2(N-k)} a_{kk},
\tag{5.11}
$$
for a matrix $A=(a_{ij})_{1\le i,j\le N}$ in $\End(V)\oxK\Uq$. 
Note that the central element $C_1$ takes the form
$$
  C_1=\sum_{1\le i,j\le N} q^{2(N-i)} L^+_{ij}S(L^-_{ji}).
\tag{5.12}
$$
Its eigenvalue on the irreducible representation $V(\ld)$ is given by
$$
\chi_\ld(C_1) = \sum_{k=1}^N q^{2 \langle \ep_k, \ld + \rho \rangle}
= \sum_{k=1}^N q^{2(\ld_k+N-k)}, 
\tag{5.13}
$$
where $\rho=\sum_{k=1}^N (N-k)\ep_k$.
\Theorem{5.2}{
On the subalgebra $\K[x_1,\cdots,x_n]^{\frak{S}_n}$ of ${\Cal H}|_\T$, 
the radial component of the central element $C_1\in\Uq$ above is 
given by the following $q$-difference operator $D_1$ :
$$
\align 
\text{Case (SO):}&\quad D_1=\sum_{k=1}^n 
\frac{\Dt(x_1,\cdots,q^2x_k,\cdots,x_n)}
{\Dt(x_1,\cdots,x_n)}
T_{q^4,x_k},
\tag{5.14}\\
\text{Case (Sp):}&\quad D_1=(1+q^2) \sum_{k=1}^n 
\frac{\Dt(x_1,\cdots,q^4x_k,\cdots,x_n)}
{\Dt(x_1,\cdots,x_n)} T_{q^2,x_k}. 
\endalign
$$
}
It is easy to show that Theorem 5.2 implies Theorem 5.1. 
From (4.67), it is clear that 
$$
	\varphi(\ld+\ell\Ld_N)|_\T=\varphi(\ld)|_\T (z_1\cdots z_N)^\ell
\tag{5.15}
$$
for any $\ell\in\Z$. 
Hence, we may assume that $\ell=0$ in the parametrization of (5.6). 
When $\ell=0$, the restriction $\varphi(\ld)|_\T$ lies in the algebra $\K[x_1,\cdots,x_n]$ 
as we mentioned in Remark 4.14.  
Note also that, 
under the parametrization of (5.6),
the eigenvalue $\chi_\ld(C_1)$ is written as 
$$
\align
\text{Case (SO):}&\quad \chi_\ld(C_1) = \sum_{k=1}^n q^{2(n-k)} q^{4\mu_k},
\tag{5.16}\\
\text{Case (Sp):}&\quad \chi_\ld(C_1) = (1+q^2)\sum_{k=1}^n q^{4(n-k)} q^{2\mu_k}. 
\endalign
$$
Since  $\varphi(\ld)$  satisfies the equation 
$(C_1 - \chi_\ld(C_1))\varphi(\ld)=0$, 
its restriction $\varphi(\ld)|_\T$ satisfies the $q$-difference equation 
$$
(D_1 - \chi_\ld(C_1))\varphi(\ld)|_\T =0. 
\tag{5.17}
$$
By Theorem 5.2 and (5.16), we see that (5.17) gives rise to 
the $q$-difference equation (5.2) with $(q,t)$ replaced by $(q^4, q^2)$ in Case (SO), 
and by $(q^2, q^4)$ in Case (Sp). 
This $q$-difference equation determines $\varphi(\ld)$ up to a scalar multiple 
since the mapping $\mu\mapsto\chi_\ld(C_1)$ is still injective 
after the specialization of $(q,t)$. 
Noting that $\varphi(\ld)|_\T$ and the corresponding $P_\mu$ has the common 
leading term $z^\ld=x^\mu$, we obtain Theorem 5.1 for the case where $\ell=0$. 
\par
We remark that this argument is also valid in the setting where 
$\K=\C$ and $q$ is a real number with $|q|\neq 0,1$, as the partitions 
$\mu$ are separated by the values of $\chi_\ld(C_1)$. 

%
%
\subsection{5.3. How to compute the radial component $C|_\T$}
%
%
Before the proof of Theorem 5.2, we explain how we are going to compute 
the radial component $C|_\T$ of a central element $C$ of $\Uq$. 
Our method is based on the duality between $\Aq$ and $\Uq$;
its spirit is the same as that of Koornwinder \cite{K2}.
\par
Recall that there exists a pairing of Hopf algebras 
$(\ ,\ ) : \Uq\times\Aq \to \K$ between $\Aq$ and $\Uq$
(see Section 1.3). 
This pairing induces a $\K$-linear mapping 
$\varphi\mapsto (\cdot,\varphi)$ from 
$\Aq$ to $\Uq^\vee = \Hom(\Uq,\K)$.  
As to the $\Uq$-bimodule structure of $\Aq$, we can easily show 
by the definition (1.18) that 
$$
	(c,a.\varphi.b) = (bca, \varphi)
\tag{5.18}
$$
for all $a,b,c\in\Uq$ and $\varphi\in\Aq$. 
This means that the natural mapping $\Aq\to\Uq^\vee$ is a 
homomorphism of $\Uq$-bimodules. 
Hence, we see that $\Aq\to\Uq^\vee$ is actually injective. 
If its kernel is nontrivial, the irreducible decomposition (1.23) 
implies that the kernel contains some irreducible component $W(\ld)$;
this leads to a contradiction since the pairing $\Uq\times W(\ld)\to\K$ 
cannot be trivial for any $\ld\in P^+$. 
\par
We now consider the subalgebra ${\Cal H}=A_q(K\backslash G/K)$ of $\Aq$. 
From (5.18), it is immediately seen that an element $\varphi\in\Aq$ 
is $\kq$-biinvariant if and only if it satisfies
$$
	(\Uq\kq,\varphi)=0 \ \ \text{and}\ \ (\kq\Uq,\varphi)=0. 
\tag{5.19}
$$
This implies that there exists a commutative diagram
$$
\CD
    A_q(K\backslash G/K) @>>> \Aq \\
    @VVV                           @VVV\\
    (\Uq/\Uq\kq+\kq\Uq)^\vee @>>> \Uq^\vee,
\endCD \tag{5.20}
$$
where the four arrows are all injective. 
It should be noted here that the left action of a
central element $C\in\Uq$ on ${\Cal H}=A_q(K\backslash G/K)$ 
corresponds by duality to the $\K$-endmorphism of the quotient space 
$\Uq/\Uq\kq+\kq\Uq$ induced from the right multiplication by $C$ in $\Uq$.
\par
The next step is to take the restriction mapping 
$A_q(K\backslash G/K) \to A(\T)$ 
into the duality argument. 
In what follows, we denote by $U_q(\frak{t})$ the commutative 
subalgebra $\K[q^h ( h\in P^*)]$ of $\Uq$, 
regarding it as the quantum analogue of the Lie algebra $\frak{t}$ 
of the diagonal subgroup $\T$. 
Then the pairing $(\ ,\ )$ between $\Uq$ and $\Aq$ induces a nondegenerate 
pairing between the subalgebra $U_q(\frak{t})$ and the quotient algebra 
$A(\T)$. 
For symmetry, we also use the notation $\zeta_k=q^{\ep_k}$ for $1\le k\le N$ 
and set 
$\zeta^h=\zeta_1^{\langle h,\ep_1\rangle}\cdots
\zeta_n^{\langle h,\ep_n\rangle}$ 
for any 
$h\in P^*$. 
With this notation, the pairing between $U_q(\frak{t})=\K[\zeta^{\pm 1}]$ 
and $A(\T)=\K[z^{\pm 1}]$ is described by 
$$
(\zeta^h, z^\ld) = q^{\langle h,\ld\rangle}\quad(h\in P^*,\ \ld\in P). 
\tag{5.21}
$$
Then the diagram (5.20) is complemented as follows:
$$
\CD
    A_q(K\backslash G/K) @>>> \Aq @>{\cdot|_\T}>> A(\T)\\
    @VVV                           @VVV           @VVV \\
    (\Uq/\Uq\kq+\kq\Uq)^\vee @>>> \Uq^\vee @>>> U_q(\frak{t})^\vee. 
\endCD \tag{5.22}
$$
Hence we have the commutative diagram
$$
\CD
    A_q(K\backslash G/K) @>{\cdot|_\T}>>  A(\T)\\
    @VVV                                  @VVV \\
    (\Uq/\Uq\kq+\kq\Uq)^\vee  @>>> U_q(\frak{t})^\vee,
\endCD \tag{5.23}
$$
where the arrow $(\Uq/\Uq\kq+\kq\Uq)^\vee \to U_q(\frak{t})^\vee$ 
is the transposition of the natural $\K$-linear mapping 
$U_q(\frak{t}) \to \Uq/\Uq\kq+\kq\Uq$ which represents the ``modulo reduction''
of an element of $U_q(\frak{t})$ by the subspace $\Uq\kq+\kq\Uq$.
\par 
The pairing (5.21) between the two algebras of Laurent polynomials
$\K[\zeta^{\pm 1}]$ and $\K[z^{\pm 1}]$ 
induces the ``multiplicative Fourier transform'' 
between $q$-difference operators acting on them. 
We denote by $\K[\zeta^{\pm 1};T_{q,\zeta}^{\pm 1}]$ the $\K$-algebra of 
$q$-difference operators in the form of finite sum 
$$
 Q(\zeta;T_{q,\zeta})= \sum_{\ld\in P} a_{\ld}(\zeta)T_{q,\zeta}^\ld\quad
(\ a_\ld(\zeta)\in\K[\zeta^{\pm 1}]\ ),
\tag{5.24}
$$
where $T_{q,\zeta}^\ld=T_{q,\zeta_1}^{\ld_1}\cdots T_{q,\zeta_N}^{\ld_N}$. 
Similarly we denote by $\K[z^{\pm 1};T_{q,z}^{\pm 1}]$ the algebra of 
$q$-difference operators in the variable $z=(z_1,\cdots,z_N)$.   
Between these two algebras of $q$-difference operators, 
there exists a unique anti-isomorphism of $\K$-algebras
$$
Q\mapsto\widehat{Q} : 
\K[\zeta^{\pm 1};T_{q,\zeta} ^{\pm 1}]\to\K[z^{\pm 1};T_{q,z}^{\pm 1}] 
\tag{5.25}
$$
such that 
$$
  \widehat{\zeta}_k=T_{q,z_k}, \ \widehat{T}_{q,\zeta_k} = z_k \ \ (1\le k\le N). 
\tag{5.26}
$$
We will call $\widehat{Q}=\widehat{Q}(z;T_{q,z})$ the {\it Fourier transform} 
of $Q=Q(\zeta;T_{q,\zeta})$. 
It is easy to see that  
$$
(Q(\zeta;T_{q,\zeta})f(\zeta),g(z)) =(f(\zeta),\widehat{Q}(z;T_{q,z})g(z)), 
\tag{5.27}
$$
for any $f(\zeta)\in \K[\zeta^{\pm 1}]$ and $g(z)\in\K[z^{\pm 1}]$.
\par
We are now ready to formulate our method to compute the radial component 
$C|_\T$ of a central element of $C\in \Uq$.
\Proposition{5.3}{
Let $C$ be an element of $\Uq$ such that $\kq C\subset\Uq\kq$. 
Then the left action of $C$ on $\Aq$ preserves the subalgebra 
${\Cal H}=A_q(K\backslash G/K)$ of $\kq$-biinvariants. 
Suppose that there exist a nonzero Laurent polynomial $a(z)\in\K[z^{\pm 1}]$ 
and a $q$-difference operator 
$Q(\zeta;T_{q,\zeta})\in\K[\zeta^{\pm 1};T_{q,\zeta}^{\pm 1}]$ 
such that 
$$
  (a(T_{q,\zeta})f)C \equiv Q(\zeta;T_{q,\zeta})f \ \mod \Uq\kq+\kq\Uq
\tag{5.28}
$$
for any Laurent polynomial $f=f(\zeta)$ in $U_q(\frak{t})=\K[\zeta^{\pm 1}]$. 
Then the radial component $C|_\T : {\Cal H}|_\T \to{\Cal H}|_\T $ 
is given by the Fourier transform $a(z)^{-1}\widehat{Q}(z;T_{q,z})$.  
Namely, one has 
$$
	(C.\varphi)|_\T = a(z)^{-1}\widehat{Q}(z;T_{q,z})(\varphi|_\T),
\tag{5.29}
$$
for any $\varphi\in {\Cal H}$. 
}
\Proof{
It is clear that the left action of $C$ preserves ${\Cal H}$. 
Note also that the right multiplication by $C$ in $\Uq$ 
preserves the subspace $\Uq\kq+\kq\Uq$.
Let $\varphi$ be an element of ${\Cal H}$. 
Then, for any element $f=f(\zeta)$  in $U_q(\frak{t})=\K[\zeta^{\pm 1}]$, 
we have
$$
\align
  (f, a(z)(C.\varphi)|_\T) 
&= (a(T_{q,\zeta})f, C.\varphi) 
= ((a(T_{q,\zeta})f(\zeta))C, \varphi) 
\tag{5.30}\\
&=(Q(\zeta;T_{q,\zeta})f,\varphi) 
=(f,\widehat{Q}(z;T_{q,z})(\phi|_\T)). 
\endalign
$$
This shows that 
$$
  a(z)(C.\varphi|_\T)=\widehat{Q}(z;T_{q,z})(\varphi|_\T)
\tag{5.31}
$$
as desired. 
}
The computation of the radial component $C|_\T$ is thus 
translated to the practical problem 
how to describe the reduction of elements in $U_q(\frak{t}) C$ 
modulo $\Uq\kq+\kq\Uq$, 
in terms of $q$-difference operators on $U_q(\frak{t})=\K[\zeta^{\pm 1}]$. 
%
%
\subsection{5.4. A recurrence formula related to Macdonald's $q$-difference operator}
%
Our computation of the radial component of $C_1$ will be carried out 
by using a recurrence formula for the ``symbol'' 
$$
\sum_{1\le k\le n} \xi_k
\frac{\Dt(x_1,\cdots,tx_k,\cdots,x_n)}{\Dt(x_1,\cdots,x_n)} 
\tag{5.32}
$$
of the $q$-difference operator in the left-hand side of (5.2).
We introduce a family of rational functions 
$F_{ij}=F_{ij}(x,\xi;t)\in \Q(t,x)[\xi] \ (1\le i\le j\le n)$ 
in $x=(x_1,\cdots,x_n)$ and $\xi=(\xi_1,\cdots,\xi_n)$.  
They are defined by the recurrence relations 
$$
\align
\text{\rm i)}&\quad F_{jj}=\xi_j \ \ (1\le j\le n), \tag{5.33}\\
\text{\rm ii)}&\quad
F_{ij}=\frac{1-t}{t(1-x_ix_j^{-1})}
\{\sum_{1\le k<j} F_{ik} - \sum_{i< k\le j} x_kx_j^{-1}F_{kj}\}
\ \ (i\le i<j\le n). 
\endalign
$$
\Lemma{5.4}{
The rational functions $F_{ij}=F_{ij}(x,\xi;t)$  
defined above have the explicit formulas
$$ 
F_{ij}(x,\xi;t) 
= t^{i-j}(1-t)^2 \sum_{i\le k\le j} \xi_k 
\frac{x_kx_j\Dt(x_i,\cdots,tx_k,\cdots,x_j)}
{(tx_k-x_i)(tx_k-x_j)\Dt(x_i,\cdots,x_j)}
\tag{5.34}
$$
for $1\le i\le j\le n$.
Furthermore, one has 
$$
\sum_{1\le i\le j\le n} t^{n-i}F_{ij}(x,\xi;t)
=\sum_{1\le k\le n} \xi_k 
\frac{\Dt(x_1,\cdots,tx_k,\cdots,x_n)}
{\Dt(x_1,\cdots,x_n)}.
\tag{5.35}
$$
}
In order to clarify the structure behind these formulas, 
we consider the upper unitriangular $n\times n$ matrix $A(t)$
defined by 
$$
A(t)=\sum_{1\le j\le n} e_{jj} + 
(1-t^{-1}) \sum_{1\le i<j\le n} e_{ij}.
\tag{5.36}
$$
We remark that the inverse of $A(t)$ is given by
$$
A(t)^{-1}= \sum_{1\le j\le n} e_{jj} + 
(1-t) \sum_{1\le i<j\le n} t^{i-j}e_{ij}.
\tag{5.37}
$$
With this matrix $A(t)$, we set $A(x;t)=A(t)\diag{x_1,\cdots,x_n}$.
Then the recurrence relations in (5.33) is equivalent to the commutativity
$$
	[A(x;t), F(x,\xi;t)]=0
\tag{5.38}
$$
for the matrix $F(x,\xi;t)=(F_{ij}(x,\xi;t))_{1\le i,j\le n}$ defined by 
setting $F_{ij}(x,\xi;t)=0$ for $i>j$. 
In other words, all the upper triangular matrices that commute with 
$A(x;t)$ are parametrized by their diagonal entries $\xi=(\xi_1,\cdots,\xi_n)$ 
and the other entries are explicitly determined by (5.34). 
We remark that the matrix $A(t)$ arise as the ``diagonal'' part of an $R$-matrix 
as we will see later. 
\par
A proof of Lemma 5.4 is given by diagonalizing the matrix $A(x;t)$.  
It is clear that there exists a unique upper unitriangular matrix 
$G(x;t)$ such that 
$$
	A(x;t)G(x;t)=G(x;t)\diag{x_1,\cdots,x_n}. 
\tag{5.39}
$$
With the matrix $G(x;t)$, $F(x,\xi;t)$ is realized as 
$$
	F(x,\xi;t)=G(x;t)\diag{\xi_1,\cdots,\xi_n} G(x;t)^{-1}.
\tag{5.40}
$$
\Lemma{5.5}{
{\rm (1)} The entries of the matrices
$G(x;t)=(G^+_{ij}(x;t))_{1\le i,j\le n}$ and
$G(x;t)^{-1}=(G^-_{ij}(x;t))_{1\le i,j\le n}$ are explicitly given by 
$$
\align
G^+_{ij}(x;t)
&=t^{i-j}\frac{(1-t)x_j\Dt(x_i,\cdots,x_{j-1},tx_j)}
{(x_i-tx_j)\Dt(x_i,\cdots,x_j)}
\tag{5.41}\\ 
G^-_{ij}(x;t)
&=t^{i-j}\frac{(t-1)x_j\Dt(tx_i,x_{i+1},\cdots,x_j)}
{(tx_i-x_j)\Dt(x_i,\cdots,x_j)}
\endalign
$$
for $1\le i\le j\le n$. 
\newline\noindent
{\rm (2)} For any $1\le i\le j\le n$, one has
$$
\align
&\sum_{i\le k\le j} t^{j-k} G^+_{kj}(x;t)=
\frac{\Dt(x_i,\cdots,x_{j-1},tx_j)}
{\Dt(x_i,\cdots,x_j)},
\tag{5.42}\\
&\sum_{i\le k\le j} G^-_{ik}(x;t)=t^{i-j}
\frac{\Dt(tx_i,x_{i+1},\cdots,x_j)}
{\Dt(x_i,\cdots,x_j)}.
\endalign
$$
}
\Proof{
Formula (5.39) is equivalent to the recurrence relations
$$
(1-x_ix^{-1}_j)G^+_{ij}(x;t)
=(1-t^{-1})\sum_{i<k\le j}x_kx^{-1}_jG^+_{kj}(x;t)
\tag{5.43}
$$
for $i<j$ with initial condition $G^+_{jj}=1$. By comparing (5.43) 
and the formula obtained from (5.43) replacing $i+1$ for $i$, we obtain
$$
G^+_{ij}(x;t)=t^{-1}\frac{x_{i+1}-tx_j}{x_i-x_j} G^+_{i+1,j}(x;t). 
\tag{5.44}
$$
Hence we compute 
$$
G^+_{ij}(x;t)=t^{i-j}
\frac{\prod_{i<k\le j}x_k-tx_j}{\prod_{i\le k<j}x_k-x_j}, 
\tag{5.45}
$$
to get (5.41) for $G^+_{ij}(x;t)$. 
On the other hand, from $A(x;t)^{-1}G(x;t)\diag{x}=G(x;t)$ we have
$$
(x_i^{-1}x_j-1)G^+_{ij}(x;t)=(1-t)\sum_{i<k\le j} t^{i-k}G^+_{kj}(x;t). 
\tag{5.46}
$$
Formula (5.42) for $G^+_{ij}(x;t)$ is an easy consequence of (5.41) and 
(5.45). 
The same method is available to prove formulas for $G^-_{ij}(x;t)$. 
}
Formulas (5.34) and (5.35) of Lemma 5.4 follows immediately from (5.41) and (5.42) of Lemma 5.5, 
respectively. 
\par 
From the above description (5.40), we also see that the matrix $F(x,\xi;t)$ 
has the multiplicative property
$$
	F(x,1;t)=\id,\  F(x,\xi;t) F(x,\eta;t)= F(x,\xi\eta;t) \ \ \text{\rm and} \ \  F(x,x;t)=A(x;t).
\tag{5.47}
$$
We remark that the formula (5.35) is also related 
to the {\it Hall-Littlewood polynomials}
$P_{(\ell)}(x;t) \ (\ell=1,2,\cdots)$ for the Young diagrams of one row.  
In fact, by the substitution $\xi_k=x^{\ell}_k \ (1\le k\le n)$, the right-hand side of (5.35) 
gives rise to
$$
	P_{(\ell)}(x;t) = \sum_{k=1}^{n} {x_k}^{\ell} 
\frac{\Dt(x_1,\cdots,tx_k,\cdots,x_n)}{\Dt(x_1,\cdots,x_n)},
\tag{5.48}
$$
for $\ell=1,2,\cdots$ (see \cite{M1}). 
%
\subsection{5.5. Proof of Theorem 5.2}
%
%
We determine the radial component of the central element $C_1\in\Uq$
in the two cases (SO) and (Sp). 
Taking an arbitrary $h\in P^*$, we will compute the reduction of the element
$q^hC_1=\zeta^hC_1$ modulo $\Uq\kq+\kq\Uq$ to find an explicit representative 
of its modulo class in $U_q(\frak{t})$. 
\par
In view of (5.12), 
we begin with looking for relations among the modulo classes of 
$q^h L^+_{ij} S(L^-_{k\ell})$. 
For this purpose we use the commutation relations 
$$
L^+_1 S(L^-_2)^t {R^+_{12}}^{t_2} ={R^+_{12}}^{t_2} S(L^-_2)^t L^+_1. 
\tag{5.49}
$$
Note also that we have 
$$
	q^h L^+ = H^{-1} L^+ H q^h, \quad q^h S(L^-) = H^{-1} S(L^-) H q^h,
\tag{5.50}
$$
for the constant matrix 
$H=\diag{q^{\langle h,\ep_1\rangle},\cdots,q^{\langle h,\ep_n\rangle}}$.
By using (5.49) and (5.50), we compute 
$$
\align
	q^h L^+_1 S(L^-_2)^t {R^+_{12}}^{t_2}
& = 	q^h {R^+_{12}}^{t_2} S(L^-_2)^t L^+_1 \tag{5.51}\\
& = 	{R^+_{12}}^{t_2} H_2 S(L^-_2)^t H^{-1}_2 q^h L^+_1 \\
&\equiv	{R^+_{12}}^{t_2} H_2 J^{-1} L^+_2 J H^{-1}_2 q^h L^+_1 \\
&\equiv	q^h {R^+_{12}}^{t_2} H_2 J^{-1}_2 H_2 L^+_2 H^{-1}_2 J_2 H^{-1}_2 L^+_1\\
&\equiv	q^h {R^+_{12}}^{t_2} H_2 J^{-1}_2 H_2 L^+_2 H^{-1}_2 J_2 H^{-1}_2 J_1 S(L^-_1)^t J^{-1}_1
\endalign
$$
modulo $\Uq\kq+\kq\Uq$. 
Hence we have 
$$
q^h L^+_1 S(L^-_2)^t {R^+_{12}}^{t_2} J_1 H_2 J^{-1}_2 H_2
\equiv	q^h {R^+_{12}}^{t_2} J_1 H_2 J^{-1}_2 H_2 L^+_2 S(L^-_1)^t. 
\tag{5.52}
$$
By applying the matrix $P_{12}$ from the right, we have 
$$
q^h L^+_1 S(L^-_2)^t {R^+_{12}}^{t_2} P_{12} J_2 H_1 J^{-1}_1 H_1 
\equiv	q^h {R^+_{12}}^{t_2} P_{12} J_2 H_1 J^{-1}_1 H_1 L^+_1 S(L^-_2)^t. 
\tag{5.53}
$$
This means that the constant matrix 
$\widetilde{R}={R^+_{12}}^{t_2} P_{12} J_2 H_1 J^{-1}_1 H_1$ 
``intertwines'' the modulo class of $q^h L^+_1 S(L^-_2)^t $. 
From this formula (5.53), we can derive necessary informations for determining the 
modulo class of $q^h C_1$, by extracting its ``diagonal part''.  
\par
Let $W=\bigoplus_{j=1}^n \K u_j$ be the $n$-dimensional $\K$-vector space 
with canonical basis $\{u_1,\cdots,u_n\}$. 
For each case, we define below two $\K$-linear mappings 
$$
\iota_W : W \to V\oxK V\ \ \text{and} \ \ \pi_W : V\oxK V\to W.  
\tag{5.54}
$$
By using these mapping, we pull back the equality (5.53) by the composition 
$$
W @>\iota_W>> V\oxK V @>>> V\oxK V @>\pi_W>>  W 
\tag{5.55}
$$
to obtain the statement for $n\times n$ matrix 
$Z=\pi_W\circ (L^+_1 S(L^-_2)^t)\circ\iota_W$ with entries in $\Uq$. 
From now on, we treat the two cases (SO) and (Sp), separately. 
\par\medpagebreak\par\noindent
Case (SO): We define $\iota_W$ and $\pi_W$ by 
$$
\iota_W(u_k)= v_k\ox v_k \ (1\le k\le n), \ \ 
\pi_W(v_i\ox v_j)=\delta_{ij}u_j \ (1\le i,j\le n).
\tag{5.56}
$$
Under this definition, the matrix $Z=\pi_W\circ (L^+_1 S(L^-_2)^t)\circ\iota_W$ takes the form 
$$
	Z=(Z_{ij})_{1\le i,j\le n}, \ \ Z_{ij} = L^+_{ij} S(L^-_{ji}) \ (1\le i,j\le n). 
\tag{5.57}
$$
\Lemma{5.6.A}{
Set $\widetilde{R}={R^+_{12}}^{t_2} P_{12} J_2 H_1 J^{-1}_1 H_1$.  Then we have 
$$ 
\align
&\pi_W\circ \widetilde{R} = 
qA(q^2)\diag{q^{2\langle h,\ep_1\rangle},\cdots,q^{2\langle h,\ep_n\rangle}}\circ\pi_W, 
\tag{5.58}\\
&\widetilde{R}\circ \iota_W = 
q\iota_W \circ A(q^2)\diag{q^{2\langle h,\ep_1\rangle},\cdots,q^{2\langle h,\ep_n\rangle}},
\endalign
$$
where $A(t)$ is the matrix defined by (5.36).
}
Lemma 5.6.A is clear from the expression of $ {R^+_{12}}^{t_2} P_{12}$:
$$
{R^+_{12}}^{t_2} P_{12} 
= \sum_{i,j} e_{ij}\ox e_{ji} q^{\delta_{ij}} + (q-q^{-1}) \sum_{i<j} e_{ij}\ox e_{ij}.
\tag{5.59}
$$
\par
By using (5.58), we obtain the following formula for $Z$ from (5.53):
$$
 q^h Z  A(q^{2\langle h,\ep\rangle};q^2) \equiv  q^h A(q^{2\langle h,\ep\rangle};q^2) Z \quad \mod \Uq\kq+\kq\Uq .
\tag{5.60}
$$
Here we used the notation 
$$
A(q^{2\langle h,\ep\rangle};q^2)=A(q^2)\diag{q^{2\langle h,\ep\rangle}}, \ 
q^{2\langle h,\ep\rangle}=(q^{2\langle h,\ep_1\rangle},\cdots, q^{2\langle h,\ep_n\rangle}).
\tag{5.61}
$$
This means that the modulo classes of $q^h Z_{ij}$ satisfy the same recurrence relation 
as that of $F_{ij}(q^{2\langle h,\ep\rangle},\xi;q^2)$ that we discussed in Section 5.4. 
Note that, as to the diagonal entries, we have 
$$
	q^h Z_{jj} \equiv q^h q^{2\ep_j} \quad \mod \Uq\kq+\kq\Uq, 
\tag{5.62}
$$
for $1\le j \le n$. 
Taking $\xi_j=\Dt(q^{2\langle h,\ep\rangle}) q^{2\ep_j} (1\le j\le n)$, we conclude inductively that 
$$
 \Dt(q^{2\langle h,\ep\rangle}) q^h Z_{ij}
\equiv q^h F_{ij}(q^{2\langle h,\ep\rangle},\Dt(q^{2\langle h,\ep\rangle})q^{2\ep};q^2) \ \ (1\le i,j\le n),
\tag{5.63}
$$
where $q^{2\ep}=(q^{2\ep_1},\cdots,q^{2\ep_n})$. 
In fact, the both sides of (5.63) have the same initial values for $j-i=0$ 
and satisfy the same recurrence relations for $j-i>0$. 
Recall that 
$$
C_1=\sum_{1\le i,j\le n} q^{2(n-i)} L^+_{ij} S(L^-_{ji}) 
=\sum_{1\le i,j\le n} q^{2(n-i)} Z_{ij}. 
\tag{5.64}
$$
Hence we obtain the following expression for the modulo class of 
$\Dt(q^{2\langle h,\ep\rangle})C_1$ by Lemma 5.4: 
$$
\align
&\Dt(q^{2\langle h,\ep_1\rangle},\cdots, q^{2\langle h,\ep_n\rangle}) q^h C_1
\tag{5.65}\\
&\quad\equiv \sum_{k=1}^n q^{h+2\ep_k} 
\Dt(q^{2\langle h,\ep_1\rangle},\cdots, q^2q^{2\langle h,\ep_k\rangle},\cdots,
q^{2\langle h,\ep_n\rangle}).
\endalign
$$
In terms of the operators in  the variables $\zeta=(\zeta_1,\cdots,\zeta_n)$, 
this formula can be rewritten as 
$$ 
	\Dt(T_{q,\zeta_1}^2,\cdots,T_{q,\zeta_n}^2)  f(\zeta) C \equiv 
\sum_{k=1}^n \zeta^2_k \ 
\Dt(T_{q,\zeta_1}^2,\cdots, q^2T_{q,\zeta_k}^2,\cdots,T_{q,\zeta_n}^2) f(\zeta)
\tag{5.66}
$$
for any $f(\zeta)\in\K[\zeta^{\pm 1}]=U_q(\frak{t})$. 
By Proposition 5.3, we finally get the explicit formula for the radial component 
$$
	C_1|_\T = \sum_{k=1}^n 
\frac{\Dt(z^2_1,\cdots,q^2z^2_k,\cdots,z^2_n)} {\Dt(z^2_1,\cdots,z^2_n)} T_{q,z_k}^2
\tag{5.67}
$$
as the Fourier transform of (5.66). 
On the subalgebra 
$\K[x_1,\cdots,x_n]^{\frak{S}_n}$ of ${\Cal H}$, 
consisting of symmetric polynomials in  $x_1=z^2_1, \cdots, x_n=z^2_n$, 
this reduces to the $q$-difference operator 
$$
D_1 = \sum_{k=1}^n 
\frac{\Dt(x_1,\cdots,q^2x_k,\cdots,x_n)} {\Dt(x_1,\cdots,x_n)}T_{q^4,x_k}
\tag{5.68}
$$
in the variables $x=(x_1,\cdots,x_n)$. 
This completes the proof of Theorem 5.2 for Case (SO). 
\par\medpagebreak\par\noindent
Case (Sp): In Case (Sp), we define the linear mapping
$\iota_W : W \to V\oxK V$ by 
$$
\iota_W(u_k) = v_{2k-1}\ox v_{2k-1} + v_{2k}\ox v_{2k} \ \ (1\le k\le n),
\tag{5.69}
$$
The linear mapping $\pi_W : V\oxK V \to W$ is defined as follows:
$$
\pi_W(v_{2k-1}\ox v_{2k-1})=q u_k,\ 
\pi_W(v_{2k}\ox v_{2k})=q^{-1} u_k \ \ (1\le k\le n)
\tag{5.70}
$$
and $\pi_W(v_i\ox v_j)=0$ for the other pairs $(i,j)$.  
With these linear mappings, we see that the entries of the matrix 
$Z=\pi_W\circ L^+_1S(L^-_2)\circ\iota_W$ are given by 
$$
\align
Z_{ij}=&
q L^+_{2i-1,2j-1}S(L^-_{2j-1,2i-1}) +
q L^+_{2i-1,2j}S(L^-_{2j,2i-1})
\tag{5.71}\\
&\ \ + q^{-1} L^+_{2i,2j-1}S(L^-_{2j-1,2i}) +
q^{-1} L^+_{2i,2j}S(L^-_{2j,2i}), 
\endalign
$$
for $1\le i,j\le n$.
We also see that the central element  $C_1$  is rewritten in terms of $Z_{ij}$ 
as follows:
$$
C_1=q\sum_{1\le i,j\le n} q^{4(n-i)} Z_{ij}. 
\tag{5.72}
$$
By a direct computation, one can prove 
\Lemma{5.6.B}{
For the matrix $\widetilde{R}={R^+_{12}}^{t_2} P_{12} J_2 H_1 J^{-1}_1 H_1$, we have 
$$ 
\align
&\pi_W\circ \widetilde{R} = 
-q^2 A(q^4)\diag{q^{\langle h,\ep_1+\ep_2\rangle},\cdots,q^{\langle h,\ep_{2n-1}+\ep_{2n}\rangle}}\circ\pi_W, 
\tag{5.73}\\
&\widetilde{R}\circ \iota_W = 
-q^2 \iota_W \circ A(q^4)\diag{q^{\langle h,\ep_1+\ep_2\rangle},\cdots,q^{\langle h,\ep_{2n-1}+\ep_{2n}\rangle}}.
\endalign
$$
}
Setting $\widetilde{\ep}_k=\ep_{2k-1}+\ep_{2k}$ for $1\le k\le n$, we will 
use below the notations
$ q^{\widetilde{\ep}}=(q^{\widetilde{\ep}_1},\cdots,q^{\widetilde{\ep}_n}) $
and 
$ q^{\langle h,\widetilde{\ep}\rangle}=
(q^{\langle h,\widetilde{\ep}_1\rangle},\cdots,q^{\langle h,\widetilde{\ep}_n\rangle})$. 
By using (5.73), we obtain the commutation relation  
$$
 q^h Z  A(q^{\langle h,\widetilde{\ep}\rangle};q^4) \equiv  q^h A(q^{\langle h,\widetilde{\ep}\rangle};q^4) Z \quad \mod \Uq\kq+\kq\Uq 
\tag{5.74}
$$
from (5.53). 
As for the diagonal entries of $Z$, we have 
$$
 q^h Z_{jj} \equiv q^h(q q^{2\ep_{2j-1}}+q^{-1}q^{2\ep_{2j}})
\equiv q^h(q+q^{-1})q^{\widetilde{\ep}_j}
\quad \mod \Uq\kq+\kq\Uq,
\tag{5.75}
$$
for $1\le j\le n$. 
Hence we have 
$$
 \Dt(q^{\langle h,\widetilde{\ep}\rangle}) q^h Z_{ij}
\equiv q^h (q+q^{-1})F_{ij}(q^{\langle h,\widetilde{\ep}\rangle},\Dt(q^{\langle h,\widetilde{\ep}\rangle})q^{\widetilde{\ep}};q^4) \ \ (1\le i,j\le n),
\tag{5.76}
$$
by a similar inductive argument as in Case (SO). 
This leads to the following expression for the modulo class of 
$\Dt(q^{\langle h,\widetilde{\ep}\rangle})C_1$ by (5.72): 
$$
\align
&\Dt(q^{\langle h,\widetilde{\ep_1}\rangle},\cdots, q^{\langle h,\widetilde{\ep}_n\rangle}) q^h C_1 
\tag{5.77}\\
&\quad\equiv (1+q^2) 
\sum_{1\le i,j\le n} q^{4(n-i)} 
F_{ij}(q^{\langle h,\widetilde{\ep}\rangle},
\Dt(q^{\langle h,\widetilde{\ep}\rangle})q^{\widetilde{\ep}};q^4) \\
&\quad\equiv (1+q^2)\sum_{k=1}^n q^{h+\widetilde{\ep}_k} 
\Dt(q^{\langle h,\widetilde{\ep}_1\rangle},\cdots, q^4 q^{\langle h,\widetilde{\ep}_k\rangle},
\cdots,q^{\langle h,\widetilde{\ep}_n\rangle}).
\endalign
$$
In terms of the operators in the variables $\zeta=(\zeta_1,\cdots,\zeta_{2n})$, 
this formula can be rewritten as 
$$ 
\align
&\Dt(\widetilde{T}_1,\cdots,\widetilde{T}_n) f(\zeta) C 
\tag{5.78}
\\
&\quad\equiv 
(1+q^2)\sum_{k=1}^n \zeta_{2k-1}\zeta_{2k} \ 
\Dt(\widetilde{T}_1,\cdots, q^4\widetilde{T}_k,\cdots,\widetilde{T}_n) f(\zeta)
\endalign
$$
for any $f(\zeta)\in\K[\zeta^{\pm 1}]=U_q(\frak{t})$, where 
$\widetilde{T}_k=T_{q,\zeta_{2k-1}}T_{q,\zeta_{2k}}$ for $1\le k\le n$. 
Hence, by Proposition 5.3, we get the explicit formula for the radial component 
$$
C_1|_\T = (1+q^2)\sum_{k=1}^n 
\frac{\Dt(z_1z_2,\cdots,q^4z_{2k-1}z_{2k},\cdots,z_{2n-1}z_{2n})} {\Dt(z_1z_2,\cdots,z_{2n-1}z_{2n})} T_{q,z_{2k-1}}T_{q,z_{2k}}
\tag{5.79}
$$
as the Fourier transform of (5.78). 
On the subalgebra $\K[x_1,\cdots,x_n]\subset{\Cal H}$
of symmetric polynomials in
$x_1=z_1z_2, \cdots, x_n=z_{2n-1}z_{2n}$, 
this reduces to the 
$q$-difference operator 
$$
D_1 = (1+q^2)\sum_{k=1}^n 
\frac{\Dt(x_1,\cdots,q^4x_k,\cdots,x_n)} {\Dt(x_1,\cdots,x_n)}T_{q^2,x_k}. 
\tag{5.80}
$$
This completes the proof of Theorem 5.2 for Case (Sp). 
%
\section{\S 6. Scalar product and orthogonality}
%
In this section, we will discuss the orthogonality relations 
for Macdonald's symmetric polynomials $P_{\mu}(x;q,t)$ with 
$t=q^{\frac{1}{2}}$ or $t=q^2$ which are obtained from their 
interpretation as zonal spherical functions on quantum 
homogeneous spaces. 
%
\subsection{6.1. Invariant functional and Schur's orthogonality}
%
Recall that there exists a unique homomorphism $h_G: \Aq \to \K$ of 
$\Uq$-bimodules with $h_G(1)=1$.  
This invariant functional, corresponding to the Haar measure of the unitary group 
$\text{\rm U}(N)$, is given as the projection $\Aq \to W(0)=\K$ 
in the decomposition (1.23).  
We remark that the invariance of $h_G$ means that, for any element $a\in\Uq$, 
one has 
$$
h_G(a.\varphi)=h_G(\varphi.a)= \vep(a)h_G(\varphi)\quad (\varphi\in\Aq).
\tag{6.1}
$$
By using the $\ast$-operation of $\Aq$, we define a hermitian form 
$\langle\ ,\ \rangle_G$ on $\Aq$ by the formula
$$
\langle\varphi,\psi\rangle_G=h_G(\varphi^*\,\psi)\quad\text{for}\ \varphi,\psi\in\Aq. 
\tag{6.2}
$$
From the invariance of $h_G$ it follows that $\langle\ ,\ \rangle_G$ 
is invariant in the sense that
$$
\langle \varphi, a.\psi\rangle_G = \langle a^*.\varphi,\psi\rangle_G 
\tag{6.3}
$$
for any $a\in\Uq$ and $\varphi,\psi\in\Aq$.
It is known that this hermitian form is nondegenerate, and induces 
a positive definite hermitian form when $q$ is specialized to a real number with $|q|\ne 0,1$ (see \cite{NYM}, for instance). 
\par
The orthogonality relations for our zonal spherical functions $\varphi(\ld)$\ 
$(\ld\in P^+_\k)$ come from Schur's orthogonality relations for the matrix 
elements of irreducible representations $V(\ld)$.  
For each dominant integral weight $\ld\in P^+$, we fix a nondegenerate 
$\text{U}_q(N)$-invariant hermitian form $\langle\ ,\ \rangle$ on $V(\ld)$ as 
in Section 4.4.  
In what follows, we set 
$$
\rho=\sum_{k=1}^N (N-k)\ep_k.
\tag{6.4}
$$
The orthogonality relations for the matrix elements $\phi_\ld(u,v)$ for 
$u,v\in V(\ld)$ are formulated as follows.  
\Proposition{6.1}{
Let $\ld,\mu$ be two dominant integral weights in $P^+$ and 
take nonzero vectors $u,v\in V(\ld)$ and $u',v'\in V(\mu)$.  
If $\ld\ne\mu$, one has
$$
	\langle \phi_\ld(u,v),\phi_\mu(u',v')\rangle_G=0.  
\tag{6.5}
$$
If $\ld=\mu$, one has,
$$
\langle \phi_\ld(u,v),\phi_\ld(u',v')\rangle_G
= \frac{1}{d(\ld)} \langle q^\rho.u',q^\rho.u\rangle \langle v,v'\rangle,
\tag{6.6} 
$$
where $d(\ld)$ is the following principal specialization of the 
Schur function $s_\ld$:
$$
  d(\ld)=s_\ld(q^{2\rho}), \quad q^{2\rho}=(q^{2(N-1)},q^{2(N-2)},\cdots,1). 
\tag{6.7}
$$
}
We omit the proof of Proposition 6.1, 
since it is a variant of Schur's orthogonality already given in 
Woronowicz \cite{W} and \cite{NYM}. 
%
%
\subsection{6.2. Orthogonality of zonal spherical functions}
%
The zonal spherical function $\varphi(\ld)$ $(\ld\in P^+_\k)$
is a matrix element of the representation $V(\ld)$.  
From this fact, it directly follows that,
if $\ld$ and $\ld'$ are two distinct weights in $P^+_\k$, one has 
the orthogonality relation 
$ \langle\varphi(\ld), \varphi(\ld')\rangle_G=0. $
This fact is also proved by the fact that the central element 
$C_1$ of $\Uq$ is self-adjoint with respect to the hermitian form 
$\langle\ ,\ \rangle_G$. 
\par
For the description of the square length 
$ \langle\varphi(\ld), \varphi(\ld)\rangle_G $
of $\varphi$, we need to specialize appropriately the parameter 
$a\in (\K^*)^n$ involved in the definition of the coideal $\kq=\kq(a)$. 
In what follows, we set
$$
\align
\text{Case (SO):}&\quad a=(q^{(n-1)},q^{(n-2)},\cdots,1),\tag{6.8}\\
\text{Case (Sp):}&\quad a=(q^{2(n-1)},q^{2(n-2)},\cdots,1). 
\endalign
$$
\Proposition{6.2}{ For the special value of $a\in(\K^*)^n$ in (6.8), 
the algebra ${\Cal H}=A_q(K\backslash G/K)$ of $\kq$-biinvariant elements 
is a $\ast$-subalgebra of $\Aq$.  
}
\Proof{
We first consider Case (SO).  
In view of Lemma 4.13, we have only to show that $\varphi(2\Ld_r)^*$ 
$(1\le r\le n-1)$ and $(\detq(T)^{\pm 1})^*$ belong to ${\Cal H}$ again. 
As to $\detq(T)^{\pm 1}$, this statement is clear since 
$\detq(T)^*=\detq(T)^{-1}$.  
Note that, under the specialization (6.8), the element $\varphi(2\Ld_r)$ 
takes the form
$$
    \varphi(2\Ld_r)= \sum_{|I|=|J|=r} (\xi^I_J)^2 q^{||I||-||J||},
\tag{6.9}
$$
where $||I||=\sum_{i\in I} i$. 
Hereafter we use the notation 
$\xi^I_J=\xi^{i_1\cdots i_r}_{j_1\cdots j_r}$
for $I=\{i_1<\cdots<i_r\}$ and $J=\{j_1<\cdots<j_r\}$. 
Recall from \cite{NYM} that, for each $I,J\subset\{1,2,\cdots,N\}$ with 
$|I|=|J|=r$, we have 
$$
(\xi^I_J)^* = (-q)^{-||I||+||J||} \xi^{I^c}_{J^c} \detq(T)^{-1},
\tag{6.10}
$$
where $I^c$ stands for the complement of $I$ in $\{1,2,\cdots,N\}$. 
Hence we compute 
$$
 \varphi(2\Ld_r)^* = \detq(T)^{-2}\sum_{|I|=|J|=r} (\xi^{I^c}_{J^c})^2 q^{-||I||+||J||}
 = \varphi(2\Ld_{n-r}-2\Ld_n)
\tag{6.11}
$$
by (4.36). 
This proves that ${\Cal H}$ is closed under the $\ast$-operation in Case (SO).
\newline
In Case (Sp), the element $\varphi(\Ld_{2r})$ can be written in the form 
$$
  \varphi(\Ld_{2r})=\sum_{|I|=|J|=r} \widetilde{\xi}^I_J q^{2(||I||-||J||)},
\tag{6.12}
$$
for each $1\le r\le n$.
In (6.12), we used the notation $\widetilde{\xi}^I_J$ to refer 
$\xi^{2i_1-1,2i_1,\cdots,2i_r-1,2i_r}_{2j_1-1,2j_1,\cdots,2j_r-1,2j_r}$ 
for two subsets 
$I=\{i_1<\cdots<i_r\}$ and $J=\{j_1<\cdots<j_r\}$ of $\{1,2,\cdots,n\}$. 
Since 
$(\widetilde{\xi}^I_J)^* 
= (-q)^{4(-||I||+||J||)} \widetilde{\xi}^{I^c}_{J^c} \detq(T)^{-1},$
we compute
$$
 \varphi(\Ld_{2r})^* = \detq(T)^{-1}\sum_{|I|=|J|=r} 
	\widetilde{\xi}^{I^c}_{J^c} q^{-2||I||+2||J||}
 = \varphi(\Ld_{2(n-r)}-\Ld_{2n}),
\tag{6.13}
$$
where $I^c$ stands for the complement of $I$ in $\{1,2,\cdots,n\}$.
Hence Lemma 4.13 implies that ${\Cal H}$ is closed under the $\ast$-operation 
also in Case (Sp).  
}
Another way to prove Proposition 6.2 is to show the left ideal 
$\Uq\kq$ and the right ideal $\kq\Uq$ are both stable under the 
involution $\tau$ (see (1.31)). 
This fact can be checked directly by using the generator system 
described in Proposition 2.4. 
\par
Under the specialization of the parameters $a=(a_1,\cdots,a_n)$ as in (6.8), 
we have an expression of the square length 
$\langle\varphi(\ld), \varphi(\ld)\rangle_G$
in terms of the principal specialization of $\varphi(\ld)$. 
\Proposition{6.3}{
The zonal spherical functions $\varphi(\ld)$ $(\ld\in P^+_\k)$
form an orthogonal basis of the algebra ${\Cal H}$. 
For each $\ld\in P^+_\k$, 
the square length of $\varphi(\ld)$ is given by the 
formula 
$$
\langle\varphi(\ld), \varphi(\ld)\rangle_G=\frac{c(\ld)^2}{d(\ld)},
\tag{6.14}
$$
where $c(\ld)=(q^\rho,\varphi(\ld))$ stands for the value 
of $\varphi(\ld)$ at the point $q^\rho\in\T$ of the diagonal subgroup. 
}
For the proof of Proposition 6.3, we describe how the 
$\ast$-operation acts on the coideal $\kq$. 
\Lemma{6.4}{ The coideal $\kq=\kq(a)$ for the parameter $a$ of (6.8), one has 
$\kq^* = q^{-\rho}\kq q^\rho$.  
Accordingly, if $w\in V(\ld)$ is a $\kq$-fixed vector, then $q^{-\rho}.w$ gives 
a $\kq^*$-fixed vector. 
}
\Proof{
Setting $D=\diag{q^{(N-1)},q^{(N-2)},\cdots,1}$, we have
$$
\align
q^{-\rho}Mq^\rho &= DL^+D^{-1} - JD^{-1}S(L^-)^tDJ^{-1}\tag{6.15}\\
&= D(L^+ - D^{-1}JD^{-1}S(L^-)^t DJ^{-1}D)D^{-1}. 
\endalign
$$
For the value (6.8) of $a$, one can check easily 
the matrix $D^{-1}JD^{-1}=D^{-1}J(a)D^{-1}$ 
is a scalar multiple of $J(a^{-1})$. 
This shows that $q^{-\rho}\kq(a)q^\rho= \kq(a^{-1}) = \kq(a)^*$ 
as desired. 
}
\Proofof{Proposition 6.3}{
Let us take the vectors $w(\ld)$, $w^*(\ld)$ as in the definition (4.4). 
Then by Lemma 6.4, we must have 
$w^*(\ld)=q^{\langle\rho,\ld\rangle} q^{-\rho}.w(\ld)$, so that 
$$
\varphi(\ld)
=q^{\langle\rho,\ld\rangle}
\frac{\phi_\ld(q^{-\rho}.w,w)}{\langle u,u\rangle}
=q^{\langle\rho,\ld\rangle}
\frac{\phi_\ld(w,w).q^{-\rho}}{\langle u,u\rangle},
\tag{6.16}
$$
where $w=w(\ld)$ and $u=u(\ld)$.
Hence we have,
$$
c(\ld)=(q^\rho,\varphi(\ld))=\vep(\varphi(\ld).q^\rho)
=q^{\langle\rho,\ld\rangle}
\frac{\langle w,w\rangle}{\langle u,u\rangle}
\tag{6.17}
$$
On the other hand, we have
$$
\align
\langle\varphi(\ld),\varphi(\ld)\rangle_G
&=q^{2\langle\rho,\ld\rangle}
\frac{\langle\phi_\ld(q^{-\rho}.w,w),\phi_\ld(q^{-\rho}.w,w)\rangle_G}
{{\langle u,u\rangle}^2} \tag{6.18}\\
&=\frac{1}{d(\ld)} q^{2\langle\rho,\ld\rangle}
\frac{\langle w,w\rangle^2}{\langle u,u\rangle^2}.
\endalign
$$
by Proposition 6.1. 
Comparing (6.17) and (6.18), we obtain the expression of (6.14).
}
%
%
\subsection{6.3. Description of orthogonality on the diagonal subgroup}
%
From this subsection on, we take the field $\K=\C$ of complex numbers 
as the ground field and assume that $q$ is a real number with $0<|q|<1$. 
\par
We now consider to describe the invariant functional $h_G: {\Cal H}\to \C$ 
on the diagonal subgroup $\T=(\C^*)^N$.  
For this purpose we use the subalgebra 
$\C[x]^{\frak{S}_n}=\C[x_1,\cdots,x_n]^{\frak{S}_n}$ of 
${\Cal H}|_\T$ of symmetric polynomials in the following variables: 
$$ 
\align
\text{Case (SO):}&\quad x_1=z^2_1,\cdots, x_n=z^2_n,\tag{6.19}\\
\text{Case (Sp):}&\quad x_1=z_1z_2,\cdots, x_n=z_{2n-1}z_{2n}. 
\endalign
$$
In what follows, we denote by ${\Cal R}$ the subalgebra of ${\Cal H}$ 
such that ${\Cal R}|_\T= \C[x]^{\frak{S}_n}$. 
Let us consider only the highest weights $\ld\in P^+_\k$ that are 
parametrized by the partitions $\mu=(\mu_1,\cdots,\mu_n)$ 
($\mu_1\ge\cdots\ge\mu_n\ge 0$) as follows:
$$ 
\align
\text{Case (SO):}&\quad \ld=\sum_{k=1}^n 2\mu_k \ep_k, \tag{6.20}\\
\text{Case (Sp):}&\quad \ld=\sum_{k=1}^n \mu_k (\ep_{2k-1}+\ep_{2k}). 
\endalign
$$ 
Note that this parametrization can be graphically described as 
the duplication of Young diagrams, in the horizontal direction in 
Case (SO) and in the vertical direction in Case (Sp).
As we noticed in Remark 4.14, the zonal spherical functions 
$\varphi(\ld)$ parametrized by partitions as above form a $\C$-basis of 
the subalgebra ${\Cal R}$. 
Furthermore, we know by Theorem 5.1 that, for each $\varphi(\ld)\in {\Cal R}$, 
its restriction to the diagonal subgroup $\T$ coincides with the Macdonald
symmetric polynomial $P_\mu(x)=P_\mu(x;q^4,q^2)$ in Case (SO), and 
with $P_\mu(x)=P_\mu(x;q^2,q^4)$ in Case (Sp), respectively.  
Note also that $P_\mu(x)$ form a $\C$-basis for the algebra 
$\C[x]^{\frak{S}_n}$ as $\mu$ ranges over all partitions. 
\par
On the subalgebra ${\Cal R}$, the invariant functional $h_G: {\Cal R}\to \C$ 
is described in terms of the following meromorphic function 
on the algebraic torus $(\C^*)^n$ : 
$$
w(x;q,t)=\prod_{1\le i<j\le n} 
\frac{(x_i/x_j;q)_\infty(x_j/x_i;q)_\infty}
{(tx_i/x_j;q)_\infty(tx_j/x_i;q)_\infty}, 
\tag{6.21}
$$
where $(a;q)_\infty=\prod_{k=0}^\infty (1-aq^k)$. 
For a holomorphic function $F(x)$ defined in a neighborhood 
of the torus 
$T=\{ x=(x_1,\cdots,x_n)\in (\C^*)^n ; |x_1|=\cdots=|x_n|=1\}$, 
we use the notation 
$$ 
[F(x)]_1 
=\left(\frac{1}{2\pi\sqrt{-1}}\right)^n\int_T F(x_1,\cdots,x_n) 
\frac{dx_1\cdots dx_n}{x_1\cdots x_n}
\tag{6.22}
$$
to refer the constant term in the Laurent expansion of $F(x)$. 
\Proposition{6.5}{
Let $\varphi$ be an element in the subalgebra ${\Cal R}$ of ${\Cal H}$ 
and $F(x)=\varphi|_\T$ the corresponding symmetric polynomial 
in $\C[x]$.  
Then the value $h_G(\varphi)$ of the invariant functional is described by
the formula 
$$
h_G(\varphi)=\frac{\left[F(x)w(x)\right]_1}{\left[w(x)\right]_1},
\tag{6.23}
$$
where the weight function $w(x)$ is given by
$$
\align
\text{Case (SO):}&\quad w(x)=w(x;q^4,q^2), \tag{6.24}\\
\text{Case (Sp):}&\quad w(x)=w(x;q^2,q^4). 
\endalign
$$
}
\Proof{
It is well known that Macdonald's symmetric polynomials $P_\mu(x)$ have 
the property
$$
\left[P_\mu(x)w(x)\right]_1=0 \quad \text{for}\ \ \mu\ne 0. 
\tag{6.25}
$$
(This is equivalent to saying that $P_\mu$ is orthogonal 
to 1 if $\mu\ne 0$.)
Since $P_\mu(x)$ form a $\C$-basis for $\C[x]^{\frak{S}_n}$, 
this property determines the functional 
$F \mapsto\left[F(x)w(x)\right]_1$ on $\C[x]^{\frak{S}_n}$ 
up to scalar multiples. 
Hence, the proof of Proposition 6.5 is reduced to show 
that the left hand side of (6.23) has the same property, when it 
is regarded as a functional on $\C[x]^{\frak{S}_n}$. 
It can be done by using the invariance of the functional $h_G$.  
In fact, property (6.1) implies that, 
for the central element $C_1$ of $\Uq$ in (5.12), we have 
$$
h_G(C_1.\varphi)=\vep(C_1) h_G(\varphi) \quad 
\text{for any}\ \ \varphi\in \Aq,
\tag{6.26}
$$
where $\vep(C_1)=\sum_{k=1}^N q^{2(N-k)}$. 
For each $\ld\in P^+_\k$, the
zonal spherical function $\varphi(\ld)$ is an eigenfunction 
of $C_1$ with eigenvalue $\chi_\ld(C_1)=\sum_{j=1}^N q^{2(\ld_k+N-k)}$. 
Hence we have 
$$
\chi_\ld(C_1)h_G(\varphi(\ld))=\vep(C_1) h_G(\varphi(\ld)). 
\tag{6.27}
$$
from (6.26).
Since $\chi_\ld(C_1)\ne\vep(C_1)$ unless $\ld=0$, we have 
$h_G(\varphi(\ld))= 0$ for any $\ld\ne 0$, as expected. 
}
In our setting, the scalar product $\langle\ ,\ \rangle'_{q,t}$ of 
Macdonald \cite{M2} takes the form
$$
\langle\ F,G\ \rangle' =\frac{1}{n!}\left[F(x)^* G(x)w(x)\right]_1. 
\tag{6.28}
$$
Here the $\ast$-operation on $\C[x]$ is given $x_k^*=x^{-1}_k$ 
for $1\le k\le n$. 
Note that, with this $\ast$-operation, $\C[x]$ is a $\ast$-subalgebra 
of $A(\T)$ endowed with the $\ast$-operation such that $z_j^*=z^{-1}_j$ 
for $1\le j\le N$. 
Recall that the subalgebra ${\Cal H}$ is closed under the $\ast$-operation 
of $\Aq$ representing the quantum unitary group (Proposition 6.2). 
Since the restriction mapping ${\Cal H}\to A(\T)$ is a $\ast$-homomorphism,
Proposition 6.5 implies that the hermitian form $\langle\ ,\ \rangle_G$ is 
related to Macdonald's scalar product through 
$$
\langle\varphi, \psi\rangle_G = 
\frac{\langle F, G\rangle'}{\langle 1, 1\rangle'}, 
\tag{6.29}
$$
if $\varphi,\psi\in{\Cal R}$ and $\varphi|_\T=F(x), \psi|_\T=G(x)$. 
\par
Proposition 6.1 and formula (6.29) implies that 
Macdonald's symmetric polynomials $P\mu$ are orthogonal under the 
scalar product $\langle\ ,\ \rangle'$. 
Furthermore, as to the square length of $P_\mu$, we have
$$
\frac{\langle P_\mu,P_\mu\rangle'}{\langle 1, 1\rangle'}
=\frac{c(\ld)^2}{d(\ld)}, 
\tag{6.30}
$$
by Proposition 6.3. 
%
%
%
\subsection{6.4. Computation of the ratio 
$\langle P_\mu,P_\mu\rangle'/\langle 1, 1\rangle'$ }
%
In the rest of this section, we will evaluate the ratio (6.30) of 
square lengths, by using a result of Macdonald on the value of 
the principal specialization of $P_\mu$. 
\par
We first recall a result of Macdonald \cite{M2} on the principal 
specialization of Macdonald's symmetric polynomials $P_\mu(x;q,t)$ : 
{\sl For each partition $\mu=(\mu_1,\cdots,\mu_n)$, one has}
$$
P_\mu(t^{n-1},t^{n-2},\cdots,1;q,t)= 
t^{\sum_{k=1}^n(k-1)\mu_k}\prod_{s\in\mu}
\frac{1-q^{a'(s)}t^{-\ell'(s)+n}}{1-q^{a(s)}t^{\ell(s)+1}}. 
\tag{6.31}
$$
In formula (6,30), the symbols 
$a(s),a'(s), \ell(s),\ell'(s)$ stand for the {\it arm-length}, 
{\it coarm-length}, {\it leg-length}, {\it coleg-length} 
of a box $s$ in the Young diagram $\mu$, respectively. 
If the box $s$ has the coordinates $(i,j)$ $(1\le i\le n, 1\le j\le\mu_i)$
in the Young diagram, 
they are given by
$$
a(s)=\mu_i-j,\ \ a'(s)=j-1,\ \ \ell(s)=\mu'_j -i, \ \ \ell'(s)=i-1,
\tag{6.32}
$$
where $\mu'$ denotes the conjugate partition of $\mu$. 
We also remark that (6.31) is a generalization of the following 
well-known formula for the Schur functions $s_\mu(x)=P_\mu(x;q,q)$:
$$
s_\mu(q^{n-1},q^{n-2},\cdots,1)= 
q^{\sum_{k=1}^n(k-1)\mu_k}\prod_{s\in\mu}
\frac{1-q^{c(s)}}{1-q^{h(s)}},
\tag{6.33}
$$
where $c(s)=a'(s)-\ell'(s)+n$ and $h(s)=a(s)+\ell(s)+1$ are the 
{\it content} and the {\it hook-length} of the box $s$. 
At present, the author does not know whether the principal 
specialization $(q^\rho,\varphi(\ld))$ of the zonal spherical function 
$\varphi(\ld)$ can be effectively evaluated as in (6.31), 
within the framework of quantum homogeneous spaces. 
\par
As to the value 
$\langle P_\mu(x;q,t),P_\mu(x;q,t)\rangle'_{q,t}$
of the scalar product, the following formula 
is proposed by Macdonald \cite{M3}:
{\sl For each partition $\mu=(\mu_1,\cdots,\mu_n)$, one has}
$$
\frac{\langle P_\mu(x;q,t),P_\mu(x;q,t)\rangle'_{q,t}}
{\langle1, 1\rangle'_{q,t}}
=\prod_{s\in\mu} 
\frac{(1-q^{a'(s)}t^{-\ell'(s)+n})(1-q^{a(s)+1}t^{\ell(s)})}
{(1-q^{a'(s)+1}t^{-\ell'(s)+n-1})(1-q^{a(s)}t^{\ell(s)+1})}. 
\tag{6.34}
$$
We will derive this formula 
for the two special cases corresponding to Cases (SO) and (Sp), 
by using our expression (6.30) and Macdonald's formula (6.31) 
for the principal specialization. 
\medpagebreak\par\noindent
Case (SO): 
In this case, the point $q^{\rho}=(q^{n-1},q^{n-2},\cdots,1)$ gives the 
value $x=(q^{2(n-1)},q^{2(n-2)},\cdots,1)$. 
Hence we have
$$
\align
c(\ld)&=P_\mu(q^{2(n-1)},q^{2(n-2)},\cdots,1;q^4,q^2) \tag{6.35}\\
&=q^{2\sum_{k=1}^n(k-1)\mu_k}\prod_{s\in\mu}
\frac{1-q^{2(2a'(s)-\ell'(s)+n)}}{1-q^{2(2a(s)+\ell(s)+1)}}, 
\endalign
$$
replacing $(q,t)$ in (6.31) by $(q^4,q^2)$.
On the other hand, we compute 
$$
\align
d(\ld)&=s_\ld(q^{2(n-1)},q^{2(n-2)},\cdots,1) \tag{6.36}\\
&=q^{2\sum_{k=1}^{n}(k-1)\ld_k}
\prod_{p\in\ld}
\frac{1-q^{2(a'(p)-\ell'(p)+n)}}
{1-q^{2(a(p)+\ell(p)+1)}}\\
&=q^{4\sum_{k=1}^{n}(k-1)\mu_k}
\prod_{s\in\mu}
\frac{(1-q^{2(2a'(s)-\ell'(s)+n)})(1-q^{2(2a'(s)-\ell'(s)+n+1)})}
{(1-q^{2(2a(s)+\ell(s)+1)})(1-q^{2(2a(s)+\ell(s)+2)})}.
\endalign
$$
The last equality follows from the simple fact that 
the Young diagram $\ld=2\mu$ 
is the duplication of $\mu$ in the horizontal direction. 
When we reparametrize the factors by the boxes in $\mu$, 
each horizontally adjacent pair of boxes $p$ in $\ld$, 
corresponding to a same $s$ in $\mu$, 
gives rise to the two factors with $a(p)$ replaced by 
$2a(s)$,$2a(s)+1$, and with $a'(p)$ replaced by 
$2a'(s)$,$2a'(s)+1$. 
Combining (6.35) and (6.36), we obtain 
$$
\frac{\langle P_\mu,P_\mu\rangle'}{\langle 1, 1\rangle'}
=\frac{c(\ld)^2}{d(\ld)}
=\prod_{s\in\mu}
\frac{(1-q^{2(2a'(s)-\ell'(s)+n)})(1-q^{2(2a(s)+\ell(s)+2)})}
{(1-q^{2(2a'(s)-\ell'(s)+n+1)})(1-q^{2(2a(s)+\ell(s)+1)})}.
\tag{6.37}
$$
This is exactly the formula obtained from (6.34) replacing 
$(q^4,q^2)$ for $(q,t)$.
\medpagebreak\par\noindent
Case (Sp): 
In this case, the point $q^{\rho}=(q^{2n-1},q^{2n-2},\cdots,1)$ corresponds 
to the value $x=q(q^{4(n-1))},q^{4(n-2)},\cdots,1)$. 
Replacing $(q,t)$ in (6.31) by $(q^2,q^4)$, we have 
$$
\align
c(\ld)
&=q^{\sum_{k=1}\mu_k}P_\mu(q^{4(n-1)},q^{4(n-2)},\cdots,1;q^2,q^4) \tag{6.38}\\
&=q^{\sum_{k=1}^n(4k-3)\mu_k}\prod_{s\in\mu}
\frac{1-q^{2(a'(s)-2\ell'(s)+2n)}}{1-q^{2(a(s)+2\ell(s)+2)}}. 
\endalign
$$
On the other hand, we compute 
$$
\align
d(\ld)&=s_\ld(q^{2(2n-1)},q^{2(2n-2)},\cdots,1) \tag{6.39}\\
&=q^{2\sum_{k=1}^{2n}(k-1)\ld_k}
\prod_{p\in\ld}\frac{1-q^{2(a'(p)-\ell'(p)+2n)}}
{1-q^{2(a(p)+\ell(p)+1)}}\\
&=q^{2\sum_{k=1}^{n}(4k-3)\mu_k}
\prod_{s\in\mu}
\frac{(1-q^{2(a'(s)-2\ell'(s)+2n)})(1-q^{2(a'(s)-2\ell'(s)+2n-1)})}
{(1-q^{2(a(s)+2\ell(s)+1)})(1-q^{2(a(s)+2\ell(s)+2)})}.
\endalign
$$
This time the Young diagram $\ld=(2\mu')'$  
is the duplication of $\mu$ in the vertical direction. 
Accordingly, $\ell(p)$ should be replaced by 
$2\ell(s)$,$2\ell(s)+1$, and $\ell'(p)$ by 
$2\ell'(s)$,$2\ell'(s)+1$. 
Combining (6.38) and (6.39), we obtain 
$$
\frac{\langle P_\mu,P_\mu\rangle'}{\langle 1, 1\rangle'}
=\frac{c(\ld)^2}{d(\ld)}
=\prod_{s\in\mu}
\frac{(1-q^{2(a'(s)-2\ell'(s)+2n)})(1-q^{2(a(s)+2\ell(s)+1)})}
{(1-q^{2(a'(s)-2\ell'(s)+2n-1)})(1-q^{2(a(s)+2\ell(s)+2)})}.
\tag{6.40}
$$
This coincides with the formula obtained from (6.34) replacing 
$(q^2,q^4)$ for $(q,t)$.
%
%

%

\Refs
\widestnumber\key{NUW2}
%
\ref\key{GK} \by A.M.\,Gavrilik and A.U.\,Klimyk
\paper $q$-Deformed orthogonal and pseudo-orthogonal algebras and 
their representations
\jour Lett. Math. Phys. \vol 21 \yr 1991 \pages 215--220
\endref
%
\ref\key{H1} \by T.\,Hayashi 
\paper Quantum deformation of classical groups
\jour Publ. RIMS \vol 28 \yr 1992 \pages 57--81
\endref
%
\ref\key{H2} \bysame
\paper Non-existence of homomorphisms between quantum groups
\jour Tokyo J. Math. \toappear
\endref
%
\ref\key{J} \by M.\,Jimbo
\paper A $q$-analogue of $U_q(\gl(N+1))$,
Hecke algebra and the Yang-Baxter equation,
\jour Lett. Math. Phys. \vol 11 \yr 1986 \pages 247--252
\endref
%
\ref\key{Koe} \by H.T.\,Koelink
\paper The addition formula for continuous $q$-Legendre polynomials
and associated spherical elements on the $\SU(2)$ quantum group 
related to Askey-Wilson polynomials
\jour SIAM J. Math. Anal. \toappear
\endref
%
%
\ref\key{K1}\by T.H.\,Koornwinder
\paper Continuous $q$-Legendre polynomials as spherical matrix elements 
of irreducible representations of the quantum $\SU(2)$ group
\jour CWI Quarterly \vol 2 \yr 1989 \pages 171--173
\endref
%
\ref\key{K2} \bysame 
\paper Orthogonal polynomials in connection with quantum groups
\inbook in ``Orthogonal Polynomials: Theory and Practice'' 
\ed P. Nevai \bookinfo NATO ASI Series
\publ Kluwer Academic Publishers \yr 1990 \pages 257--292
\endref
%
\ref\key{K3} \bysame 
\paper Askey-Wilson polynomials as zonal spherical functions on the 
$\SU(2)$ quantum group
\jour SIAM J. Math. Anal. \vol 24 \yr 1993 \pages 795--813
\endref
%
\ref\key{Ku} \by P.P.\,Kulish
\paper Quantum groups and quantum algebras as symmetries of 
dynamical systems
\paperinfo preprint YITP/K-959 \yr 1991
\endref
%
\ref\key{L} \by G.\,Lusztig
\paper Quantum deformations of certain simple modules over enveloping
algebras
\jour Advances in Math. \vol 70 \yr 1988 \pages 237--249
\endref
%
\ref\key{M1} \by I.G.\,Macdonald 
\book Symmetric functions and Hall polynomials,
\publ Oxford University Press \yr 1979
\endref
%
\ref\key{M2} \bysame 
\paper A new class of symmetric functions 
\inbook Actes $20^e$ S\'{e}minaire Lotharingien
\publ Publ. I. R. M. A. Strasbourg \yr 1988 \pages 131--171
\endref
%
\ref\key{M3} \bysame 
\paperinfo Draft of Chapter VI, for the new edition of the book 
``Symmetric functions and Hall polynomials''
\endref
%
\ref\key{N1} \by M.\,Noumi 
\paper Quantum groups and $q$-orthogonal polynomials --- Towards 
a realization of Askey-Wilson polynomials on $\SU_q(2)$
\inbook in ``Special Functions''
\bookinfo ICM-90 Satellite Conference Proceedings 
\eds M. Kashiwara and T. Miwa
\publ Springer-Verlag \yr 1991 \pages 260--288
\endref
%
\ref\key{N2} \bysame
\paper A remark on semisimple elements in $U_q(\frak{sl}(2);\Bbb{C})$
\inbook in ``Combinatorial Aspects in Representation Theory and Geometry"
\bookinfo RIMS Kokyuroku \vol 765 \yr 1991 \pages 71--78
\endref
%
\ref\key{N3} \bysame
\paper A realization of Macdonald's symmetric polynomials on 
quantum homogeneous spaces
\inbook in Proceedings of the XXI International Conference on 
Differential Geometric Methods in Theoretical Physics, 
Tianjin, China, 5--9 June 1992
\bookinfo Int. J. Mod. Phys. A(Proc. Suppl.) \vol 3A \yr 1993
\pages 218--223
\endref
%
\ref\key{NM1} \by M.\,Noumi and K.\,Mimachi
\paper Quantum 2-spheres and big $q$-Jacobi polynomials
\jour Commun. Math. Phys. \vol 128 \yr 1990 \pages 521--531
\endref
%
\ref\key{NM2} \bysame
\paper Spherical functions on a family of quantum 3-spheres
\jour Compositio Mathematica \vol 83 \yr 1992 \pages 19--42
\endref
%
\ref\key{NM3} \bysame
\paper Rogers's $q$-ultraspherical polynomials on a quantum 2-sphere
\jour Duke Math. J. \vol 63 \yr 1991 \pages 65--80
\endref
%
\ref\key{NM4} \bysame
\paper Askey-Wilson polynomials and the quantum group $\SU_q(2)$
\jour Proc. Japan Acad.  \vol 66 \yr 1990 \pages 146--149
\endref
%
\ref\key{NM5} \bysame
\paper Askey-Wilson polynomials as spherical functions on $\SU_q(2)$
\inbook in ``Quantum Groups''
\ed P.P. Kulish
\bookinfo Proceedings of Workshops held in the Euler International 
Mathematical Institute, Leningrad, Fall 1990, Lecture Notes in Math.
\vol 1510 \publ Springer Verlag \yr 1992 \pages 98--103
\endref
%
\ref\key{NUW1} \by M.\,Noumi, T.\,Umeda and M.\,Wakayama
\paper A quantum analogue of the Capelli identity and
an elementary differential calculus on $\GL_q(n)$
\paperinfo preprint J-Tokyo-Math 91-16 \yr 1991
\endref
%
\ref\key{NUW2} \bysame
\paper A quantum dual pair $(\frak{sl}_2, \frak{o}_n)$ 
and the associated Capelli identity
\jour Lett. Math. Phys. \toappear
\endref
%
\ref\key{NYM} \by M.\,Noumi, H.\,Yamada and K.\,Mimachi
\paper Finite dimensional representations of the quantum group  $\GL_q(n;\Bbb{C})$ 
and the zonal spherical functions on
${\text{\rm U}}_q(n-1)\backslash {\text{\rm U}}_q(n)$
\jour Japan. J. Math. \vol 19 \yr 1993 \pages 31--80
\endref
%
\ref\key{O} \by G.I.\,Olshanski
\paper Twisted Yangians and infinite dimensional classical Lie algebras
\inbook in ``Quantum Groups''
\ed P.P. Kulish
\bookinfo Proceedings of Workshops held in the Euler International 
Mathematical Institute, Leningrad, Fall 1990, Lecture Notes in Math.
\vol 1510 \publ Springer-Verlag \yr 1992 \pages 104--120
\endref
%
\ref\key{RTF} \by N.Yu.\,Reshetikhin, L.A.\,Takhtajan and L.D.\,Faddeev 
\paper Quantization of Lie groups and Lie algebras 
\jour Algebra and Analysis \vol 1 \yr 1989 \pages 178--206 
\transl\nofrills English transl. in 
\jour Leningrad Math. J. \vol 1 \yr 1990 \pages 193--225
\endref
%
\ref\key{R} \by M.\,Rosso
\paper Finite dimensional representations of the quantum analog 
of the enveloping algebra of a complex simple Lie algebra
\jour Commun. Math. Phys. \vol 117 \yr 1988 \pages 581-593
\endref
%
\ref\key{TT} \by E.\,Taft and J.\,Towber
\paper Quantum deformation of flag schemes and Grassmann schemes. 
I. A $q$-deformation of the shape-algebra for $\text{\rm GL}(n)$
\jour Journal of Algebra \vol 142 \yr 1991 \pages 1--36
\endref
%
\ref\key{T} \by T.\,Tanisaki
\paper Killing forms, Harish-Chandra isomorphisms, and universal
$R$-matrices for quantum algebras
\jour International Journal of Modern Physics A 
Vol.7 Suppl.1B \yr 1992 \pages 941--961
\paperinfo in the Proceedings of the RIMS Research Project 1991 ``Infinite Analysis''
\endref
%
\ref\key{UT} \by K.\,Ueno and T.\,Takebayashi
\paper Zonal spherical functions on quantum symmetric spaces and 
Macdonald's symmetric polynomials
\inbook in ``Quantum Groups''
\ed P.P. Kulish
\bookinfo Proceedings of Workshops held in the Euler International 
Mathematical Institute, Leningrad, Fall 1990, Lecture Notes in Math.
\vol 1510 \publ Springer-Verlag \yr 1992 \pages 142--147
\endref
%
\ref\key{W} \by S.L.\,Woronowicz
\paper Compact matrix pseudogroups
\jour Commun. Math. Phys. \vol 111 \yr 1987 \pages 613--665
\endref
%
\endRefs
\enddocument